\documentclass[11pt]{article}

\usepackage{amsmath,amssymb,amsfonts,mathtools}
  \usepackage{paralist}
  \usepackage{graphics} 
  \usepackage{epsfig} 
\usepackage{graphicx}  \usepackage{epstopdf}
 \usepackage[colorlinks=true]{hyperref}
 \usepackage[superscript, biblabel, nomove]{cite}

\usepackage[percent]{overpic}
\usepackage{tikz}   
\usepackage{enumitem}
\usepackage{pict2e}
\usepackage[position=bottom]{subfig}   
\usepackage{lipsum}
\usepackage{wrapfig}
\usepackage{graphicx,overpic}
\usepackage{siunitx}
\usepackage{textcomp}
\usepackage{titlesec}
\usepackage{bbm}

\usepackage{ntheorem} 
\theorembodyfont{\upshape}
\newtheorem{remark}{Remark}[section]
\numberwithin{equation}{section}

\usepackage{graphicx}
\graphicspath{{img/}} 

\usepackage{multirow}
\usepackage{hhline}

\def\b{{\boldsymbol{b}}}

\def\f{{\boldsymbol{f}}}
\def\g{{\boldsymbol{g}}}

\def\u{{\boldsymbol{u}}}
\def\v{{\boldsymbol{v}}}
\def\x{{\boldsymbol{x}}}

\def\J{{\boldsymbol{J}}}

\def\F{{\boldsymbol{F}}}

\def\blam{{\boldsymbol{\lambda}}}
\def\bLam{{\boldsymbol{\Lambda}}}

\def\cl {\nonumber \\}
\def\el {\nonumber}

\DeclareMathOperator*{\argmin}{arg\,min}

\begin{document}

\title{Data driven learning to enhance a kinetic model of distressed crowd dynamics}

\maketitle

{\centering
\author{Daewa Kim \\
\footnotesize{
Department of Mathematical Sciences, University of Delaware,
501 Ewing Hall
Newark, DE 19716, USA}\\
\footnotesize {daewakim@udel.edu}
\vspace{0.5cm}}\\

\author{Demetrio Labate, Kamrun Mily, and Annalisa Quaini\\
\footnotesize{
Department of Mathematics, University of Houston, 3551 Cullen Blvd, Houston TX 77204, USA}\\
\footnotesize{ dlabate2@central.uh.edu, kmily@cougarnet.uh.edu, aquaini@central.uh.edu}
}\\
}


\begin{abstract}
The mathematical modeling of crowds is complicated 
by the fact that crowds possess the behavioral ability to develop and adapt
moving strategies in response to the context.
For example, in emergency situations, people tend to alter their walking strategy in response to fear. To be able to simulate these situations,
we consider a kinetic model of crowd
dynamics that features the
level of stress as a parameter and propose to
estimate this key parameter
by solving an inverse crowd dynamics problem.
This paper states the mathematical problem and presents a method for its numerical solution. 
We show some preliminary results based on a synthetic data set, i.e.,
test cases where the exact stress level is known and the crowd density data are generated numerically by solving a forward crowd dynamics problem.\\
\smallskip

\noindent {$Keywords$: Crowd dynamics; Boltzmann-type kinetic model; Complex systems; Inverse problems; Parameter estimation.}\\
\smallskip

\end{abstract}

\section{Introduction}

\begin{figure}[htb!]
\begin{center}
    \begin{overpic}[percent,width=0.32\textwidth]{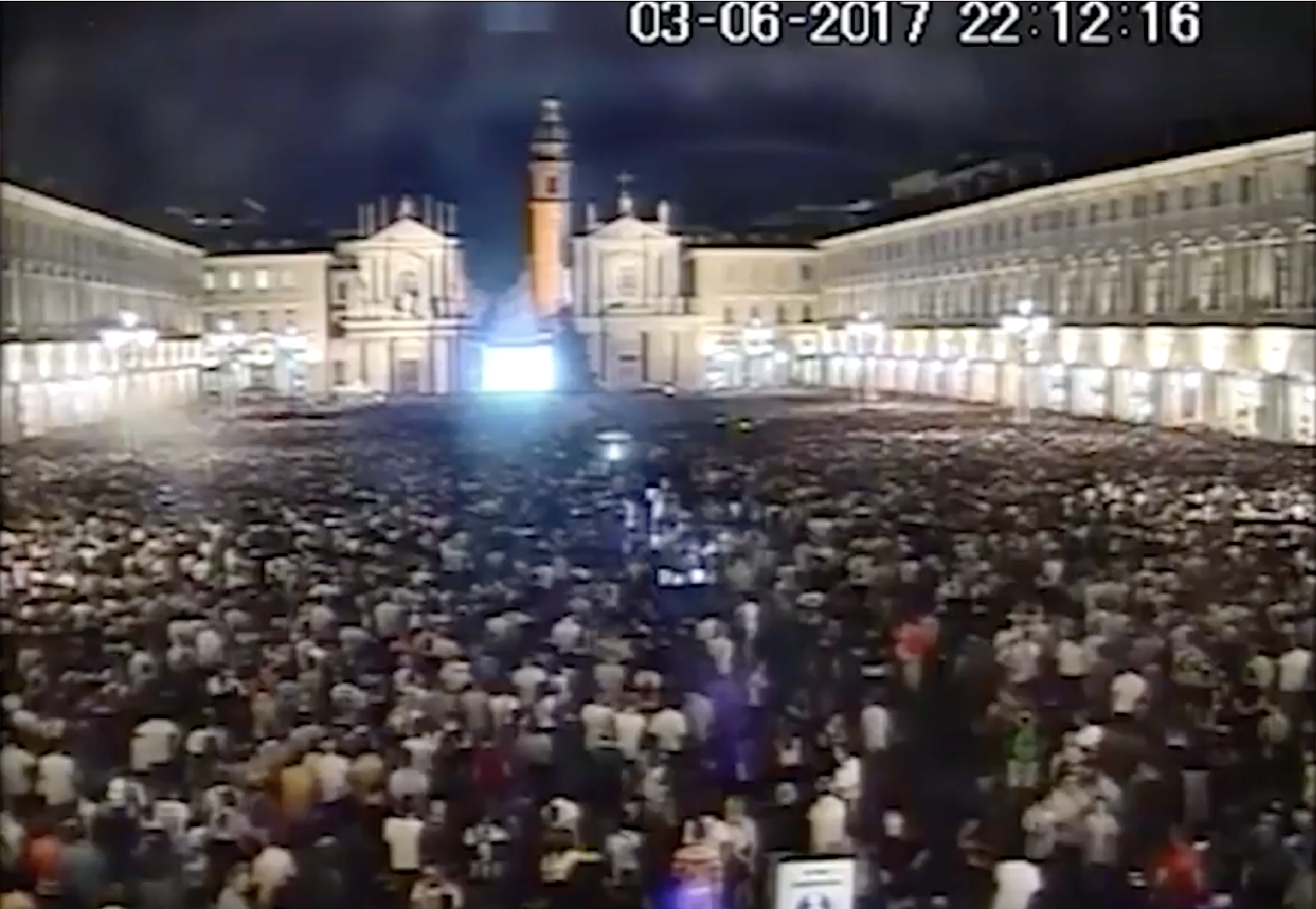}
    \end{overpic}
    \begin{overpic}[percent,width=0.32\textwidth]{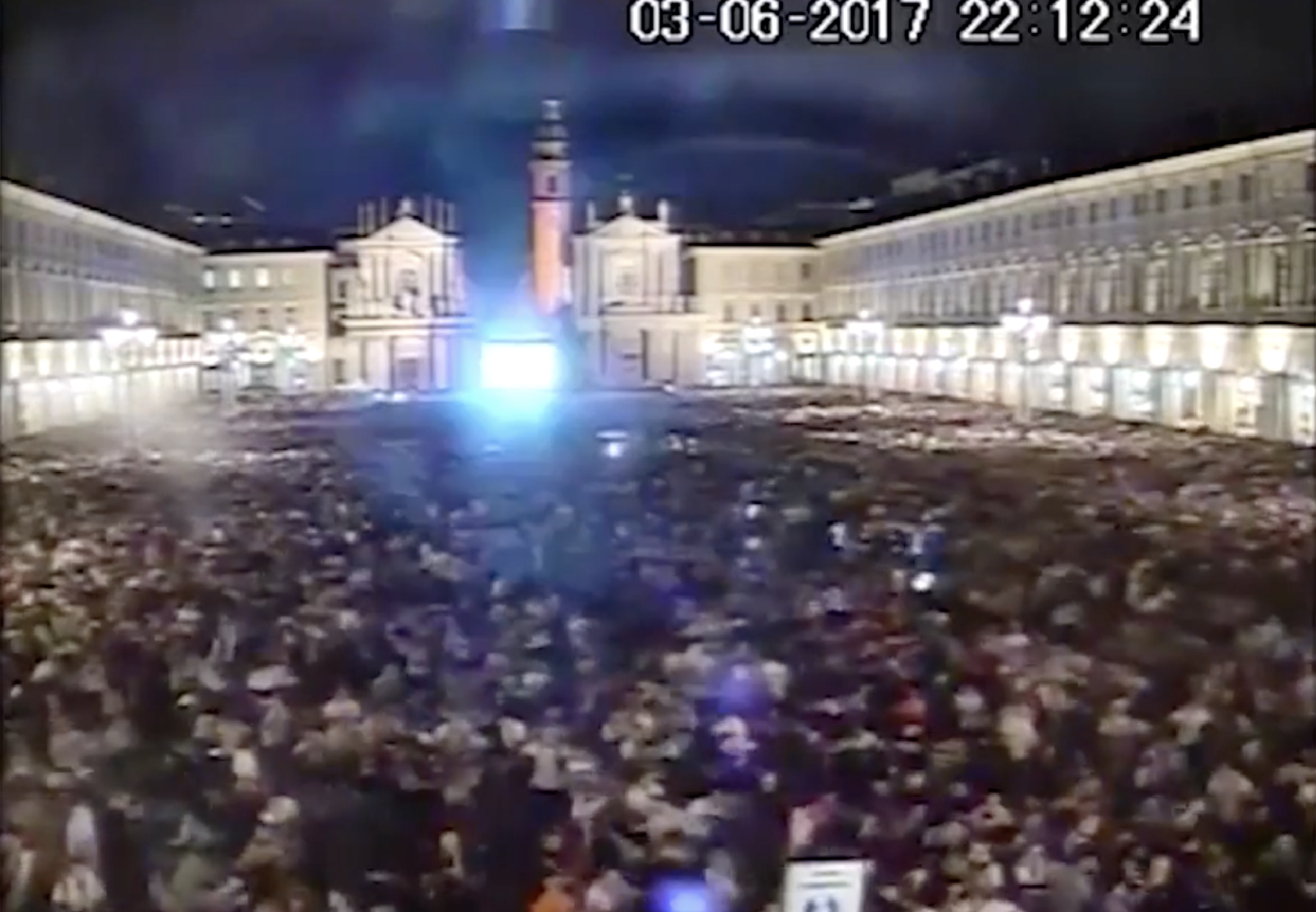}
    \end{overpic}
    \begin{overpic}[percent,width=0.32\textwidth]{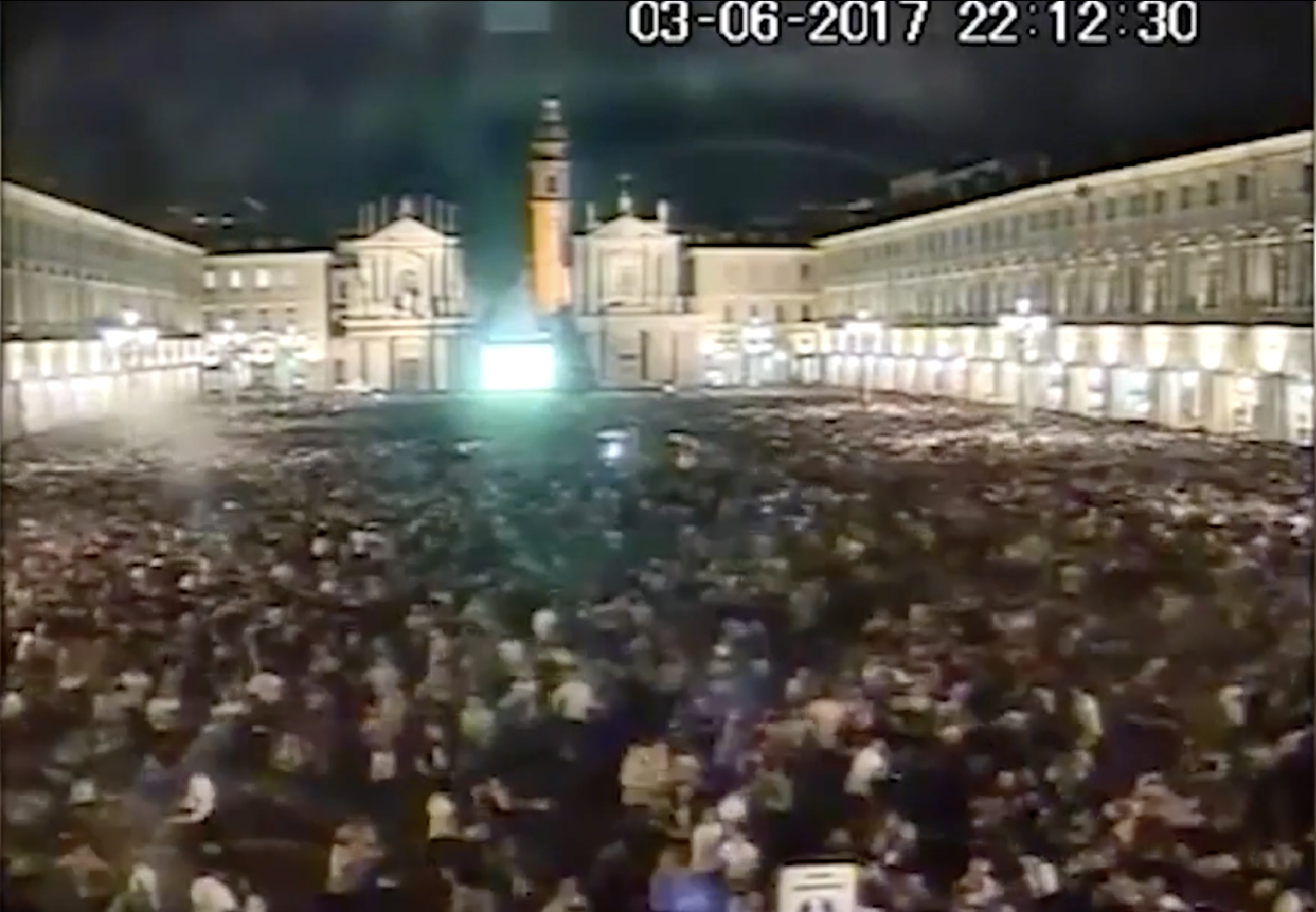}
    \end{overpic}
    \caption{Security camera frames of a surge in the crowd gathered in Torino (Italy) to watch a soccer game on June 3rd, 2017. The stampede 
    led to 3 people dead and over 1600 injured. Video available at \url{https://youtu.be/yuqcNgcgzIA?feature=shared}.}
    \label{fig:Torino}
\end{center}
\end{figure}

Stampede-related injuries and deaths continue to happen at 
events like concerts and religious, festive, or sport 
gatherings. See Fig.~\ref{fig:Torino} for an example captured on security
cameras. Since experiments on panicked people are 
unethical and potentially harmful, the most practical 
and effective way to improve crowd management is via modeling and computational simulations.

Unlike inert matter, crowds possess the behavioral ability to develop and adapt moving strategies in response to the context, which significantly complicates their mathematical modeling. See Refs.~\citenum{Bellomo2022,Bellomo2023} for recent reviews on the topic.
A consequence of this complication is that, 
although the mathematical modeling of crowd motion has fascinated researchers
since the 1950s, the large majority of existing mathematical models assume that pedestrians do not alter their, typically
rational, walking behavior. However, in emergency situations,
people may behave irrationally in response to fear.
Thus, an accurate model of crowd dynamics in stressful situations cannot be fully captured by classical mechanics and must include heterogeneous behavior of individual entities and their interactions, accounting
for the effects of the individuals’ levels of stress 
and its propagation.~\cite{PhysRevE.75.046109,Helbing2009,BELLOMO20161,WIJERMANS2016142,HAGHANI201724}
While the widely used agent-based models~\cite{Helbing2000,WANG2015396,Zou2016,PhysRevE.98.032312,8667341} successfully account for heterogeneous individual behavior (each agent can develop a specific walking strategy), they fail to reproduce non-local and nonlinearly additive interactions, becoming progressively less accurate as the number of individuals increases. Thus, in this paper we focus on an alternative model that can reproduce such complex interactions.

Taking inspiration from the kinetic theory of gases,
kinetic models of crowd dynamics derive a Boltzmann-type evolution equation for the statistical
distribution function of the position and velocity 
of the pedestrians. In these models, people
are seen as \emph{active} particles, as opposed
to classical particles in gases. This leads
to a key difference: the interactions in 
kinetic models for gas dynamics are
conservative and reversible, while the interactions in the kinetic models for crowd dynamics are irreversible, 
non-conservative and, in some cases, nonlocal and nonlinearly additive. Hence, kinetic models can overcome the limitations of agents-based models in terms of complexity
of interactions.
An important consequence of the aforementioned difference 
is that often for active particles the
Maxwellian equilibrium does not exist.~\cite{Aristov_2019}
The reader interested in introductory concepts on the kinetic theory
of active particles is referred to Refs.~\citenum{Bellomo2013_new} and \citenum{BellomoBellouquid2015}. More recent developments can be found in 
Refs.~\citenum{Bellomo2017_book,Bellomo2015_new,Bellomo2019_new} and \citenum{kim_quaini}. Refs.~\citenum{Bellomo2022,Bellomo2023} focus on
aspects specific to the modeling of crowds taking into account behavioral dynamics.
A recent interesting spin on kinetic models as possible 
mathematical models for artificial intelligence, i.e., 
collective behavior of artificial (instead of natural) self-organizing systems, is discussed in
Ref.~\citenum{BELLOMO20241}.

The kinetic model considered in this paper was initially proposed in~\cite{Bellomo2013_new} and further developed in \cite{Bellomo2015_new,Agnelli2015,kim_quaini}. It assumes 
discrete walking directions, a continuous deterministic
walking speed, and it accounts for the presence of walls
and other kinds of obstacles in the walking domain. 
The kinetic models in \cite{Bellomo2013_new,Bellomo2015_new,Agnelli2015,kim_quaini}
features the level of stress as a parameter, that in the
numerical experiment is set to be constant in space
and time. This homogeneity of the stress level in a crowd is obviously not realistic. Hence, some models have introduced
the emotional state as a variable that, in response to interactions with other people, can change in space
and time and, in turn, alter the walking strategy~\cite{Bellomo2019_new,Bertozzi2014ContagionSI}. 
In fact, it is known that emotional contagion significantly 
affects the overall crowd dynamics in complicated real-life situations such as emergency evacuations~\cite{HAGHANI201724}.
A multiscale (from agent-based to continuum)
approach to crowd dynamics with emotional contagion is 
introduced in~\cite{Bertozzi2015}. Therein, 
the stress level is propagated by a Bhatnagar–Gross–Krook-like model
and results are limited to one space dimension.
Extension to 2D is presented in \cite{kim_quaini2020b}.
While treating the emotional state as a variable
leads to more realistic simulation results than representing
it as a constant parameter~\cite{Bellomo2019_new}, the computational cost of an additional variable can become significant as the domain size increases.

Since the overall goal
of crowd dynamics simulations is to contribute to crowd management in emergency situations~\cite{Banda2020},
the computational cost of realistic simulations must be contained. Thus, in this paper we propose a third alternative to treat the emotional state. The computational cost
is contained by considering the emotional state, i.e., 
the stress level, as a parameter that is optimized in 
space and time to match density data extracted from video frames
for realism. With a large database of videos such as the one in 
Fig.~\ref{fig:Torino}, one could split 
the database into validation and training
datasets and use the videos in the training
set to learn the evolution of the stress level via an optimization procedure that minimizes the difference between the density data extracted from the video and the simulated crowd density. In this paper, we present a possible optimization approach for this purpose. Parameter learning techniques have been proposed for agent-based models.~\cite{Lee2007,BERA201668,YAO2019314} 
However, to the best of our knowledge, it is the first time that such a technique is proposed for a kinetic model. 

\begin{figure}[htb!]
\centering
    \begin{overpic}[percent,width=0.305\textwidth, grid=false]{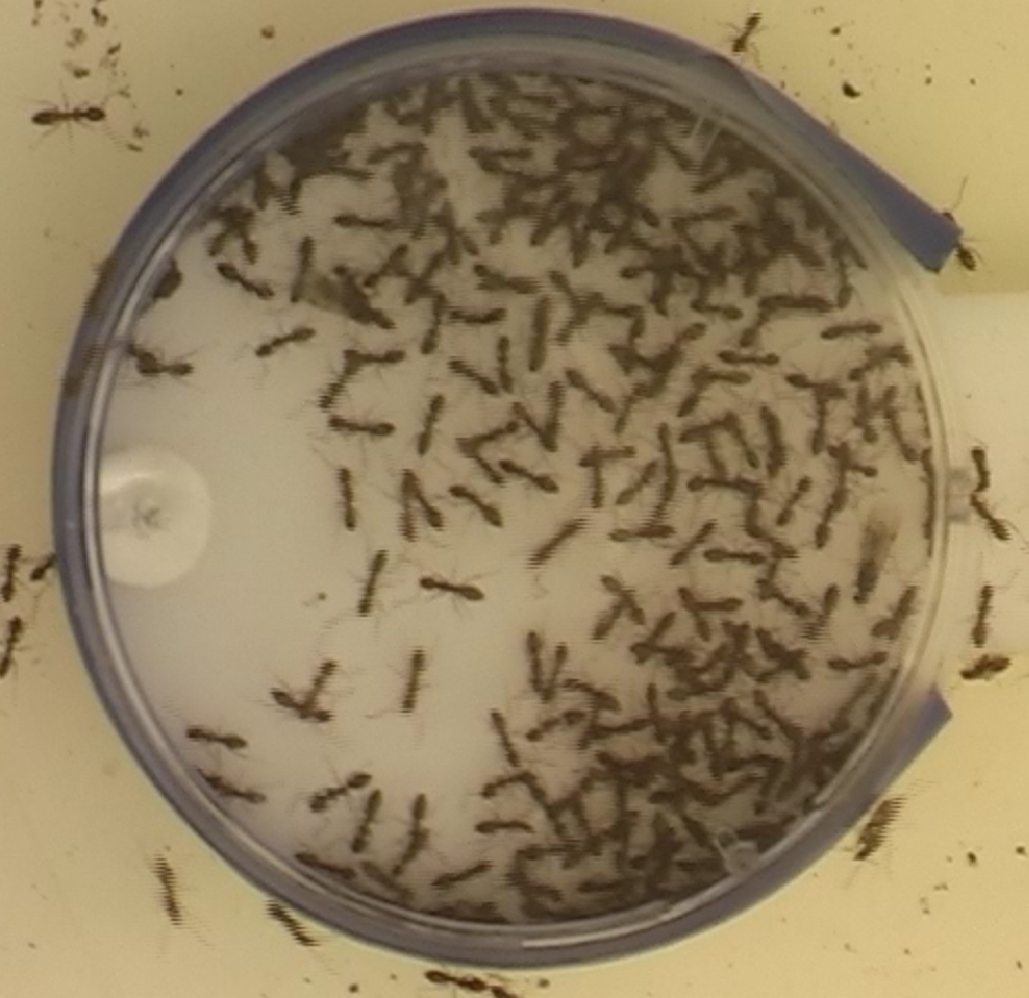}
    \put(110,48){exit}
    \put(107,50){\textcolor{red}{\vector(-1,0){13}}}
    \end{overpic}
    \quad \quad \quad \quad \quad 
    \begin{overpic}[percent,width=0.305\textwidth, grid=false]{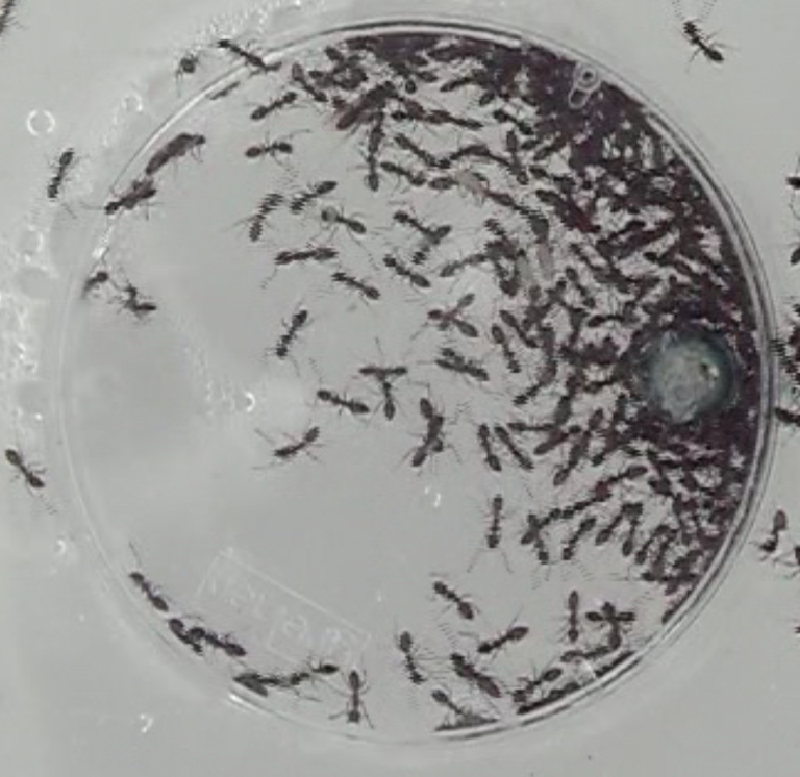}
    \put(110,60){exit}
    \put(107,60){\textcolor{red}{\vector(-1,-1){10}}}
    \put(110,37){column}
    \put(107,40){\textcolor{red}{\vector(-2,1){20}}}
    \end{overpic}
        \caption{Snapshots of experiments with panicking ants from Ref.~\citenum{SHIWAKOTI20111433}: circular chamber without (left) and with (right) column obstructing the exit.}
    \label{fig:ants}
\end{figure}
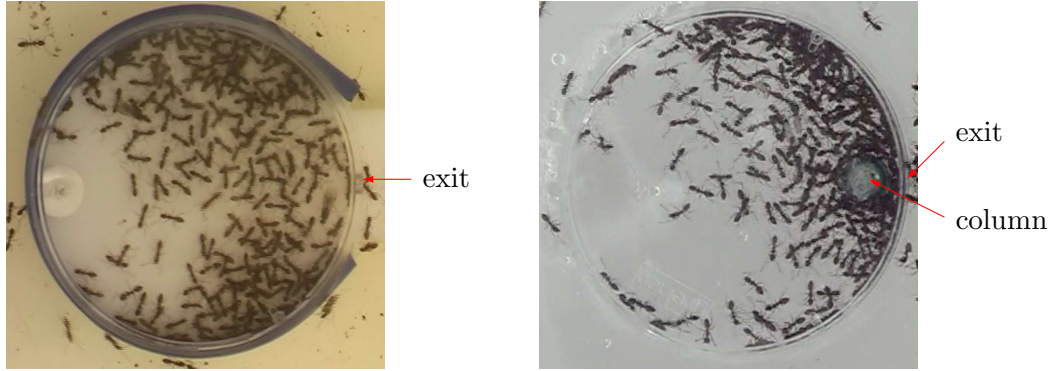

Extracting crowd density data from security camera footage is
challenging due to poor resolution and, typically, 
reduced visibility. Since these are challenges
related to the extraction of data from images
and not to the learning-from-data technique, 
in this paper we considered a simplified setting:
videos of experiments where a few hundred ants are enclosed in a simple chamber and an insect repellent is injected at the opposite end of the exit to induce panic. \cite{SHIWAKOTI20111433}
Despite the obvious differences between ants and human crowds, these experiments exhibit complex interactions which are comparable to those observed in a human crowd responding to stress or panic. \cite{Lee2007}
In the recorded experiments, the geometric setting is simple: a circular or square chamber with one exit and, in some cases,
a column obstructing the exit. Fig.~\ref{fig:ants} shows the circular chamber with and without obstructing column. It is clear that the image quality of the 
frames in Fig.~\ref{fig:ants} is superior when compared to  
Fig.~\ref{fig:Torino}.

There exist several methodologies to estimate the crowd density from images or video frames, including classical methods relying on texture analysis \cite{marana1998} and a multiplicity of more recent strategies based on deep learning architectures \cite{fu201,chen2022counting,khan2023revisiting}.
In this work, we assumes that such estimated crowd density, denoted by $\rho_v$,  
is available and we minimize the difference between $\rho_v$
and the crowd density computed by our kinetic model, 
which depends on the stress level parameter. This
minimization problem yields an estimation of the stress level, 
so it can be considered an inverse problem. 
In fact, inverse problems can be used for parameter
identification of
mathematical models, also called forward problems, whose state variables are observed through available data. 
In our case, the forward problem is our kinetic model for crowd 
dynamics, the state variables are the probability 
distribution functions, from which crowd density 
can be computed, and the observed data are the crowd densities
extracted from video images. In this paper, we solve
the inverse problem as a constrained optimization
problem, where the kinetic model is the constraint and the functional to be minimized is a measure of the mismatch between the density data and the computed density solution.

In our approach to the constrained minimization problem, 
the kinetic model is first discretized in time and then the
optimization step is performed. We present a
numerical method for its solution. In this first
work, we consider the simple geometrical settings 
used in the experiments on ants \cite{SHIWAKOTI20111433} 
and synthetic data, i.e., obtained by a forward crowd dynamics
problem. The comparison with real data extracted
from the videos of the ants experiments, generously
shared with us by the authors of Ref.~\citenum{SHIWAKOTI20111433}, 
will be presented in a follow-up paper. 

The rest of the paper is structured as follows. Sec.~\ref{sec:forward}
presents the time continuous inverse crowd dynamics problem
and Sec.~\ref{sec:td_model} describes its time discrete counterpart. 
Numerical methods for the solution are presented in Sec.~\ref{sec:num_sol}, while numerical results
are presented and discussed in Sec.~\ref{sec:numericalresults}. Sec.~\ref{sec:concl} draws the conclusions and 
mentions possible future perspectives.

\section{The time continuous inverse crowd dynamics problem}\label{sec:cont_model}

\subsection{The forward problem}\label{sec:forward}

Let  $\Omega \subset \mathbb{R}^2$ denote a bounded walkable domain. 
We assume that the boundary $\partial \Omega$ includes an exit $E$,
which could be the finite union of disjoint sets, and walls $W$. Here, $\overline{E}\cup \overline{W}
=\overline{\partial \Omega}$ and $E \cap W =\emptyset$. 
Let $\x=(x,y) \in \Omega$ denote position and ${\v}=v (\cos \theta, \sin \theta) \in \Omega_{\v}$ denote  
velocity, where $v$ is the velocity modulus, $\theta$ is the velocity direction, and
$\Omega_{\v} \subset \mathbb{R}^2$ is the velocity domain. 
For a system composed by a large number of people 
distributed inside  
$\Omega$, the distribution function is given by 
\[ f= f(t, \x, \v)\quad \text{for all} \,\,\, t \ge 0,  \,\, \x \in \Omega, \,\, \v \in \Omega_{\v}. \]
Since we use polar coordinates for the velocity, we can write the distribution function as $f= f(t, \x, v, \theta)$, where
$v$ is the walking speed and $\theta$ is the walking direction. 
Under suitable integrability conditions, $f(t, \x, v, \theta)d \x d \v$ represents the number of individuals who, at time $t$, 
are located in the infinitesimal rectangle $[x, x+dx] \times [y, y+dy]$ and have a velocity belonging to $[v, v + dv] 
\times [\theta, \theta+d\theta]$. 

For simplicity, 
we assume $\theta$ can take values in the set:
\[ I_{\theta}=\left \{ \theta_{i}= \frac{i-1}{N_d} 2\pi : i = 1, \dots, N_d \right \}, \]
where $N_d$ is the maxim number of possible directions. 
The scale of observation for a kinetic model 
can be considered between the scale of agent-based models (i.e., 
microscopic) and the scale of continuum models (i.e., macroscopic). Hence, the assumption of discrete
walking directions is reasonable because if there
was a large enough crowd in every possible walking direction
then a macroscopic model would be more suited. 
To further simplify the problem, we model the walking speed
$v$ as a continuous deterministic variable 
that evolves in time and space depending on the crowd density. 
This is also a reasonable assumption since the walking speed of a person typically depends on on the level of congestion around
them.

Since $v$ is now a deterministic variable, we can reduce
the distribution function to
\begin{equation}\label{eq:f}
f(t, \x, \theta)= \sum_{i=1}^{N_d} f^i(t, \x)\delta(\theta - \theta_i),
\end{equation}
where $f^i(t, \x)=f(t, \x, \theta_i)$ is related to the 
number of individuals that, at time $t$ and position $\x$, walk with direction $\theta_i$. In equation~\eqref{eq:f}, $\delta$ denotes the Dirac delta function.


From now on, we will assume that all variables are dimensionless.
See, e.g., Ref.~\citenum{kim_quaini} for details on how the non-dimensionalization process is carried
out. 
Due to the normalization of $f$ and each $f^i$, the dimensionless local density is obtained by summing the distribution functions over the set of directions:
\begin{align}\label{eq:rho}
\rho(t, \x)=\sum_{i=1}^{N_d}f^i(t, \x) .
\end{align}

\begin{figure}[h!]
\centering
\begin{overpic}[width=0.47\textwidth,grid=false]{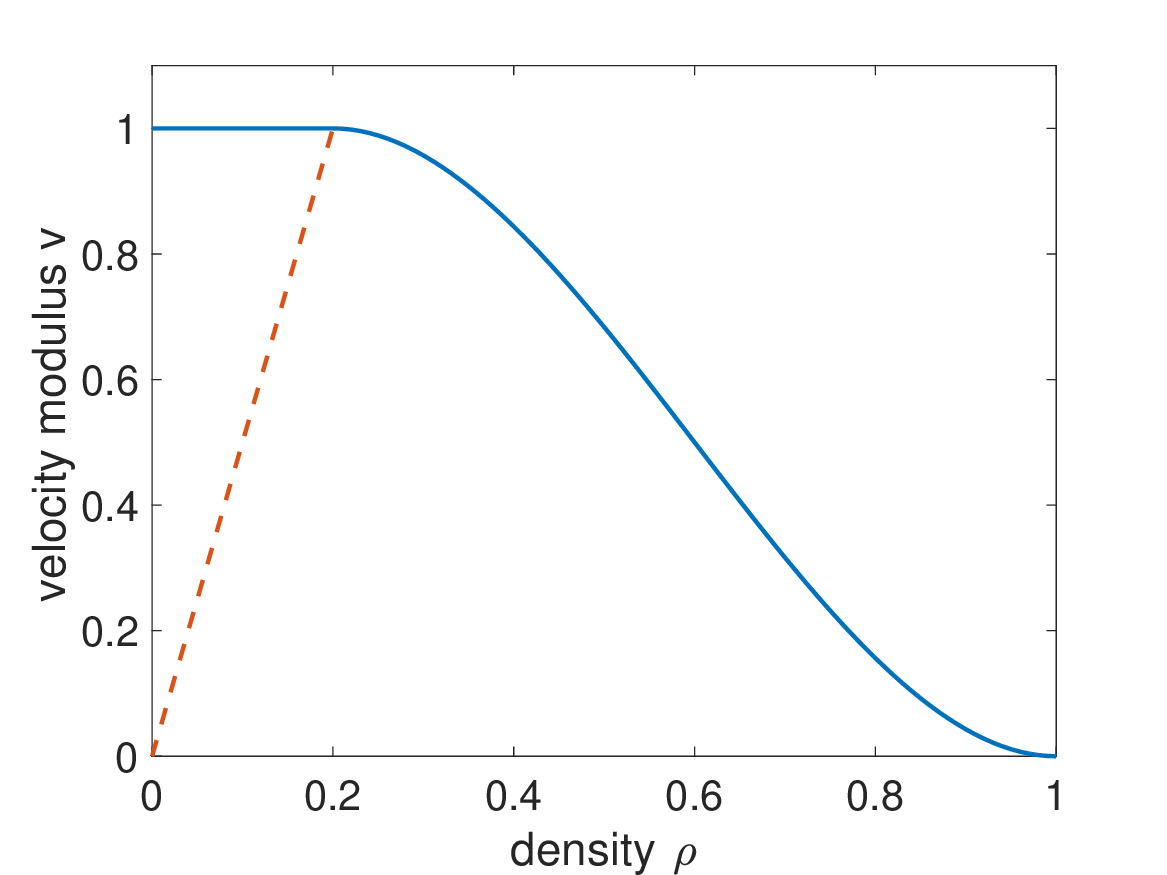}
\end{overpic}

\caption{
Dependence of the dimensionless velocity modulus $v$ on the dimensionless density $\rho$. The dashed red line connects the origin
of the axis to the critical value $\rho_c = 1/5$.
}\label{velocity}
\end{figure}

As mentioned above, the walking speed is a function of $\rho$
in \eqref{eq:rho}. 
For simplicity, we assume that the quality of the walkable domain
is high everywhere, meaning that a person can walk at the the maximum dimensionless speed
(i.e., 1) everywhere in $\Omega$.
This maximum walking speed
can be kept so long as the density around a person does not
reach a critical value $\rho_c$, past which the person
is forced to slow down. 
Following experimental estimates published in the literature, \cite{Schadschneider2011545} we set $\rho_c = 1/5$. In the slowdown zone, i.e., for
$\rho > \rho_c$, we assume the walking speed is a cubic polynomial
in $\rho$. See, e.g., Ref.~\citenum{kim_quaini} for the exact definition.
Figure~\ref{velocity} reports $v$ as a function of $\rho$.

\begin{remark}
The assumption of high quality of the walkable domain is reasonable 
for the experiments on ants, but it would make perhaps less sense in a real
life context such as the one in Fig.~\ref{fig:Torino}. An easy way to incorporate
the quality of the domain in the model is through a parameter. The reader interested
in this is referred to, e.g., Refs.~\citenum{Bellomo2022,Bellomo2023}.
\end{remark}

A model governing the dynamics of a 
crowd can be obtained by a balance of people, or
active particles, in an elementary volume
of the space of microscopic states, considering the net flow into such a volume due to
transport and interactions with other people and the environment:
\begin{align}
\frac{\partial f^i}{\partial t} &+ \nabla \cdot \left( \v [\rho] (t, \x) f^i(t, \x) \right) = \mathcal{I}^i[f](t, \x) \label{eq:model}
\end{align}
for $i= 1,2, \dots, N_d$. By $[\cdot]$, we mean nonlinear dependence on the argument. 
Functional $\mathcal{I}^i[f]$ represents the net balance of people 
that walk with direction $\theta_i$. One way to look at $\mathcal{I}^i[f]$ could be as follows:
\[
\mathcal{I}^i[f] = \mathcal{G}^i[f] - \mathcal{L}^i[f],
\]
where $\mathcal{G}$ and $\mathcal{L}$ represent gain and loss, both acting nonlinearly on $f$, of people walking with direction $\theta^i$ in the elementary volume. However, 
for our purpose, it is more convenient to write $\mathcal{I}^i[f]$
as: 
\begin{align}
\mathcal{I}^i[f](t, \x) = \mathcal{I}^i_G[f](t, \x) + \mathcal{I}^i_P[f](t, \x), \el
\end{align}
where $\mathcal{I}^i_G$ is an interaction between \emph{candidate} people and the environment, and 
 $\mathcal{I}^i_P$ is an interaction between \emph{candidate} and \emph{field} people.
This terminology is borrowed from gas dynamics: a \emph{test} individual $i$ is representative
of the whole system, and a \emph{candidate} individual $h$ reaches in probability 
the state of the test individual after interactions with the environment or with 
\emph{field} people.

We define $\mathcal{I}^i_G$ and $\mathcal{I}^i_P$ as
\begin{align}
\mathcal{I}^i_G[f](t, \x)& = (1 - \rho) \left( \sum_{h = 1}^{N_d} \mathcal{A}_h^i f^h(t, \x) - f^i(t,\x) \right), \cl
\mathcal{I}^i_P[f](t, \x) &= \rho \left( \sum_{h,k = 1}^{N_d} \mathcal{B}_{hk}^i [\rho] f^h(t, \x)f^k(t, \x) - f^i(t,\x) \rho\right), \el
\end{align}
 where
 \begin{itemize}
 \item[-] $\mathcal{A}_h^i$ is the probability that a candidate individual $h$, 
 i.e.~with direction  $\theta_h$, adjusts their direction into that 
 of the test individual $\theta_i$ due to the presence of walls, 
 obstacles and/or an exit. The following constraint for $\mathcal{A}_h^i$ 
 has to be satisfied:
\[
\sum_{i=1}^{N_d} \mathcal{A}_h^i=1 \quad \text{for all} \,\, h \in \{1, \dots, N_d\}.
\]
 \item[-] $\mathcal{B}_{hk}^i$ is the probability that a candidate individual
$h$ modifies their direction $\theta_h$ into that of the test individual $i$, 
i.e., $\theta_i$, due to the interaction with a field individual $k$ 
that moves with direction $\theta_k$. 
The following constraint for $\mathcal{B}_{hk}^i$ has to be satisfied:
\begin{equation}\label{eq:const_B}
\sum_{i=1}^{N_d} \mathcal{B}_{hk}^i=1 \quad \text{for all} \,\, h, k  \in \{1, \dots, N_d\}.
\end{equation}
 \end{itemize}

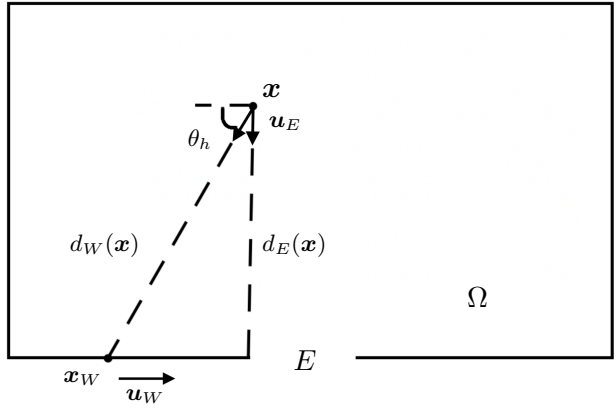
\begin{figure}[h!]
\centering

\begin{overpic}[width=0.5\textwidth,grid=false]{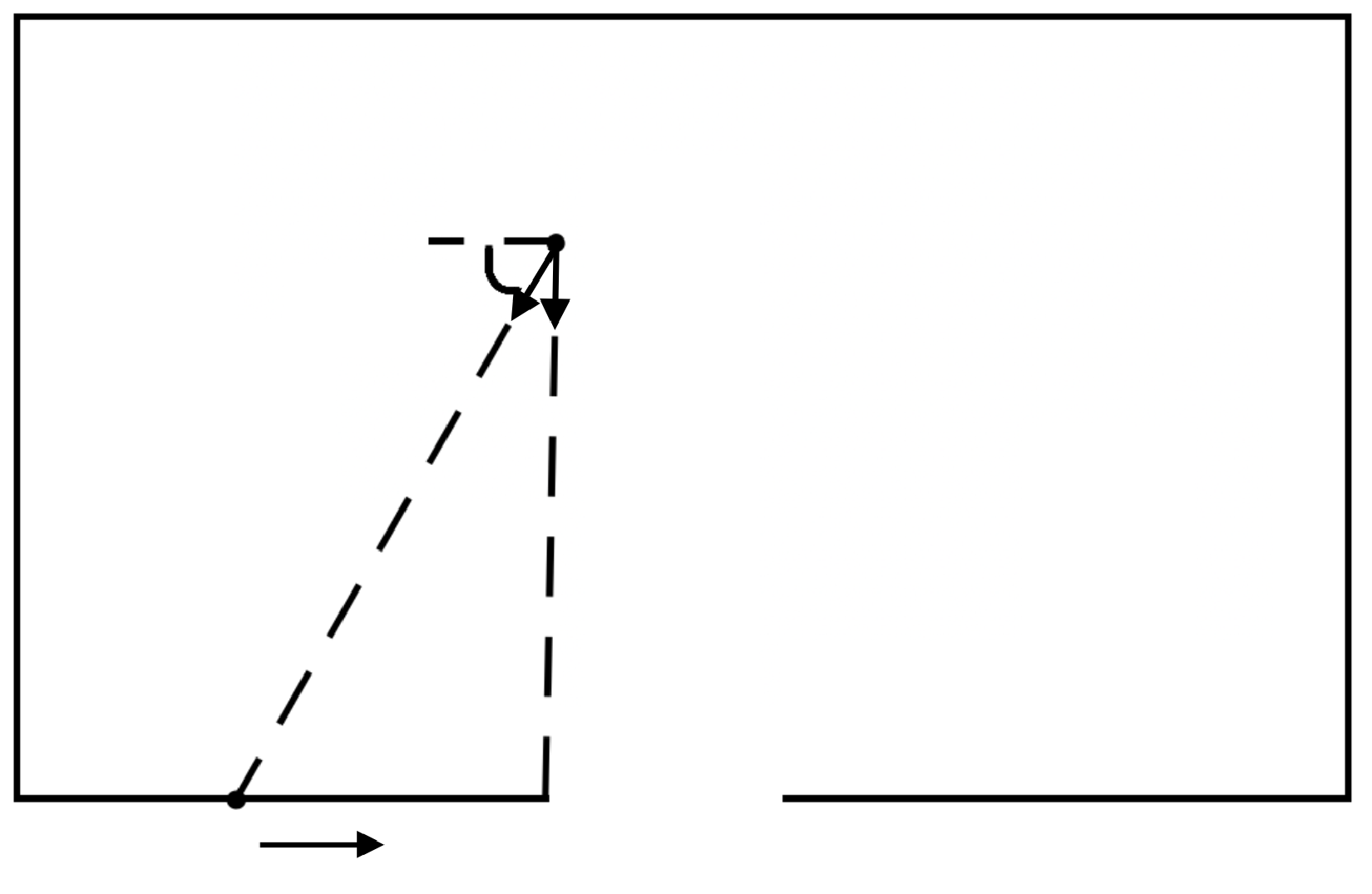}
\put(47, 5){{\text{$E$}}}
\put(75, 15){{\text{$\Omega$}}}
\put(42, 49){{\text{$\x$}}}
\put(30, 41){\footnotesize{\text{$\theta_h$}}}
\put(10, 3){\footnotesize{\text{$\x_{W}$}}}
\put(20, 0){\footnotesize{\text{$\u_{W}$}}}
\put(43, 44){\footnotesize{\text{$\u_{E}$}}}
\put(11, 24){\footnotesize{\text{$d_W(\x)$}}}
\put(42, 24){\footnotesize{\text{$d_E(\x)$}}}
\end{overpic}
\caption{
Sketch of computational domain $\Omega$ with exit $E$ and a person located at $\x$, 
walking
with direction $\theta_h$. The person should choose direction $\u_E$ to reach the exit
and direction $\u_W$ is to avoid collision with the wall. The person's distances form the exit and from the wall where they would collide
are $d_E$ and $d_W$, respectively.
}\label{domain}
\end{figure}

The transition probabilities $\mathcal{A}_h^i$ and $\mathcal{B}_{hk}^i$ 
are meant to account for the process 
through which a person decides the direction to take. 
Let us start from defining $\mathcal{A}_h$, which accounts for the following two factors in deciding the walking direction:
 \begin{itemize}
 \item[-] \textit{The goal to reach the exit.}\\
Given a candidate individual at point $\x$, we define their distance to the exit as
\[d_E(\x)= \min_{\x_E \in E} || \x-\x_E ||,\]
and we consider the unit vector $\u_E(\x)$, pointing from $\x$ to the exit. See Figure~\ref{domain}. 
\item[-] \textit{The desire to avoid collisions with the walls.}\\
Given a candidate individual at the point $\x$ walking 
with direction $\theta_h$, we define the distance 
$d_W(\x, \theta_h)$ from the individual to a wall at a point $\x_W(\x, \theta_h)$, 
where the person is expected to collide with the wall.
The unit tangent vector $\u_W(\x, \theta_h)$ to $\partial \Omega$ at $\x_W$
points to the direction of the the exit. Vector $\u_W$ is used to avoid a collision with the wall. See Figure~\ref{domain}.
 \end{itemize} 
Let $\theta_G$ be the \textit{geometrical preferred direction}, which accounts for the above 
two factors:
 \begin{align}\label{eq:uG}
\u_G(\x, \theta_h) &= \frac{(1-d_E(\x))\u_E(\x) + (1-d_W(\x, \theta_h))\u_W(\x, \theta_h)}{|| (1-d_E(\x))\u_E(\x) + (1-d_W(\x, \theta_h))\u_W(\x, \theta_h) ||} \cl
&= (\cos \theta_G, \sin \theta_G).
\end{align}
A candidate individual $h$ will update their direction by choosing
the angle closest to $\theta_G$. Then, we  define 
$\mathcal{A}_{h}^i$ as: 
\begin{equation}\label{eq:A}
\mathcal{A}_{h}^i=
\begin{cases}
1 - \frac{4}{\pi}|\theta_G - \theta_i|  & \text{if} \,\,\, |\theta_G - \theta_i | \leq \frac{\pi}{4}, \\
0  & \text{if} \,\,\, |\theta_G- \theta_i| > \frac{\pi}{4},            
\end{cases}
\end{equation}
with
\begin{equation}\label{eq:distance}
 d(\theta_p, \theta_q)=
\begin{cases}
|\theta_p - \theta_q|  & \text{if} \,\,\, |\theta_p - \theta_q | \leq \pi, \\
2\pi - |\theta_p- \theta_q|  & \text{if} \,\,\, |\theta_p- \theta_q| > \pi .             
\end{cases}
\end{equation} 

Transition probability  $\mathcal{B}^i_{hk}$ accounts for the following factors in deciding the walking direction:
\begin{itemize}
\item[-] \textit{The tendency to look for less congested areas.}\\
A candidate person $(\x, \theta_h)$ may decide to change direction in order to avoid congested areas. One way to do this is by considering the minimal directional derivative of the density at the point $\x$. 
We denote such a direction by unit vector $\u_C(\theta_h, \rho)=(\cos\theta_C, \sin\theta_C)$, with 
\[C=\argmin_{j \in \{h-1, h, h+1\}}\{\partial_j\rho(t, \x)\}.\]
In conditions of no-stress, people typically choose this walking strategy. 
\item[-] \textit{The tendency to follow the stream.}\\
A candidate person $h$ interacting with a field person $k$  
may decide to follow the field person and thus walk with their direction, 
denoted with unit vector $\u_F=(\cos\theta_k, \sin\theta_k)$. Following the stream is what people tend to do in stressful conditions, like in an emergency evacuation, and is a competing walking strategy to the search for less congested areas. 
\end{itemize}
To weight between these two walking strategies,  we introduce a parameter $\varepsilon \in [0,1]$: $\varepsilon=0$ corresponds
to searching for less congested areas only, 
while $\varepsilon=1$ corresponds to only following the stream.
Then, the \textit{interaction-based preferred direction} is given by:
\[\u_P(\theta_h, \theta_k, \rho)= \frac{\varepsilon\u_F+(1-\varepsilon)\u_C(\theta_h, \rho)}
{||\varepsilon\u_F+(1-\varepsilon)\u_C(\theta_h, \rho)||} = (\cos \theta_P, \sin \theta_P),
\]
and
\begin{equation}\label{eq:B}
\mathcal{B}_{hk}^i=
\begin{cases}
1 - \frac{4}{\pi}|\theta_P - \theta_i|  & \text{if} \,\,\, |\theta_P - \theta_i | \leq \frac{\pi}{4}, \\
0  & \text{if} \,\,\, |\theta_P- \theta_i| > \frac{\pi}{4}.            
\end{cases}
\end{equation}
We recall that $d(\cdot, \cdot)$ is defined in \eqref{eq:distance}. Since $\theta_P- \theta_i=\arccos{(\u_P \cdot \u_i)}$, \eqref{eq:B} can be easily computed from:
\begin{equation}\label{eq:2}
\mathcal{B}_{hk}^i=
\begin{cases}
1 - \frac{4}{\pi}|\arccos{(\u_P \cdot \u_i)}|  & \text{if} \,\,\, |\arccos{(\u_P \cdot \u_i)} | \leq \frac{\pi}{4}, \\
0  & \text{otherwise},            
\end{cases}
\end{equation}
where 
\[
\u_P \cdot \u_i = \frac{\varepsilon\u_F \cdot \u_i+(1-\varepsilon)\u_C \cdot \u_i}
{\sqrt{\varepsilon^2 + (1 - \varepsilon)^2 + 2 \varepsilon (1 - \varepsilon) \u_F \cdot \u_C}}.
\]
So, \eqref{eq:B} is piecewise differentiable in $\varepsilon$.

\begin{remark}
We define $\mathcal{A}_{h}^i$ and $\mathcal{B}_{hk}^i$ 
differently from, e.g., 
Refs.~\citenum{Agnelli2015,Bellomo2015_new,kim_quaini}. The change is dictated by the need for $\mathcal{B}_{hk}^i$ to be piecewise differentiable for the approach presented in 
Sec.~\ref{sec:KKT}. Then, the definition of 
$\mathcal{A}_{h}^i$ is changed accordingly.
The main difference between our definitions in \eqref{eq:A} and \eqref{eq:B} and those in 
the above-mentioned references is that we allow 
people to switch their walking direction to the closest
to the geometrical preferred direction, in the case of $\mathcal{A}_{h}^i$, or the
interaction-based preferred direction, in the case 
of $\mathcal{B}_{hk}^i$, while previously 
the change in direction for the candidate individual $h$ was restricted to $\theta_{h-1}, \theta_h, \theta_{h+1}$.
\end{remark}



Finally, let us introduce notation to write model \eqref{eq:model} in compact form.
Let $\f \in \mathbb{R}^{N_d}$ be the vector function whose $i$-th entry is $f^i$, $\boldsymbol{\mathcal{A}} \in \mathbb{R}^{N_d \times N_d} $ the matrix whose $ih$ 
entry is $\mathcal{A}_h^i$, $\boldsymbol{\mathcal{B}}^i \in \mathbb{R}^{N_d \times N_d}$ 
the matrix whose $hk$ entry is $\mathcal{B}_{hk}^i$, and $\b \in \mathbb{R}^{N_d}$ 
the vector function whose $i$-th entry is $\f^T\boldsymbol{\mathcal{B}}^i \f$.
Furthermore, we define $\J = \v \otimes \f$. Then, model \eqref{eq:model} 
can be written in vector form as follows:
\begin{align}
\frac{\partial \f}{\partial t} &+ \nabla \cdot \J  = (1 - \rho) (\boldsymbol{\mathcal{A}} \f - \f) + \rho (\b - \rho \f). \label{eq:model_vec}
\end{align}

\subsection{The time-continuous inverse problem}\label{sec:IP}

Let us assume that we can extract the density data from the video frames every $\tau_k$, $k = 1, \dots, N$, within the time interval of interest $[0, T]$, where $t =0$ is the time when the insect repellent is injected, and $T = \tau_N$. In addition, we assume that the time step $\Delta \tau$ to acquire the data is constant. 
Although simplifying, this assumption is neither unrealistic nor restrictive. We denote by $\rho_{v}(\tau_k, \x)$ the density extracted from the video footage at time $\tau_k$ for $\x \in \Omega$.

We introduce the following functional:
\begin{align}
\mathcal{J}_c = \frac{1}{2} \sum_{k = 1}^N \int_\Omega (\rho(\tau_k, \x) - \rho_{v}(\tau_k, \x))^2 d\x, \label{eq:J}
\end{align}
where $ \rho(\tau_k,\x)$ is computed with \eqref{eq:rho} from the solution of \eqref{eq:model_vec}.
Usually, a regularization term is added to functional $\mathcal{J}_c$:
\begin{align}
\mathcal{J}^\mathcal{R}_c = \mathcal{J}_c  + \mathcal{R}_c(\varepsilon), \label{eq:JR}
\end{align}
where $\mathcal{R}_c$ could be Tikhonov regularization, which is expected to improve the mathematical
and numerical properties of the minimization process. It could be defined as:
\begin{align}
\mathcal{R}_c(\varepsilon) = \frac{\xi}{2} \sum_{k = 1}^N \int_\Omega (\varepsilon(\tau_k, \x) - \varepsilon_{ref}(\tau_k, \x))^2 d\x, \label{eq:R}
\end{align}
where $\varepsilon_{ref}$ is a reference distribution of the stress level. With this choice of
regularization, we are forcing $\varepsilon$ to be close to $\varepsilon_{ref}$. 
Parameter $\xi$, which weights the importance of the regularization term in the minimization 
procedure, can be tuned empirically. 

Let $\mathcal{E}$ be an admissible set for the control variable $\varepsilon$. We define it as
\begin{align}
\mathcal{E} = \{ \varepsilon~:~\varepsilon \in L^\infty(\Omega), ~0 \leq \varepsilon \leq 1 \}. \label{eq:E}
\end{align}
Then, we can define the inverse crowd dynamics problem: 
\vskip 2mm
\noindent \underline{Time-continuous problem}:~\emph{For $\x \in \Omega$, find $\varepsilon = \varepsilon(t,\x) \in \mathcal{E}$ 
that minimizes functional $\mathcal{J}^\mathcal{R}_c$ \eqref{eq:JR} under constraints \eqref{eq:rho}, \eqref{eq:model_vec}.}
\vskip 2mm

This is a time-dependent minimization problem. 
A possible approach for its solution is the introduction of a Lagrange multiplier for the constraint given by \eqref{eq:model_vec} and a Lagrangian functional to be minimized with no constraint. Then, the Lagrangian functional is differentiated with respect to the state variable $\f$, the Lagrange multiplier, and the control variable $\varepsilon$ to obtain the classical Karush–Kuhn–Tucker (KKT) system, i.e., state problem, adjoint problem, and optimality condition. 

While feasible, this approach has two important
drawbacks, both related to the fact that 
the adjoint of an initial value problem is a final value problem. 
The first is complexity because the solution of a final
value problem is nontrivial. The second is the need
to store the solution at all the time steps 
due to the back-in-time nature of the adjoint problem. 
Since we prefer to avoid these drawbacks, 
we follow an alternative approach: we first discretize in time and then minimize the time-discrete problem. The time-discrete problem is presented in Sec.~\ref{sec:td_FP} and the new minimization problem that arises from this \textit{discretize-then-optimize} approach is stated in Sec.~\ref{sec:KKT}.

\section{The time discrete inverse crowd dynamics problem}\label{sec:td_model}

\subsection{The time discrete forward problem}\label{sec:td_FP}

Let $\Delta t$ be a constant time step. 
It is reasonable to assume that $\Delta \tau$,
i.e.,~the time step used to extract the data 
from the video (see Sec.~\ref{sec:IP}), is larger than $\Delta t$, which might have to fulfill stability and accuracy requirements. We assume that $\Delta \tau$ is a multiple $m$ of $\Delta t$, i.e. $\Delta \tau = m \Delta t$.
Should this assumption not hold true, suitable interpolation procedures can be used to interpolate the data at the time discretization nodes. For simplicity of exposition, from now on we will assume that $m = 1$. The case of $m > 1$ would not require a fundamental change in the procedure presented below.

We set $t^n = n \Delta t$, $n = 0, \dots, N$ with $N \Delta t = T$. We denote by $y^n$ the approximation of a generic
time-dependent function $y$ at time $t^n$. For the time discretization of problem \eqref{eq:model_vec}, we consider an explicit method. Although a CFL condition needs to be satisfied for stability, we choose this method for its simplicity and to build upon 
our previous work,\cite{kim_quaini,kim_quaini2020,kim_quaini2020b} which uses this method. 
Of course, the approach presented below can be adapted to other time-discretization schemes.

Problem \eqref{eq:model_vec} discretized in time with a first order
explicit method reads: at time $t^{n+1}$, with $n = 0, \cdots, N-1$, given $\f^n$, find $\f^{n+1}$ such that
\begin{align}
&\frac{\f^{n+1}}{\Delta t} - \rho^n \b(\f^n,\varepsilon^{n+1}) =  \g^n\label{eq:model_vec_td} \\
& \g^n = \frac{\f^{n}}{\Delta t} -\nabla \cdot \J^n + (1 - \rho^n) (\boldsymbol{\mathcal{A}} \f^n - \f^n) - (\rho^n)^2 \f^n, \el
\end{align}
where $\b(\f^n,\varepsilon^{n+1})$ is the vector whose $i$-th entry is $(\f^n)^T\boldsymbol{\mathcal{B}}^i(\varepsilon^{n+1}) \f^n$ and $\J^n = \v^n \otimes \f^n$. Notice that
we have made explicit the dependence of vector $\b$ on $\f^n$, i.e., the probability distribution function at the previous time step, and $\varepsilon^{n+1}$, i.e., the fear level (our control variable)
at the current time step. The right-hand-side 
in \eqref{eq:model_vec_td} depends only on quantities computed at the previous time step. 
Once $\f^{n+1}$ is obtained from \eqref{eq:model_vec_td}, the density at time $t^{n+1}$ is found by summing the components of $\f^{n+1}$:
\begin{align}\label{eq:rhon}
\rho^{n+1}=\sum_{i=1}^{N_d}f^{i, n+1} .
\end{align}

\subsection{The time discrete inverse problem and the KKT conditions}\label{sec:KKT}

Let us consider functional:
\begin{align}
    \mathcal{J}^n = \frac{1}{2} \int_\Omega (\rho(t^n, \x) - \rho_{v}(t^n, \x))^2 d\x. \label{eq:td_J}
\end{align}
As mentioned in  
Sec.~\ref{sec:IP}, it is typical to add a a regularization term $\mathcal{R}$. 
If one has chosen regularization \eqref{eq:R},
it could be defined as follows for functional \eqref{eq:td_J}:
\begin{align}
\mathcal{R}^{n}(\varepsilon) = \frac{\xi}{2} \int_\Omega (\varepsilon(t^{n}, \x) - \varepsilon_{ref}(t^{n}, \x))^2 d\x. \label{eq:td_J_R}
\end{align}
Then, at every time $t^{n+1}$, for $n = 0, \dots, N -1$, we have to solve the following new minimization problem:
\vskip 2mm
\noindent \underline{Time-discrete variational problem}:~\emph{Find $\varepsilon \in \mathcal{E}$ that minimizes functional
\begin{align}
\mathcal{J}^{\mathcal{R},n+1} = \mathcal{J}^{n+1}  + \mathcal{R}^{n+1}(\varepsilon), \label{eq:JR_dt}
\end{align}
under constraints \eqref{eq:model_vec_td}-\eqref{eq:rhon}. 
}
\vskip 2mm

To find the minimizer of \eqref{eq:td_J_R} under constraints \eqref{eq:model_vec_td}-\eqref{eq:rhon}, we adopt the classical Lagrange multiplier approach. In this method, the 
problem of finding the minimizer of a functional subject to one or more constraints is reformulated as a saddle point problem. The solution corresponding to the original constrained 
minimization problem is a saddle point of the so-called Lagrangian functional. For problem \eqref{eq:td_J_R}, 
the Lagrangian functional is given by:
\begin{align}
    \mathcal{L}(\F, \bLam, E) = \mathcal{J}^{n+1}+\mathcal{R}^{n+1}(E)
    + \left\langle \frac{\F}{\Delta t} - \rho^n \b(\f^n,E) -  \g^n, \bLam \right\rangle, \label{eq:L}
\end{align}
where $\langle \cdot, \cdot \rangle$ denotes
an inner product and $\bLam$ is the Lagrange multiplier. Since in our previous work\cite{kim_quaini, kim_quaini2020, kim_quaini2020b} we have used a Finite Difference (FD) method for space discretization, we use the dot product 
as the inner product and \eqref{eq:L} becomes
\begin{align}
    \mathcal{L}(\F, \bLam, E) = \mathcal{J}^{n+1}+\mathcal{R}^{n+1}(E)
    + \left( \frac{\F}{\Delta t} - \rho^n \b(\f^n,E) -  \g^n\right) \cdot \bLam. \label{eq:L_dot}
\end{align}
If one were to use a different method for space
discretization (e.g., a Finite Element method), 
a different inner product would be used 
(e.g., the $L^2$ inner product). 

A saddle point of the Lagrangian functional can be identified among the stationary points. To find such stationary 
points, all the partial derivatives need to vanish:
\begin{align}
    &\frac{\partial \mathcal{L}}{\partial \bLam} \Big|_{\F = \f^{n+1}}  = \frac{\f^{n+1}}{\Delta t} - \rho^n \b(\f^n,\varepsilon^{n+1}) -  \g^n = \boldsymbol{0}, \label{eq:state} \\
    &\frac{\partial \mathcal{L}}{\partial \F} \Big|_{\bLam = \blam^{n+1}} = \int_\Omega  (\rho^{n+1} - \rho_{v}^{n+1}) \frac{\partial \rho^{n+1}}{\partial \F} [ \f^{n+1}] ~d\x +
    \frac{\blam^{n+1}}{\Delta t} = \boldsymbol{0}, \label{eq:adjoint} \\
    &\frac{\partial \mathcal{L}}{\partial E} \Big|_{E = \varepsilon^{n+1}} =\xi \int_\Omega (\varepsilon^{n+1} - \varepsilon_{ref}(\tau_{n+1}, \x)) d\x - \rho^n \frac{ \partial \b}{\partial E} [\varepsilon^{n+1}] \cdot \blam^{n+1} = 0. \label{eq:opt_cond} 
\end{align}
Eq.~\eqref{eq:state} represents the state problem, 
eq.~\eqref{eq:adjoint} the adjoint problem, 
and \eqref{eq:opt_cond} is the optimality condition. So, at time $t^{n+1}$, we have the following coupled KKT system: find  variable $\f^{n+1}$, Lagrange multiplier $\blam^{n+1}$, and parameter $\varepsilon^{n+1}$
such that eq.~\eqref{eq:state}-\eqref{eq:opt_cond}
hold, with density $\rho^{n+1}$ related to $\f^{n+1}$ through \eqref{eq:rhon}. Note that the state problem \eqref{eq:state}
is exactly problem \eqref{eq:model_vec_td}. 
If one chooses not to regularize the minimization of
\eqref{eq:td_J}, which is equivalent to setting $\xi = 0$,
the optimality condition \eqref{eq:opt_cond} becomes
a requirement of orthogonality for
$\frac{ \partial \b}{\partial E} [\varepsilon^{n+1}]$ and $\blam^{n+1}$.

\section{Numerical solution of the inverse problem}\label{sec:num_sol}

A possible algorithm for the 
solution of coupled nonlinear system \eqref{eq:state}-\eqref{eq:opt_cond} is a steepest descent method. Let $k$ denote the steepest descent method iteration. 
The algorithm reads as follows: at time $t^{n+1}$, iteration $k + 1$, assuming that $\f^n$ and $\varepsilon_{k}$ are known, perform the following steps
\begin{itemize}
    \item[-] Step 1: Solve state problem
    \begin{align}\label{eq:statek}
 \frac{\f_{k+1}}{\Delta t} - \rho^n \b(\f^n,\varepsilon_{k}) =  \g^n.
\end{align}
With $\f_{k+1}$, one reconstructs the approximation for the density $\rho_{k+1}$.
\item[-] Step 2: Solve adjoint problem
\begin{align}\label{eq:adj}
    \frac{\blam_{k+1}}{\Delta t} = -\int_\Omega  (\rho_{k+1} - \rho_{v}^{n+1}) \frac{\partial \rho_{k+1}}{\partial \F} [ \f_{k+1}] ~d\x. 
\end{align}
\item[-] Step 3: Given a suitable $\delta_{k+1} \in \mathbb{R}$, set 
\begin{align}\label{eq:update}
\varepsilon_{k+1} - \varepsilon_{k}  = \delta_{k+1} \left(\xi \int_\Omega (\varepsilon_{k} - \varepsilon_{ref}(\tau_{k}, \x)) d\x - \rho^n \frac{ \partial \b}{\partial E} [\varepsilon_k] \cdot \blam_{k+1} \right). 
\end{align}
\item[-] Step 4: Check stopping criterion 
\begin{align}\label{eq:stop}
\frac{1}{| \Omega |}\left\|\xi \int_\Omega (\varepsilon_{k} - \varepsilon_{ref}(\tau_{k}, \x)) d\x - \rho^n \frac{ \partial \b}{\partial E} [\varepsilon_k] \cdot \blam_{k+1} \right\|_{L^2(\Omega)} < tol,
\end{align}
where $tol$ is a given stopping tolerance. If satisfied, set $\f^{n+1}_h = \f_{k+1}$, $\rho^{n+1}_h = \rho_{k+1}$, and $\varepsilon^{n+1} = \varepsilon_{k+1}$. Otherwise, repeat steps 1 to 3.
\end{itemize}

\subsection{Space discretization} \label{sec:discretization}

For simplicity, we present space discretization using a 
rectangular computational $[0, L] \times [0, H]$, for given $L$ and $H$. 
We generate a structured mesh of this domain by choosing $\Delta x$ and $\Delta y$ 
to partition interval $[0, L]$
and $ [0, H]$, respectively. Let $N_x = L/\Delta x$ and $N_y = H/\Delta y$.
We define the discrete mesh points $\x_{pq} = (x_p, \, y_q)$ by
\begin{align}
x_p &=p \Delta x, \quad p= 0, 1, \dots, N_x,    \cl
y_q  &=q \Delta y, \quad q= 0, 1, \dots, N_y.    \el
\end{align}
We also need to define:
\begin{align}
x_{p+1/2}=x_{p}+\Delta x/2=\Big(p+\frac{1}{2}\Big)\Delta x, \cl
y_{q+1/2}=y_{q}+\Delta y/2=\Big(q+\frac{1}{2}\Big)\Delta y. \el
\end{align}

The space discretization of the state problem \eqref{eq:statek}
is thoroughly described in Ref.~\citenum{kim_quaini}. Here, we present the 
space discretization of the adjoint problem \eqref{eq:adj} and
the update of the stress level from \eqref{eq:update}. Let us start form the
problem \eqref{eq:adj}.

Recalling that the 
density is computed from the probability functions 
through \eqref{eq:rhon}, it is easy to see that 
\begin{align}
\frac{\partial \rho_{k+1}}{\partial \F} [ \f_{k+1}] = \mathbbm{1},
\end{align}
where $\mathbbm{1} \in \mathbb{R}^{N_d}$ is a vector with all entries equal to 1. 
Then, we have
\begin{align}\label{eq:adj_simpl}
    \frac{\lambda_{k+1}^i}{\Delta t} = -\int_\Omega  (\rho_{k+1} - \rho_{v}^{n+1})~d\x, \quad i = 1, \dots, N_d, 
\end{align}
i.e., all entries of $\blam_{k+1}$ are equal. 
This is due to the choice of an explicit time discretization scheme. If we were to, e.g., treat the 
convective in \eqref{eq:model_vec} implicitly, 
the components of $\blam_{k+1}$, which are
associated to the different walking directions, would
be different. 
The FD method we use produces an approximation of $\lambda_{k+1}^i$ with cell average
$\Lambda^{i,pq}_{k+1} \in \mathbb{R}$:  
\[
\Lambda^{i,pq}_{k+1} \approx \dfrac{1}{\Delta x \, \Delta y} \int_{y_{q-1/2}}^{y_{q+1/2}}  \int_{x_{p-1/2}}^{x_{p+1/2}}
\lambda_{k+1}^i(x, y) dx\, dy, 
\]
where $1 \leq p \leq N_x-1$ and $1 \leq q \leq N_y-1$. From \eqref{eq:adj_simpl},
we get:
\begin{align}
    \Lambda^{i,pq}_{k+1} = - \Delta t  \int_{y_{q-1/2}}^{y_{q+1/2}}  \int_{x_{p-1/2}}^{x_{p+1/2}}
(\rho_{k+1} - \rho_{v}^{n+1}) dx\, dy.
\end{align}
The FD approximation of $\lambda_{k+1}^i$, denoted by $\Lambda_{k+1}^i$,  is given by a piece-wise constant function 
in space that takes value $\Lambda^{i,pq}_{k+1}$ over cell $[x_{p-1/2}, x_{p+1/2}]\times [y_{q-1/2}, y_{q+1/2}]$.

Now, let us take care of Step 3, which updates the value of the fear level. 
The FD approximation of eq.~\eqref{eq:update}, using the same FD method used for \eqref{eq:adj} reads:
\begin{align}
E_{k+1}^{pq} - E_{k}^{pq}  = &\delta_{k+1} \xi^{pq} \int_{y_{q-1/2}}^{y_{q+1/2}}  \int_{x_{p-1/2}}^{x_{p+1/2}}(\varepsilon_{k} -   \varepsilon_{ref}(\tau_{k}, \x)) dx\, dy  \cl
&- \delta_{k+1}P^{n,pq} \frac{ \partial \b^{pq}}{\partial E} [\varepsilon_k] \cdot \bLam^{pq}_{k+1}, \label{eq:update_disc}
\end{align}
where $E_{k+1}^{pq}$ and $E_{k}^{pq}$ are (subsequent) cell averages of $\varepsilon_{k+1}$ and 
$\varepsilon_{k}$ for the generic $pq$ cell, $P^{n,pq}$ is the cell average of $\rho^{n}$, and 
$\bLam^{pq}_{k+1} \in \mathbb{R}^{N_d}$ is the vector whose generic $i$-th component 
is $\Lambda^{i,pq}_{k+1}$. Vector $\frac{ \partial \b^{pq}}{\partial E} [\varepsilon_k]$ in \eqref{eq:update_disc} requires more
explanation. Recall that $\b(\f^n, E) \in \mathbb{R}^{N_d}$ is the vector whose $i$-th entry is $(\f^n)^T\boldsymbol{\mathcal{B}}^i (E) \f^n$, where 
$\boldsymbol{\mathcal{B}}^i (E) \in \mathbb{R}^{N_d \times N_d}$ is the matrix whose $hk$ entry 
is $\mathcal{B}_{hk}^i$ defined in \eqref{eq:B}. So, the $i$-th entry of FD approximation of
$\frac{ \partial \b(\f^n, E)}{\partial E} [\varepsilon_k]$ is:
\begin{align}
    (\F^{n,pq})^T \frac{\partial \boldsymbol{\mathcal{B}}^i}{\partial E}[\varepsilon_k] \F^{n,pq},
\end{align}
where $\F^{n,pq}$ is the cell average of $\f^{n}$. 
From 
\eqref{eq:2}, it is easy to see that
entry $hk$ of matrix $\frac{\partial \boldsymbol{\mathcal{B}}^i}{\partial E}[\varepsilon_k]$ is
\begin{equation}\label{eq:2a}
\frac{\partial \mathcal{B}_{hk}^i}{\partial E} [\varepsilon_k]=
\begin{cases}
\frac{4}{\pi}\frac{1}{\sqrt{1-(\u_P \cdot \u_i)^2}} \frac{\partial (\u_P \cdot \u_i)}{\partial E} [\varepsilon_k] & \text{if} \,\,\, |\arccos{(\u_P \cdot \u_i)} | \leq \frac{\pi}{4}, \\
0  & \text{otherwise},            
\end{cases}
\end{equation}
with 
\begin{align*}
\frac{\partial (\u_P\cdot\u_i)}{\partial E} [\varepsilon_k] = \frac{\varepsilon_k(\u_F\cdot \u_C-1)(\u_F\cdot \u_i+\u_C\cdot \u_i)+\u_F\cdot \u_i-(\u_C\cdot \u_i ) (\u_F\cdot \u_C)}
{\left( -2\varepsilon^2_k(\u_F \cdot \u_C)+ 2\varepsilon_k\u_F \cdot \u_C+ 2\varepsilon^2_k-2\varepsilon_k +1 \right)^{3/2}}.
\end{align*}

\section{Numerical results}\label{sec:numericalresults}
The numerical experiments described in this section are inspired from 
the lab experiments on ants reported in Ref.~\citenum{SHIWAKOTI20111433}, which were
aimed at understanding the effects of structural features on collective movement patterns
during rapid egress. Of the four scenarios studied in Ref.~\citenum{SHIWAKOTI20111433}, 
we consider the following three:
\begin{itemize}
    \item[-] Ants escaping from a circular chamber without partial obstruction near exit in Sec.~\ref{sec:circular}.
    \item[-] Ants escaping from a circular chamber with partial obstruction (via a column) near exit is
    Sec.~\ref{sec:circular_col}.
    \item[-] Ants escaping from a square chamber with exit at the corner of the walls in Sec.~\ref{sec:square}.
\end{itemize}
The exit in each chamber is such that it allows unimpeded passage of a single ant
or somewhat encumbered passage of two ants at the same time. 

Several repetitions of these scenarios were performed by the authors of Ref.~\citenum{SHIWAKOTI20111433}.
In each experiment, the ants (numbering 200–250) were allowed to nest naturally in
the chamber. Then,  citronella (an insect repellent liquid) 
was injected through a small hole in the chamber ceiling to create
panic. Upon the injection of the liquid, ants rushed to the exit in 
a manner reminiscent of humans in a crowd panic.

We start our simulations right after the citronella has been injected. 
Ref.~\citenum{SHIWAKOTI20111433} reports the mean time for the
first 50 ants to escape the chamber. The data were obtained by manual counting from playback of digital video recordings. Our simulations try to match those averages by correcting 
the fear level through the optimization
procedure described in Sec.~\ref{sec:num_sol}.

As mentioned in Sec.~\ref{sec:forward}, we work with 
nondimensional quantities. For this purpose, we need 
to define: 
\begin{itemize}
\item[-] a characteristic length $D$, which is the largest distance an active particle can cover in domain $\Omega$;
\item[-] a maximum moving speed $V_M$;
\item[-] a reference time $T = D/V_M$;
\item[-] a maximum admissible number of active particles per unit area $\rho_M$.
\end{itemize}
For all the tests, we set $V_M = 2$ mm/s
and $\rho_M = 0.5$ ants/mm$^2$, while the values of $D$ and thus $T$ vary depending on the geometry. 
For the circular chamber, we set $D =35$ mm, while
for the square chamber we take $D =31 \sqrt{2}$ mm.
Then, the characteristic times are found by plugging the given $D$ in $T= D/V_M$.

\subsection{Circular chamber without partial obstruction near exit}\label{sec:circular}

\begin{figure}[htb!]
\centering
\begin{overpic}[width=0.47\textwidth,grid=false]{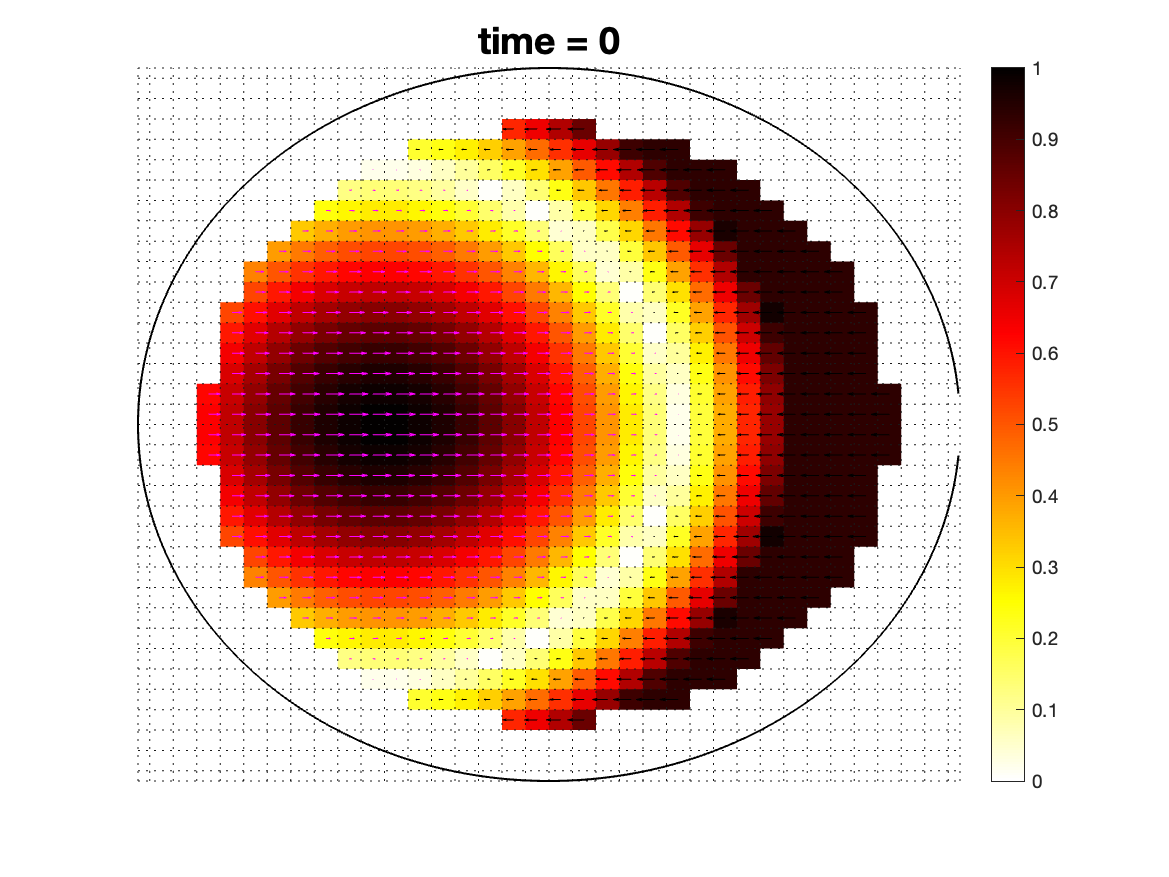}
\end{overpic}
\begin{overpic}[width=0.33\textwidth,grid=false]{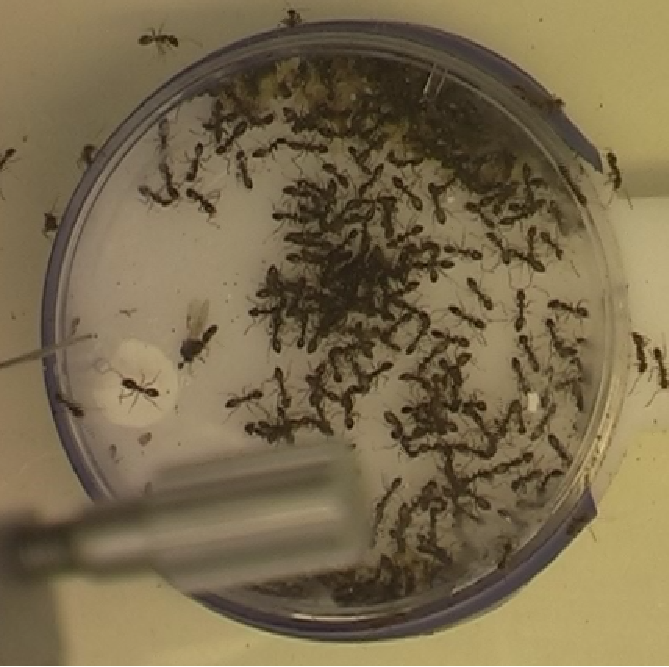}
\end{overpic}
\caption{Left: Computational domain for the circular chamber without column, with initial density and direction for the numerical simulation. Right:
snapshot from repetition 2 out of 30.}
\label{fig:initialcircle1}
\end{figure}

The  circular chamber has a diameter of 
35 mm with one exit located in the middle of the right side. The exit size is 2.5 mm. See  Fig.~\ref{fig:initialcircle1} (left).
In this chamber, we introduce 200 ants. 
As shown in Fig.~\ref{fig:ants}, ants tend to congregate
in parts of the chamber, rather than being
uniformly distributed in it. We initially distribute
the 200 ants into two distinct groups: a circular group 
located towards the opposite end to the exit
and a crescent-shaped group near the exit. See Fig.~\ref{fig:initialcircle1} (left).
A similar distribution is observed in one of
the experiments from Ref.~\citenum{SHIWAKOTI20111433}, see Fig.~\ref{fig:initialcircle1} (right), 
where we also see the device used to inject the citronella.
We assign initial moving direction  $\theta_1$
to the ants in the circular group and the opposite direction 
(i.e., $\theta_5$) to the ants in the crescent-shaped group.

\begin{figure}[htb!]
\centering
\begin{overpic}[width=0.32\textwidth]{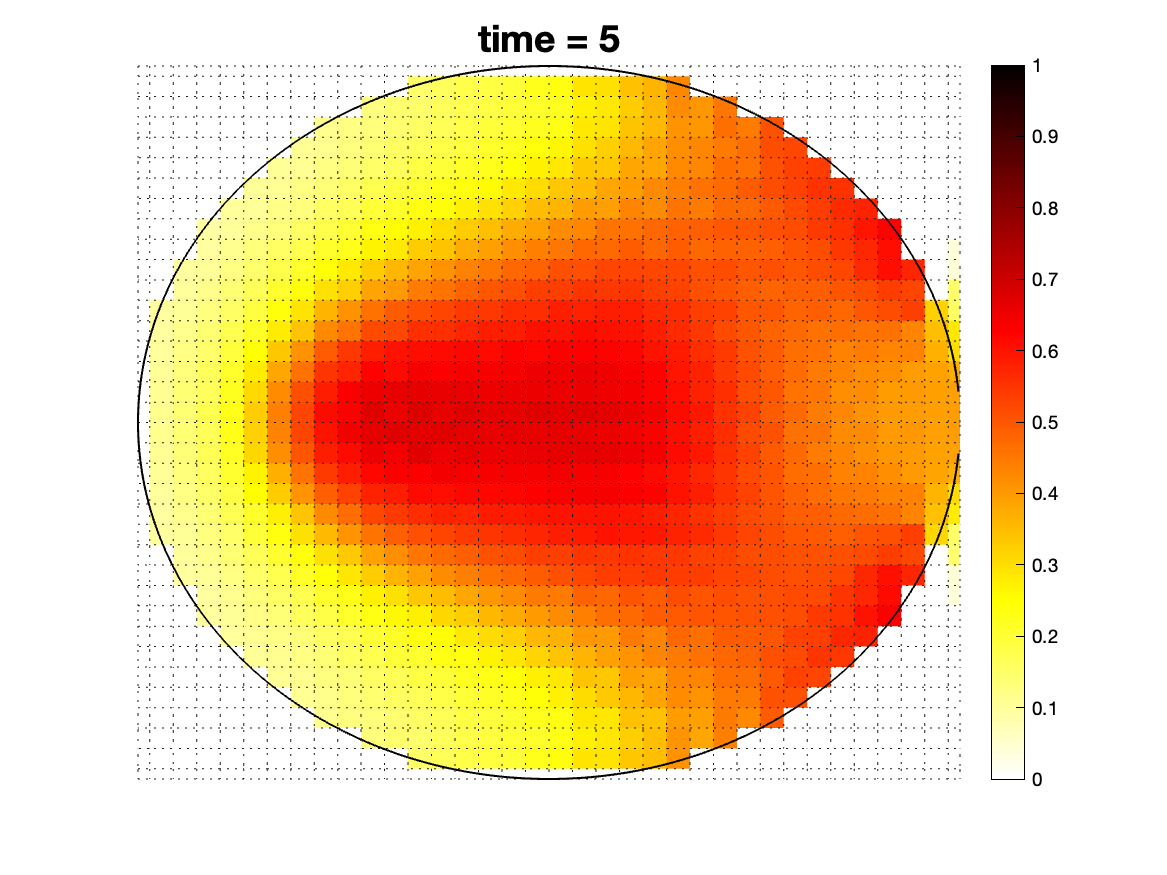}
\end{overpic} 
\begin{overpic}[width=0.32\textwidth,grid=false]{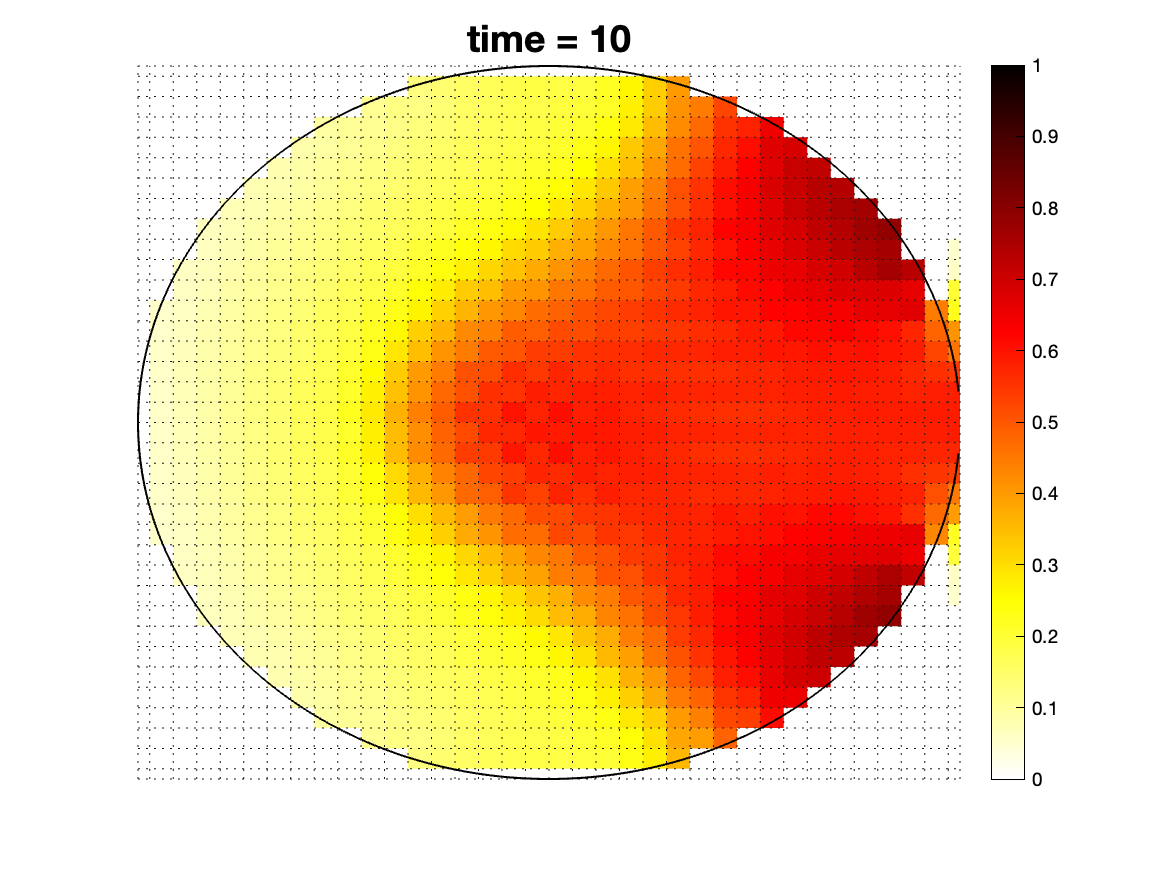}
\put(30,78){\textcolor{black}{$\varepsilon=0.05$}}
\end{overpic} 
\includegraphics[width=0.32\textwidth]{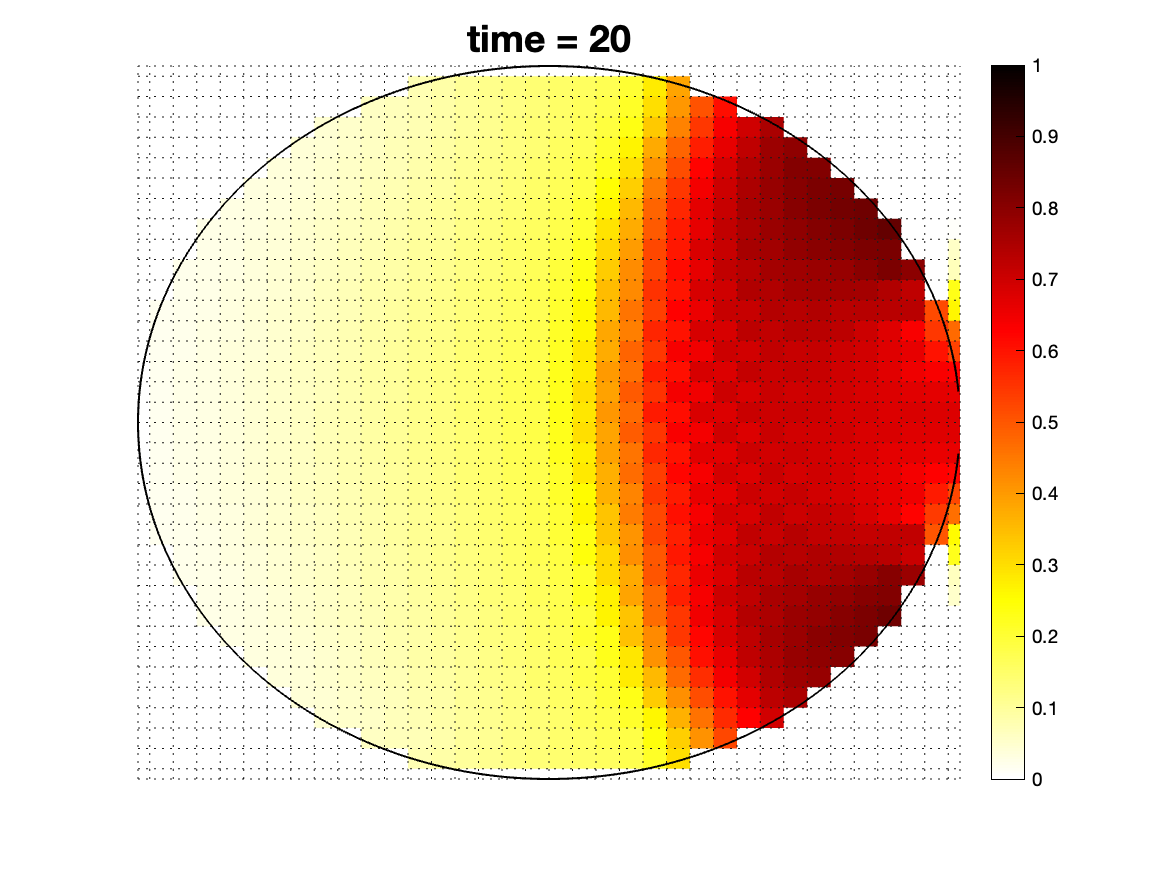} \\
\vskip .4cm
\begin{overpic}[width=0.32\textwidth]{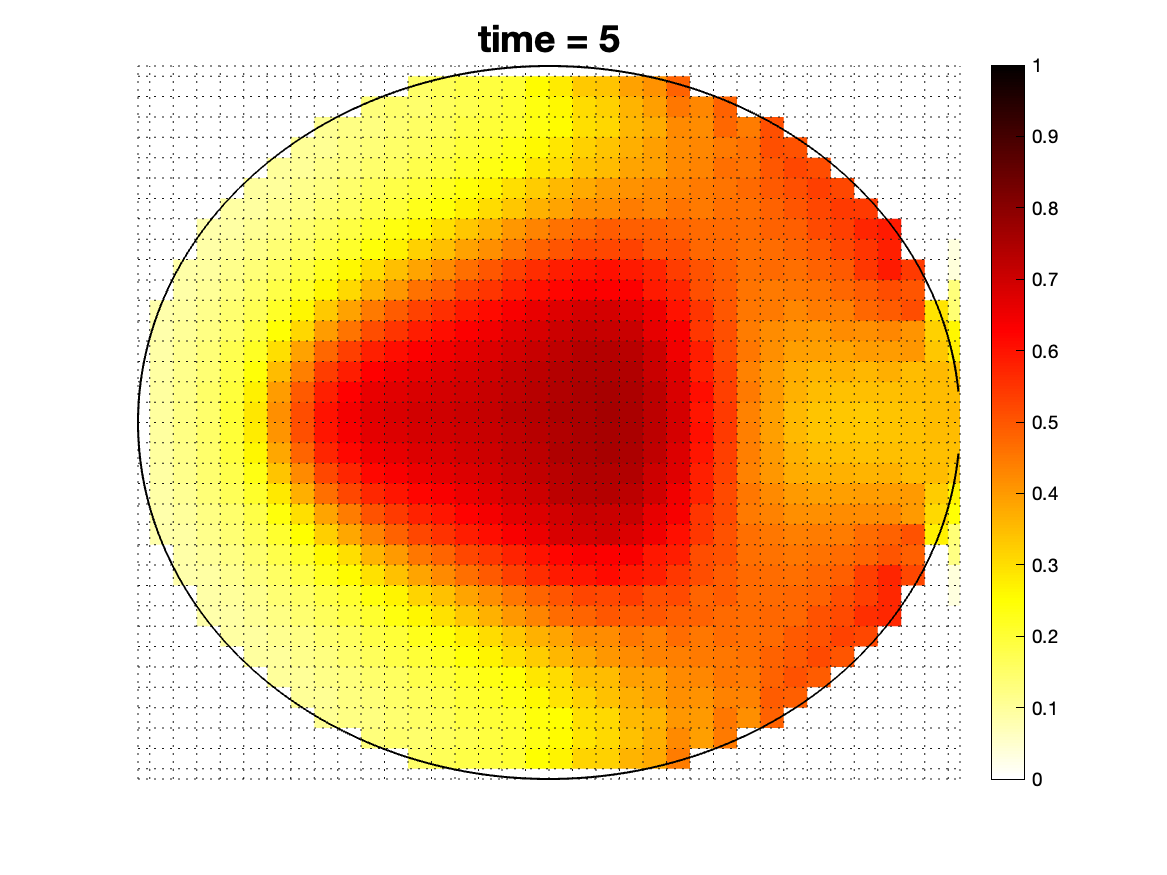}
\end{overpic} 
\begin{overpic}[width=0.32\textwidth]{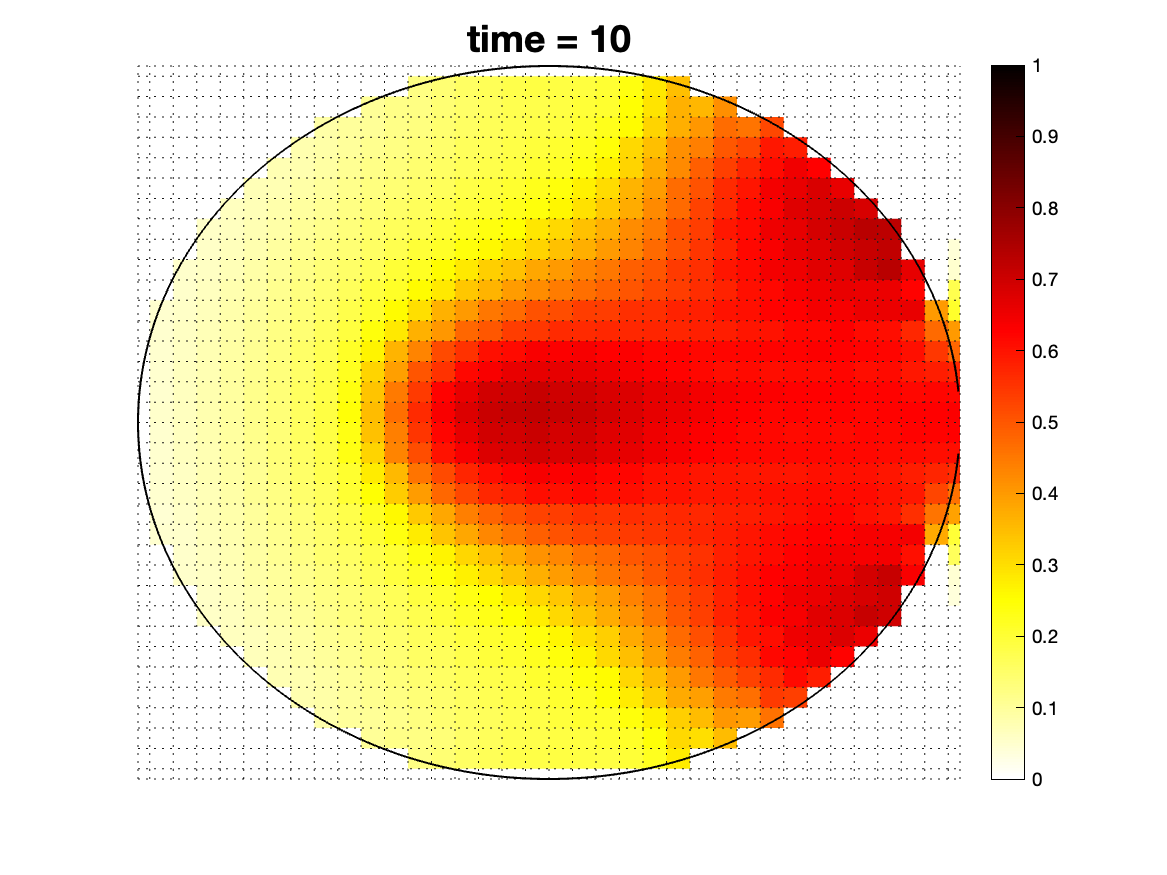}
\put(30,78){\textcolor{black}{$\varepsilon=0.95$}}
\end{overpic} 
\includegraphics[width=0.32\textwidth]{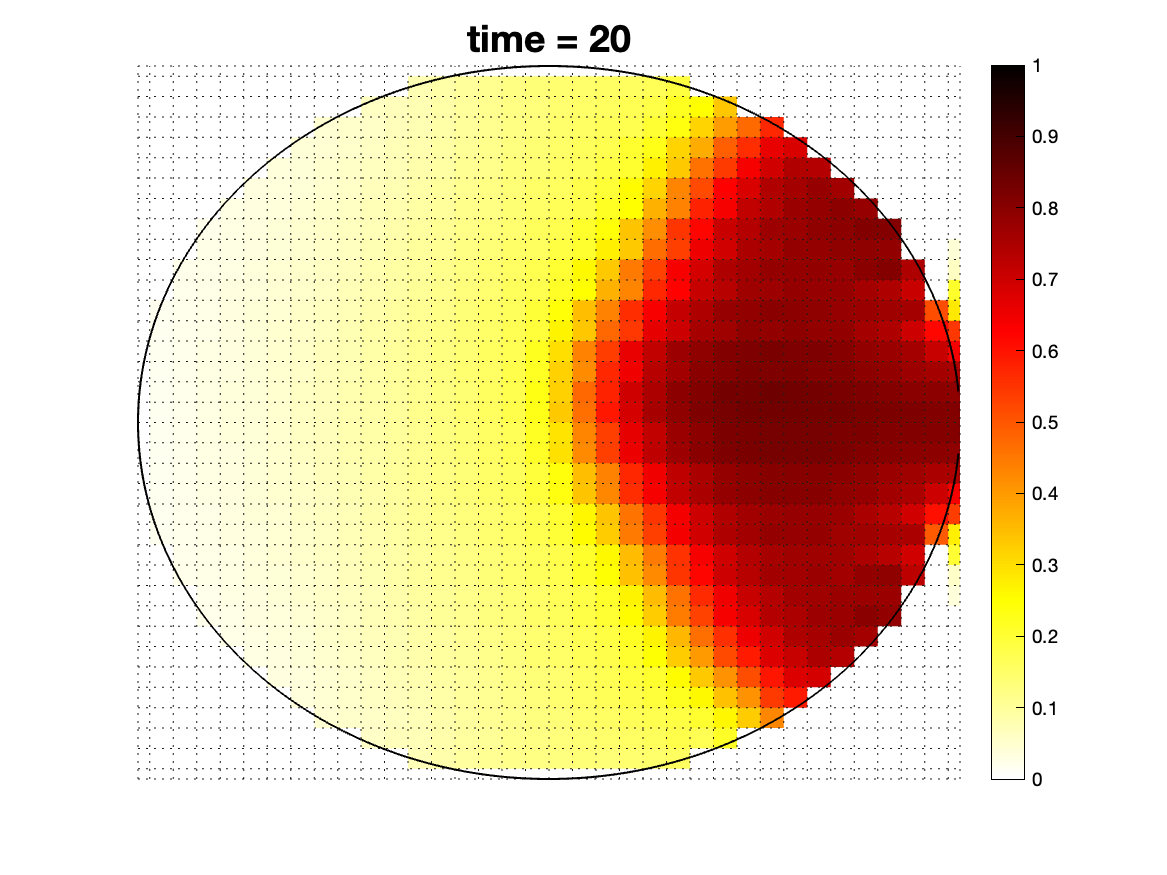}
\caption{Computed density at $t = 5$ s (left), $t = 10$ s (center), and $t = 20$ s (right)
given by the forward problem 
for 200 ants initially placed as
in Fig.~\ref{fig:initialcircle1} (left)
for $\varepsilon=0.05$ (top row) and $\varepsilon=0.95$ (second row).}
\label{fig:circular_2eps}
\end{figure}

We chose the initial configuration described above because it leads to
visibly different evacuation dynamics
in conditions of low stress ($\varepsilon = 0.05$)
vs high stress ($\varepsilon = 0.95$). See 
Fig.~\ref{fig:circular_2eps}.
We see that when the stress level is higher, 
regions of high density form around the center
of the chamber at $t = 5, 10$ s. Additionally, 
at $t = 20$ the crowd density is high in an arrow-shaped region for $\varepsilon=0.95$, while 
for $\varepsilon=0.05$ the edge of the higher-density region is smoother. 

\begin{figure}[htb!]
\centering
\begin{overpic}[width=0.32\textwidth]{Circle_Room_Ped_Eps0_95_Ants_circle_time_10.eps}
\end{overpic} 
\begin{overpic}[width=0.32\textwidth]
{Circle_Room_Ped_Eps0_95_Ants_circle_time_20.eps}
\put(-10,78){\textcolor{black}{synthetic video data ($\varepsilon=0.95$)}}
\end{overpic} 
\includegraphics[width=0.32\textwidth]{Circle_Room_Ped_Eps0_95_Ants_circle_time_40.eps}
\\
\vskip .4cm
\begin{overpic}[width=0.32\textwidth]{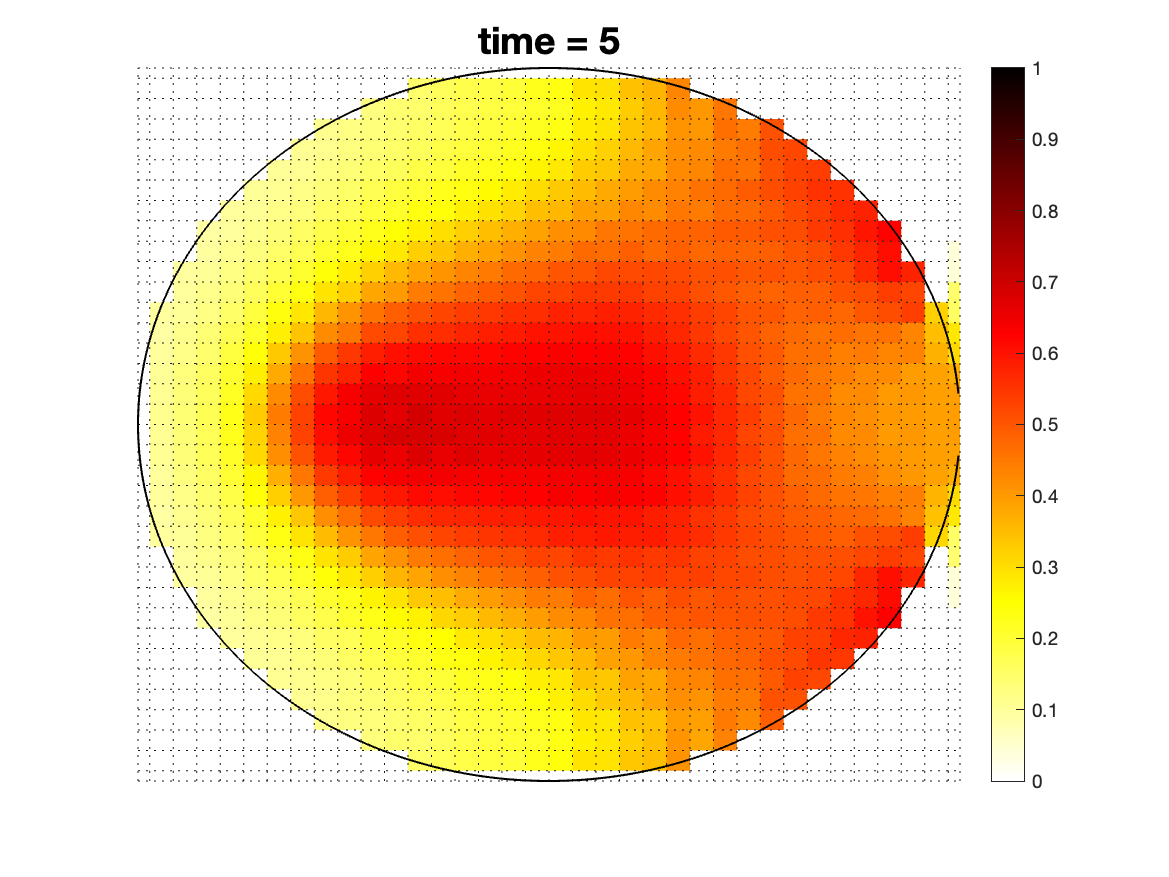}
\end{overpic} 
\begin{overpic}[width=0.32\textwidth]{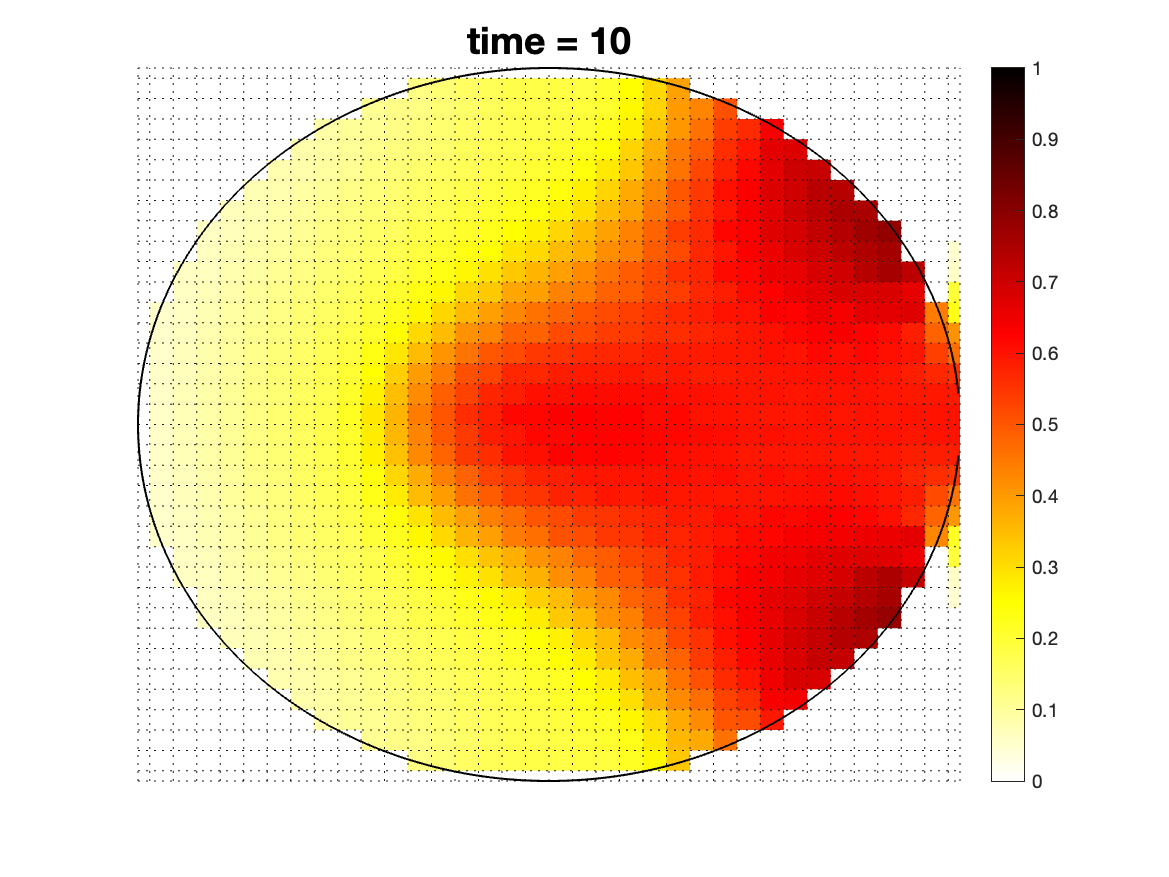}
\put(-15,78){\textcolor{black}{optimized density (from $\varepsilon=0.05$)}}
\end{overpic} 
\includegraphics[width=0.32\textwidth]{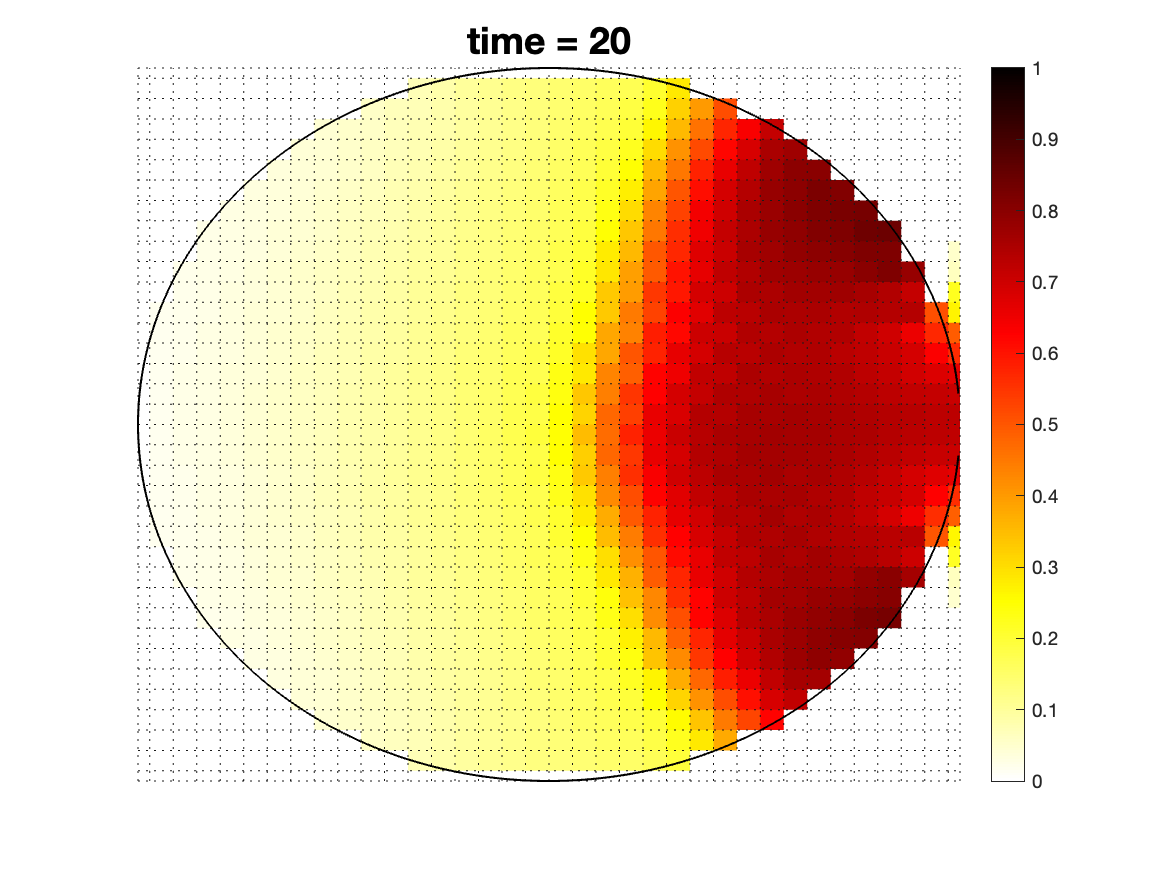}
\\
\vskip .4cm
\begin{overpic}[width=0.32\textwidth]{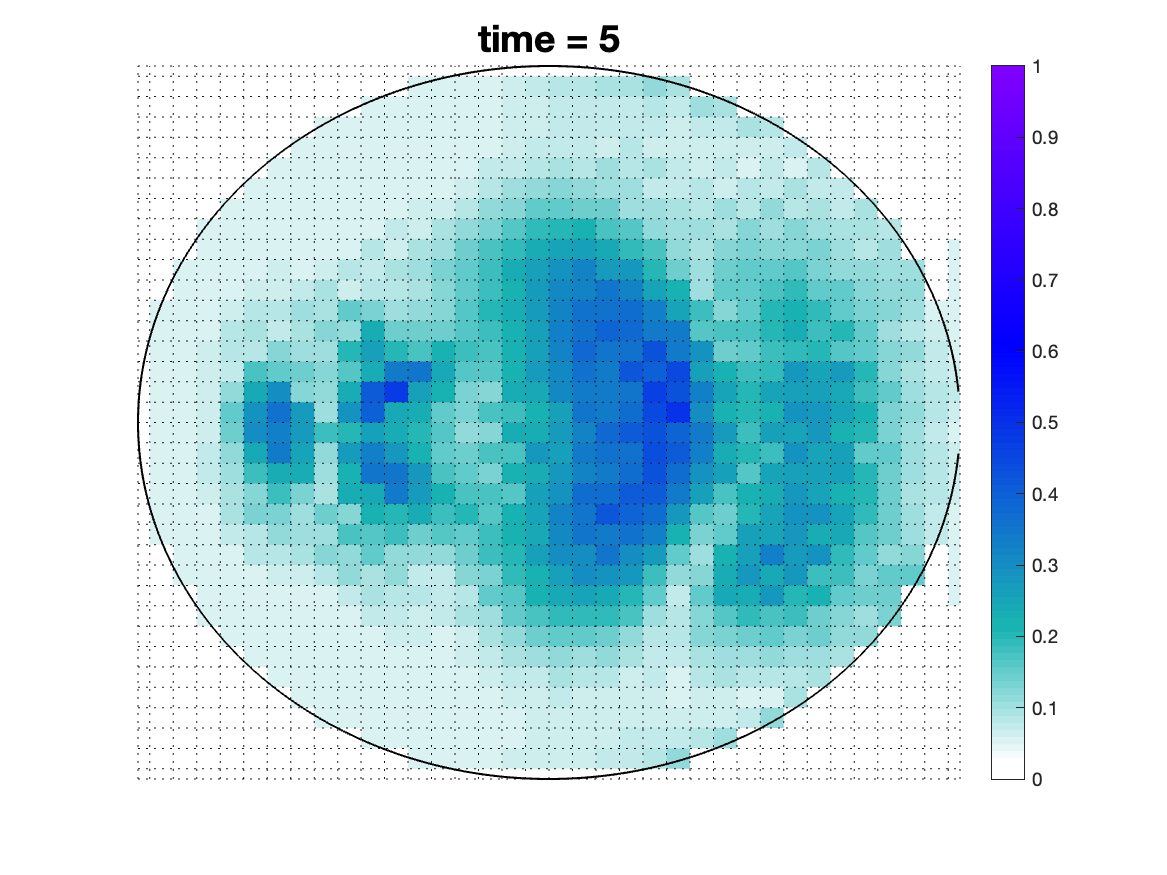}
\end{overpic} 
\begin{overpic}[width=0.32\textwidth]{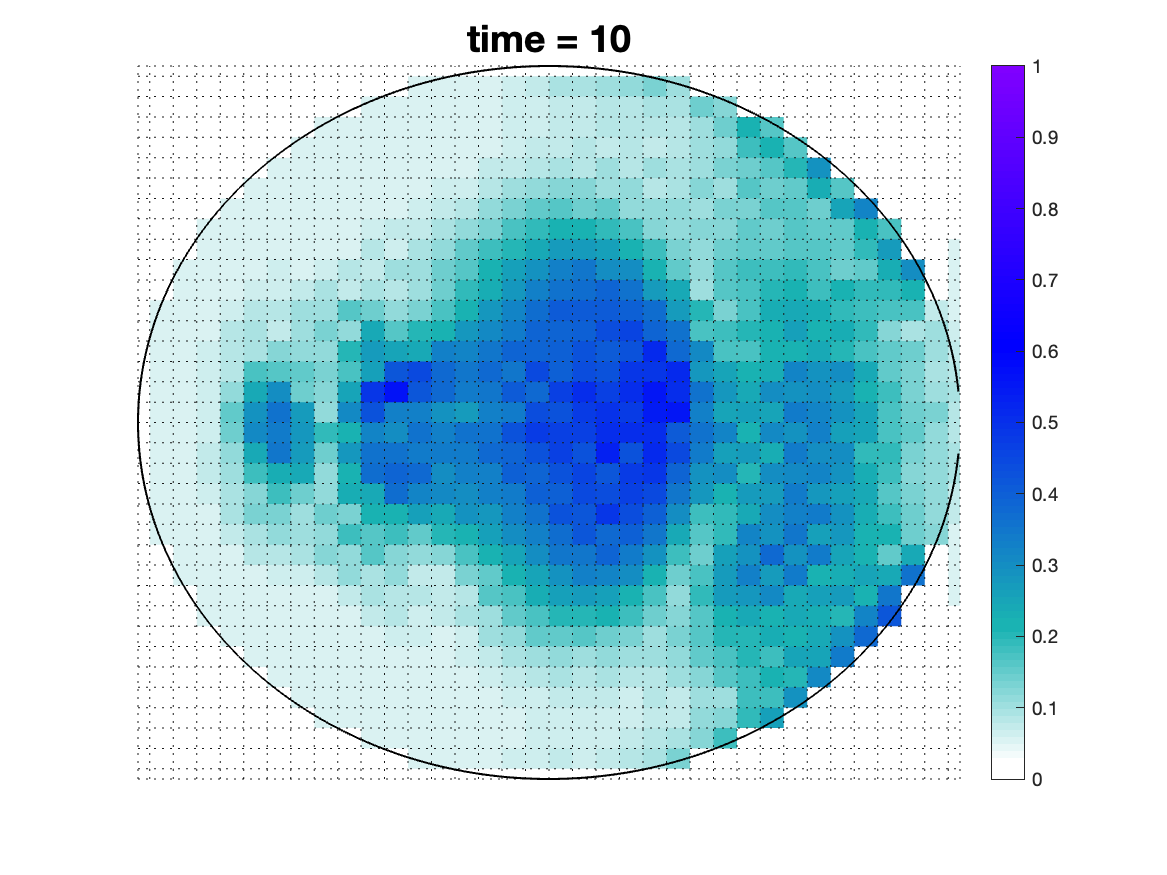}
\put(-20,78){\textcolor{black}{optimized stress level (from $\varepsilon=0.05$)}}
\end{overpic}
\includegraphics[width=0.32\textwidth]{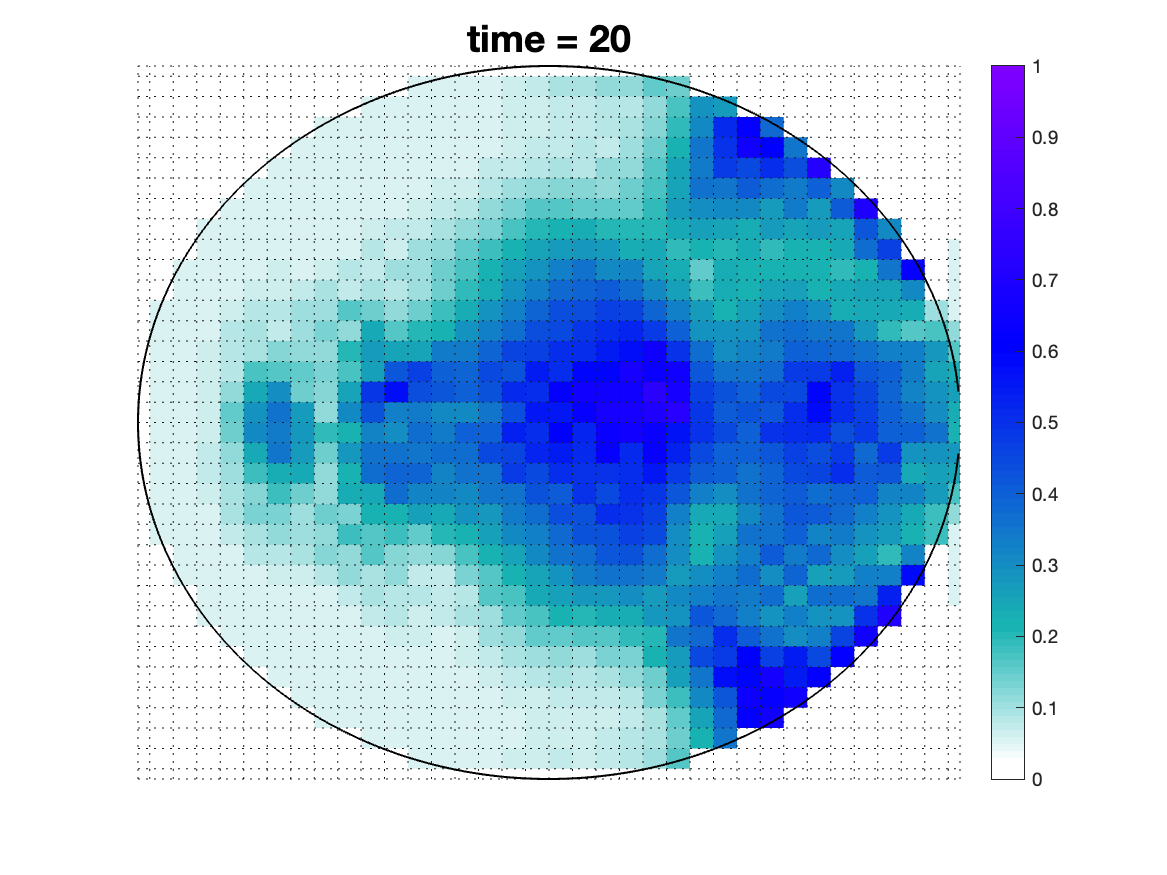}
\caption{Synthetic density data (top), 
optimized density (center), optimized
stress level (bottom) at $t = 5$ s (left), $t = 10$ s (center), and $t = 20$ s (right). Ants are initially placed as
in Fig.~\ref{fig:initialcircle1} (left) and no regularization is used.}
\label{fig:circleprocess}
\end{figure}

For the numerical simulations in Fig.~\ref{fig:circular_2eps}, we set $\Delta x$ = $\Delta y$ = 1 mm and $\Delta t =0.5$ s. For the inverse problem, we will consider the computed results for $\varepsilon=0.95$ as synthetic video data and start the optimization procedure described in Sec.~\ref{sec:num_sol} with $\varepsilon=0.05$. We set $\delta_{k+1} = \delta=50$ in \eqref{eq:update} for all $k$, $tol= 1.0e-05$ in \eqref{eq:stop},
and choose not to use regularization.
The comparison of synthetic data and optimized
density is shown in Fig.~\ref{fig:circleprocess}, together with the corresponding optimized stress level. We see that the optimized 
density does not capture the high density
regions around the center of the chamber
at $t = 5, 10$ s, but does capture the 
arrow-shaped region at $t = 20$ s.
Note that the optimized stress level
does not vary from the initial value 
$\varepsilon = 0.05$ in the regions of the domain where the computed density
matches the synthetic data, i.e., the low-density region in yellow, while it is increased 
where there is a mismatch, i.e., around the center of the chamber. However, while
the procedure is correctly identifying where in the domain the fear level needs to be increased, the optimized value of the fear level is not 
large enough to match well the synthetic
data around the center of the chamber.


\begin{figure}[htb!]
\centering
\begin{overpic}[width=0.47\textwidth,grid=false]{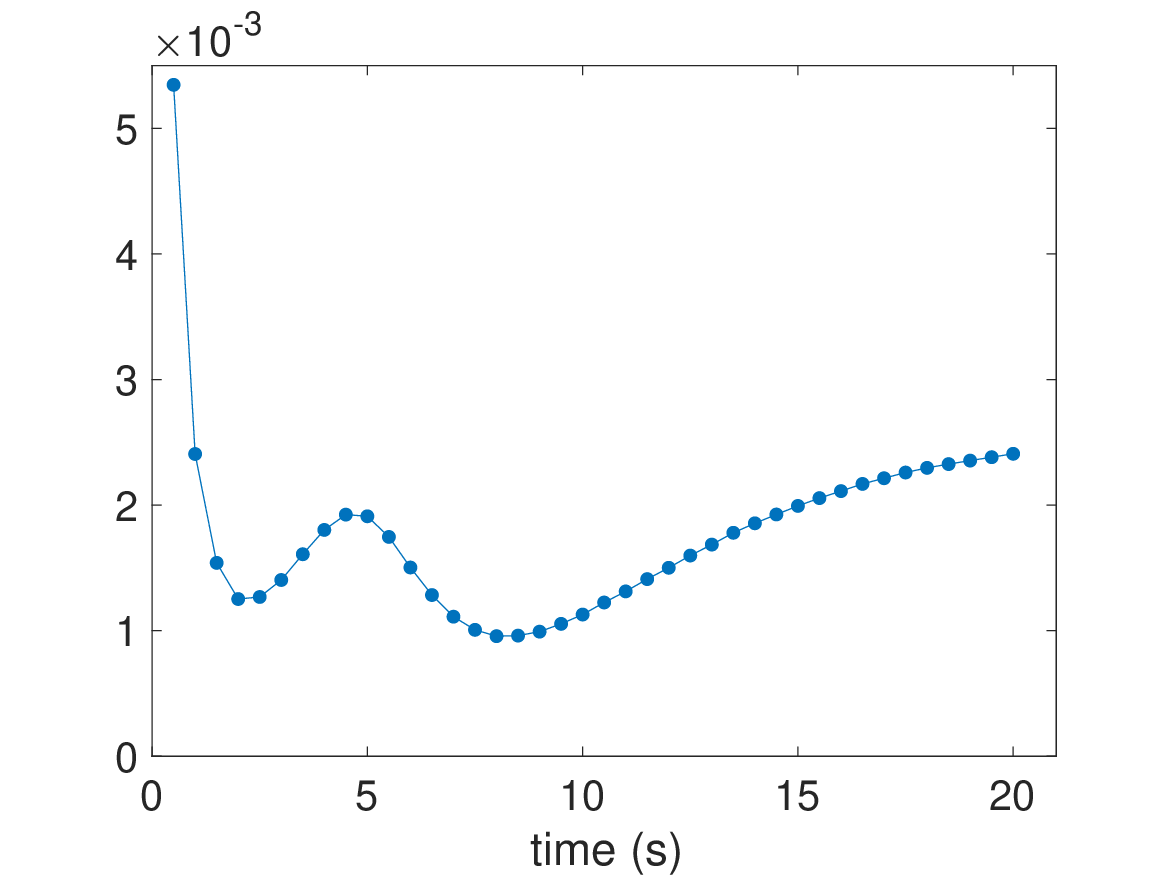}
\end{overpic}
\begin{overpic}[width=0.47\textwidth,grid=false]{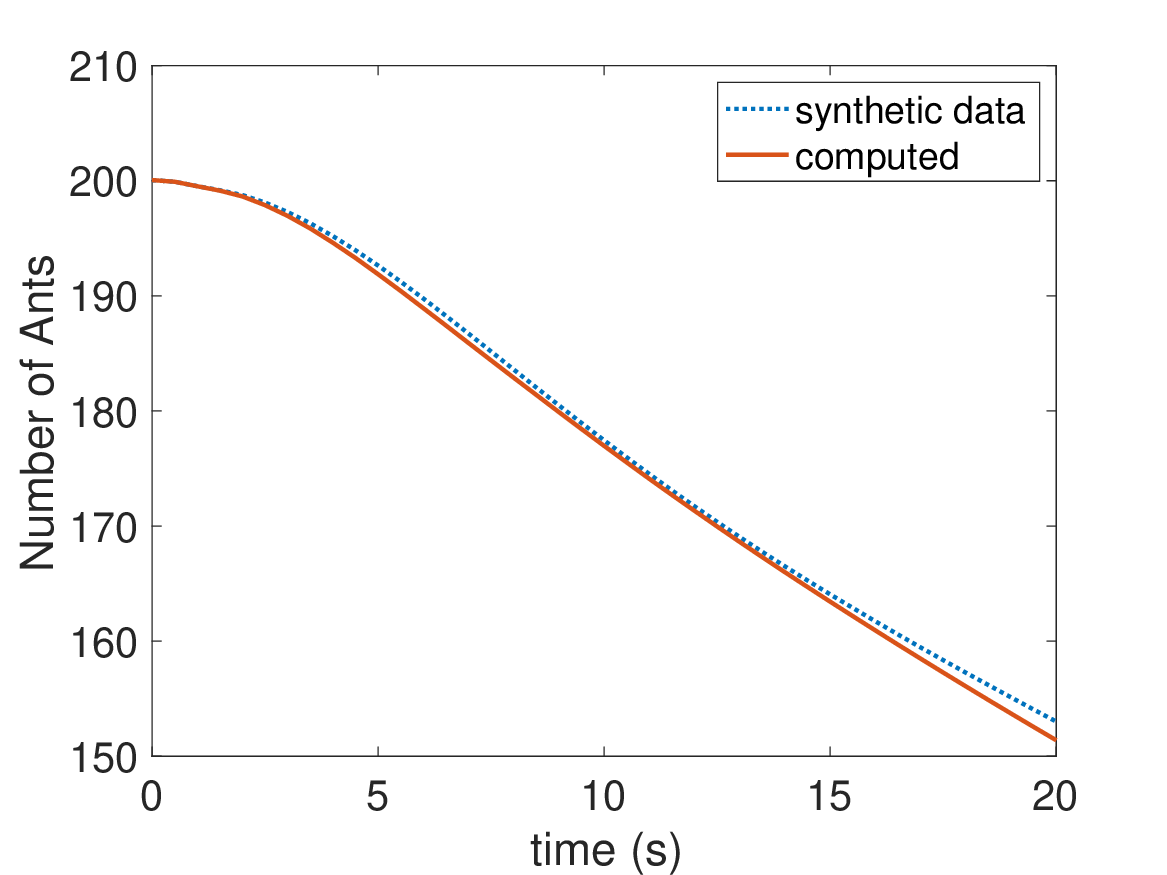}
\end{overpic}
\caption{Left: Functional \eqref{eq:td_J}
over time for the circular chamber without column.
Right: Number of ants inside the chamber over time. 
}
\label{fig:circle_J_ants}
\end{figure}

Fig.~\ref{fig:circle_J_ants} (left) shows the time evolution of functional \eqref{eq:td_J}, which quantifies the mismatch between
synthetic density data and optimized
density in the $L^2$ norm. Such mismatch is between $10^{-3}$ and $2.5 \cdot 10^{-3}$ for most of the first 20 s.

Fig.~\ref{fig:circle_J_ants} (right) 
plots the number of ants inside the chamber over time. In 20 s, almost 50 ants have left the room. This is consistent with
the findings from Ref.~\citenum{SHIWAKOTI20111433}, where
it is reported that
the measured mean escape time (out of 30 repetitions) for the first 50 ants is 21.1 s ($\pm$ 2 s standard deviation). 

\begin{figure}[htb!]
\centering
\begin{overpic}[width=0.32\textwidth]{Circle_Room_Ped_Eps0_95_Ants_circle_time_10.eps}
\end{overpic} 
\begin{overpic}[width=0.32\textwidth]
{Circle_Room_Ped_Eps0_95_Ants_circle_time_20.eps}
\put(-10,78){\textcolor{black}{synthetic video data ($\varepsilon=0.95$)}}
\end{overpic} 
\includegraphics[width=0.32\textwidth]{Circle_Room_Ped_Eps0_95_Ants_circle_time_40.eps}
\\
\vskip .4cm
\begin{overpic}[width=0.32\textwidth]{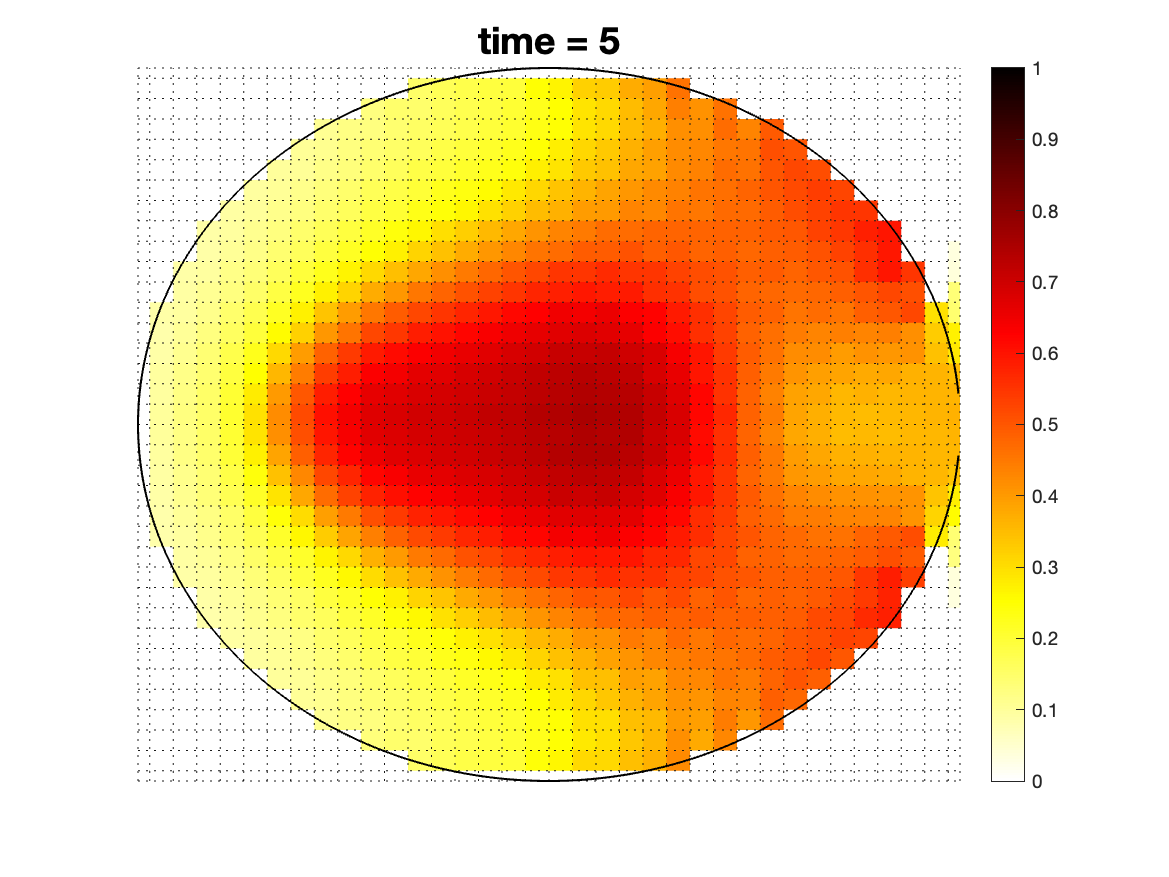}
\end{overpic} 
\begin{overpic}[width=0.32\textwidth]{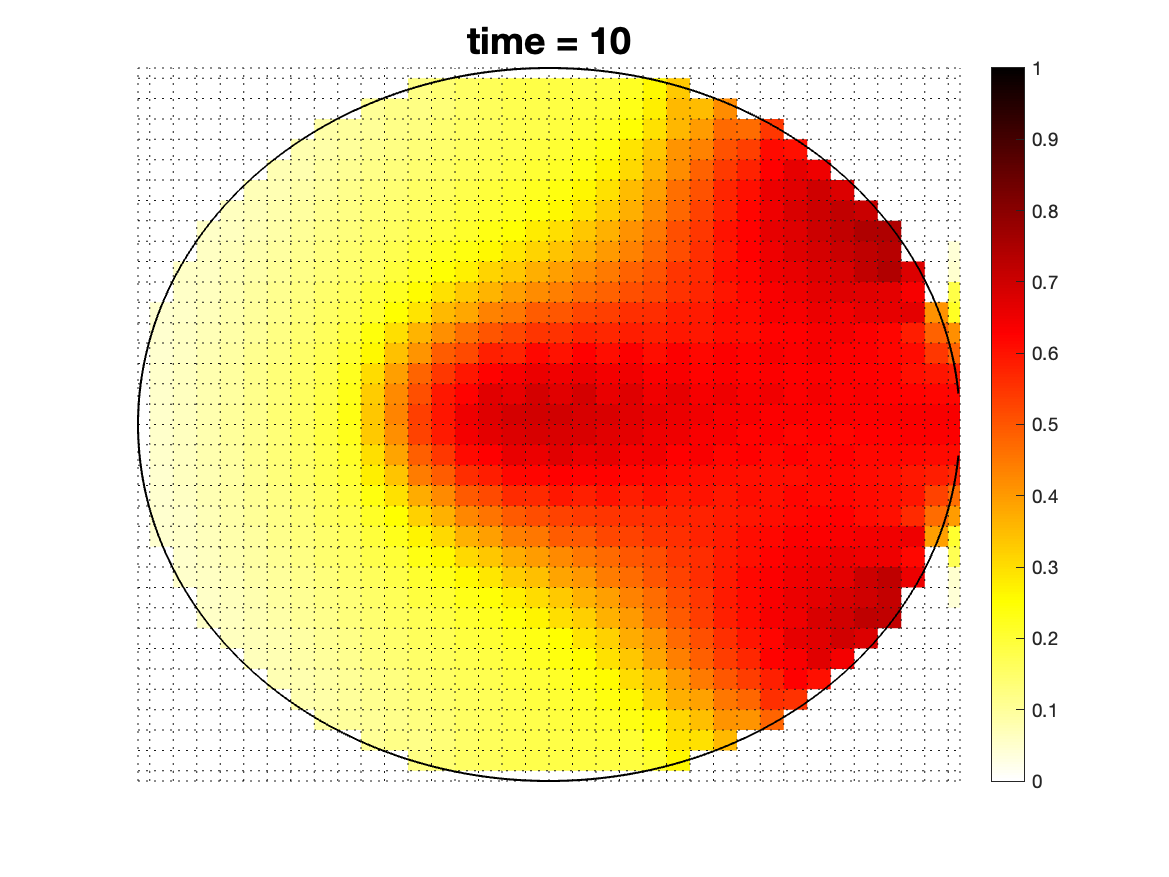}
\put(-15,78){\textcolor{black}{optimized density (from $\varepsilon=0.05$)}}
\end{overpic} 
\includegraphics[width=0.32\textwidth]{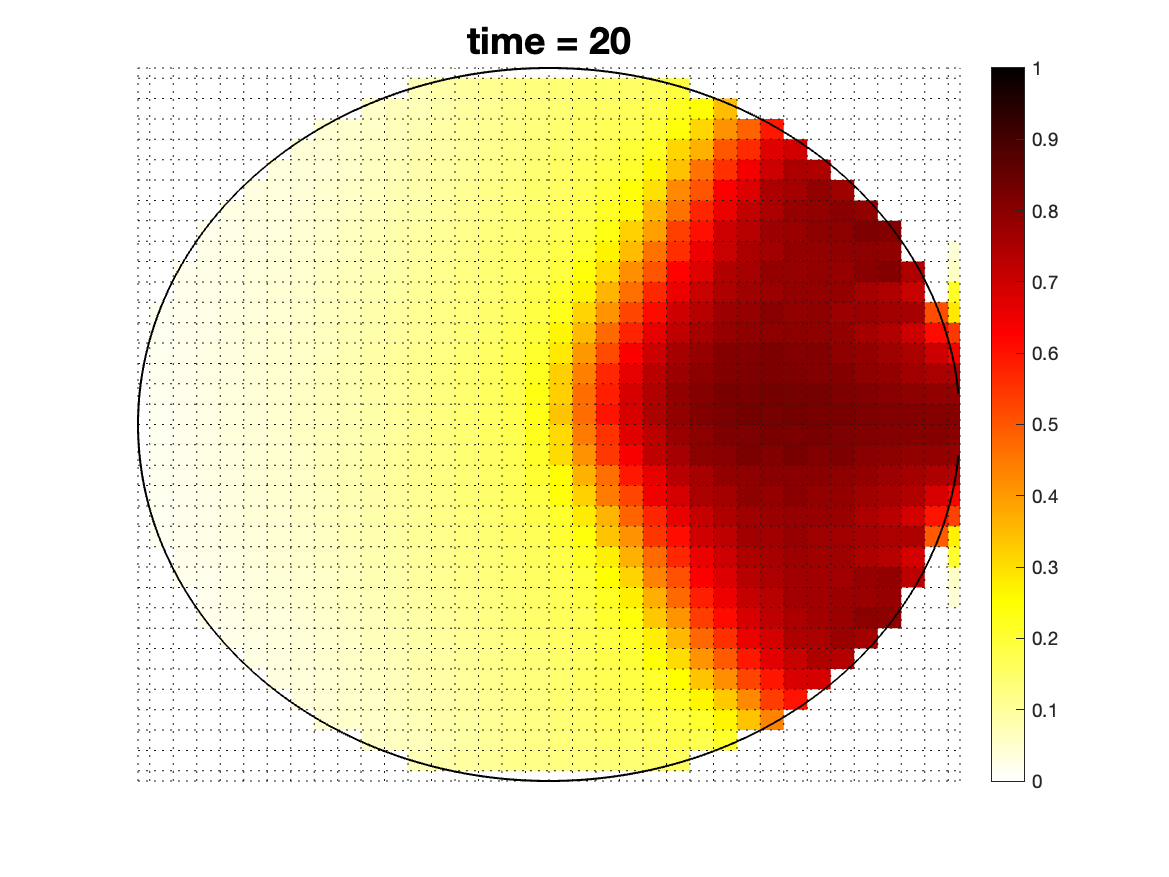}
\\
\vskip .4cm
\begin{overpic}[width=0.32\textwidth]{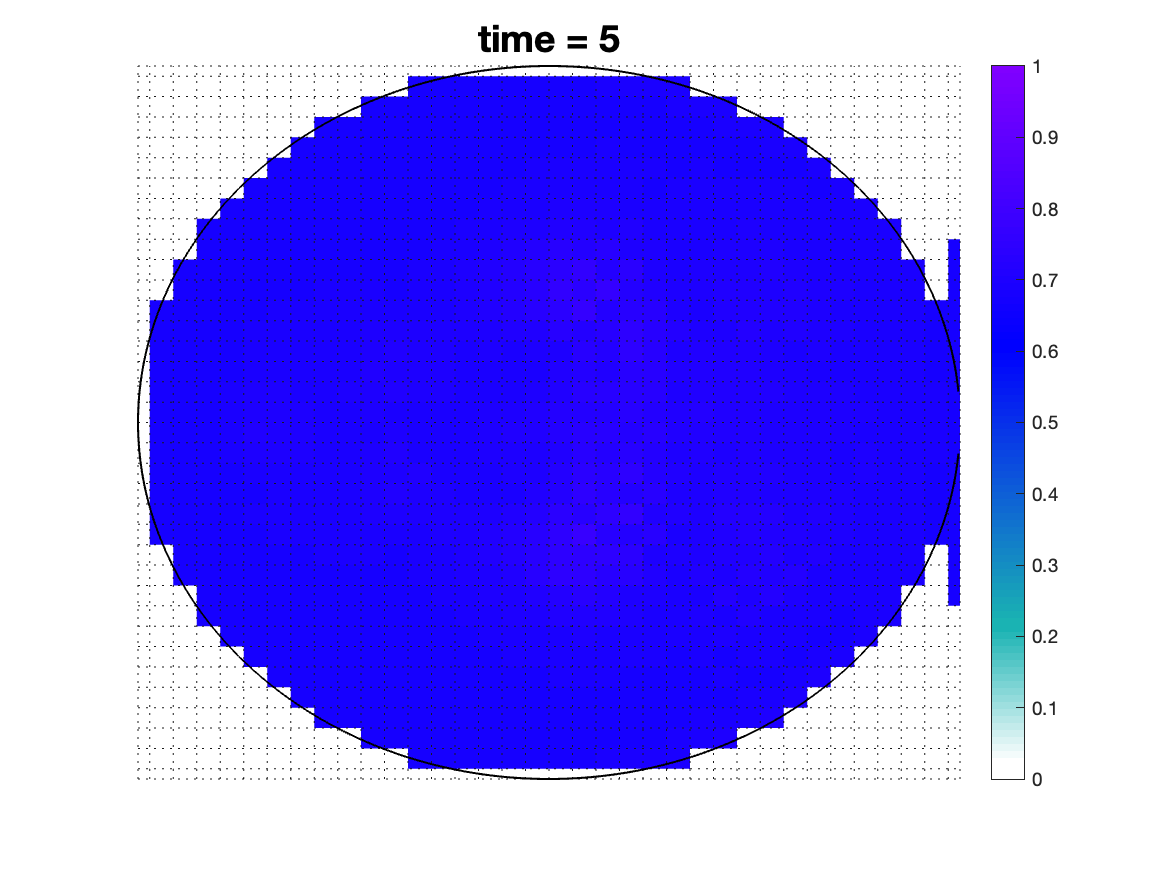}
\end{overpic} 
\begin{overpic}[width=0.32\textwidth]{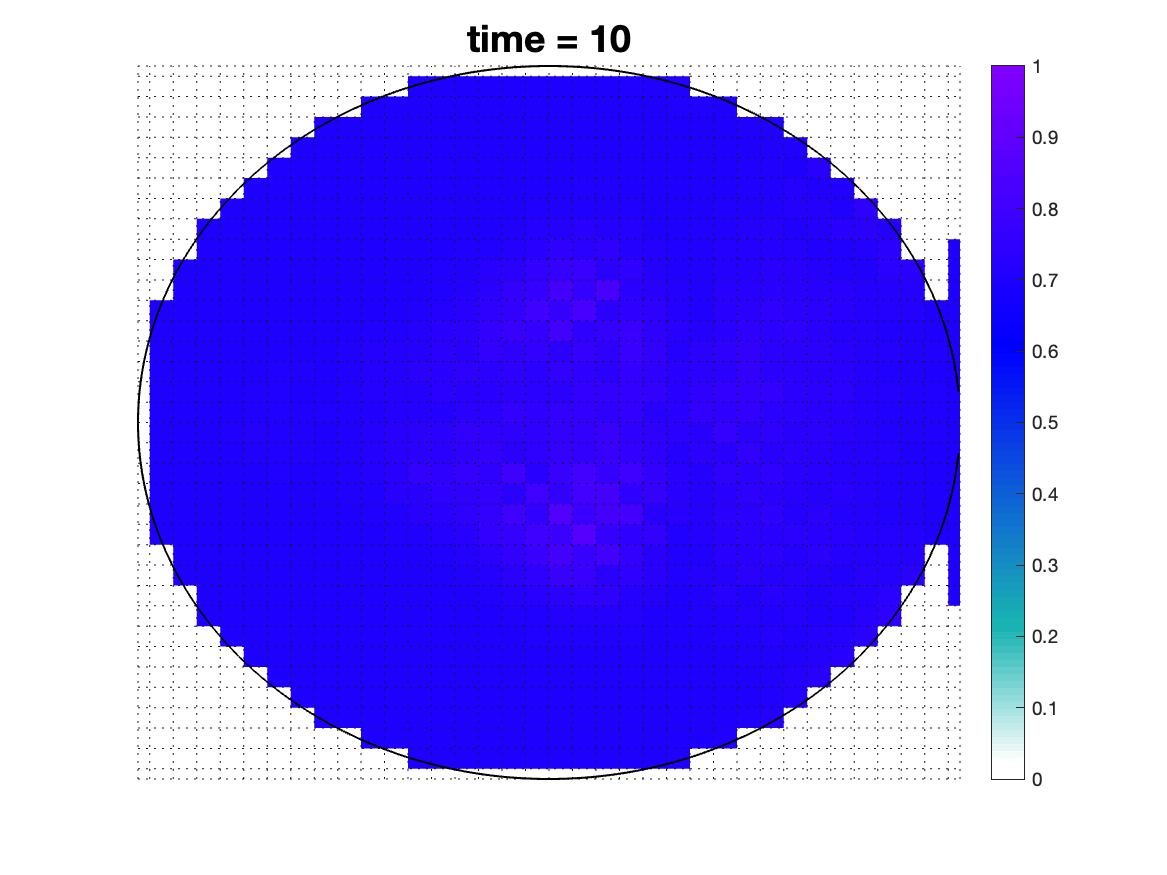}
\put(-20,78){\textcolor{black}{optimized stress level (from $\varepsilon=0.05$)}}
\end{overpic}
\includegraphics[width=0.32\textwidth]{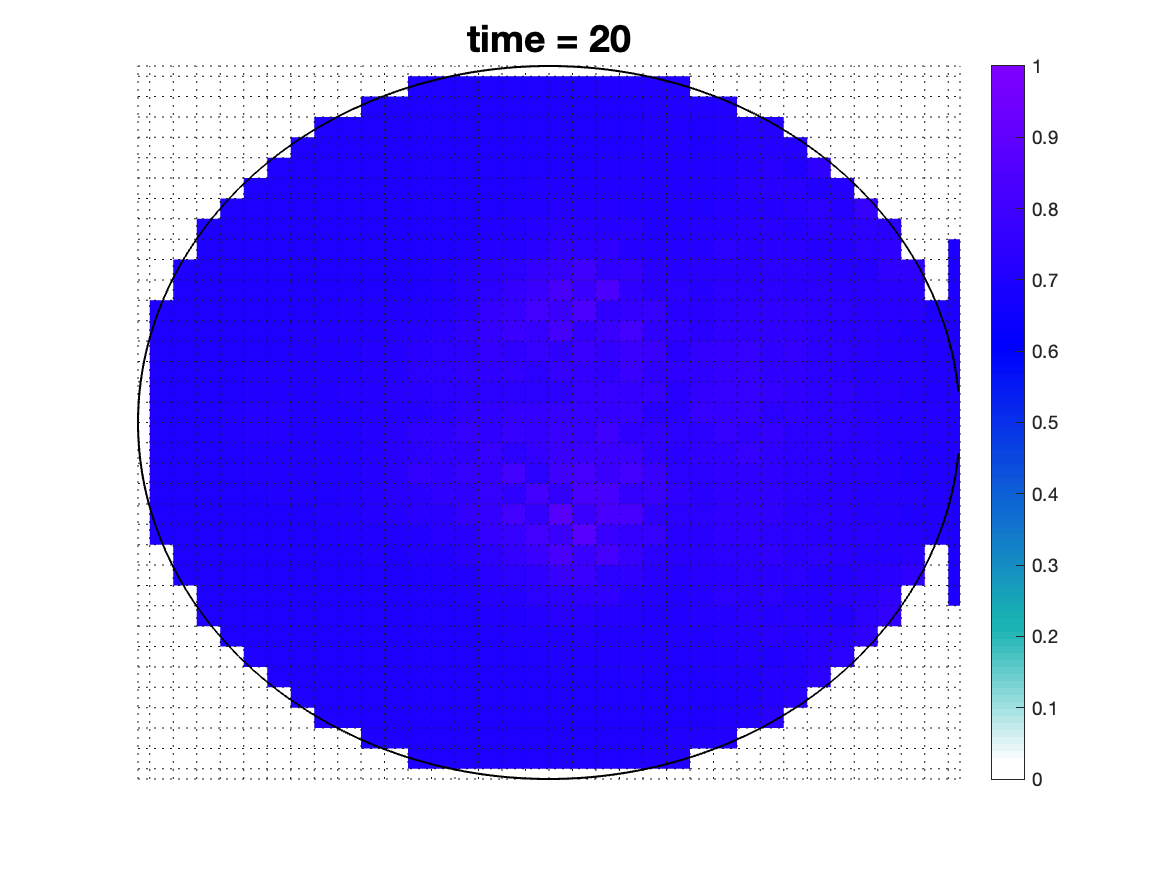}
\caption{Synthetic density data (top), 
optimized density (center), optimized
stress level (bottom) at $t = 5$ s (left), $t = 10$ s (center), and $t = 20$ s (right). Ants are initially placed as
in Fig.~\ref{fig:initialcircle1} (left).
The optimized results were obtained using regularization \eqref{eq:td_J_R} with
$\varepsilon_{ref} = 0.75$.}
\label{fig:circleprocess_reg}
\end{figure}

\begin{figure}[htb!]
\centering
\begin{overpic}[width=0.47\textwidth,grid=false]{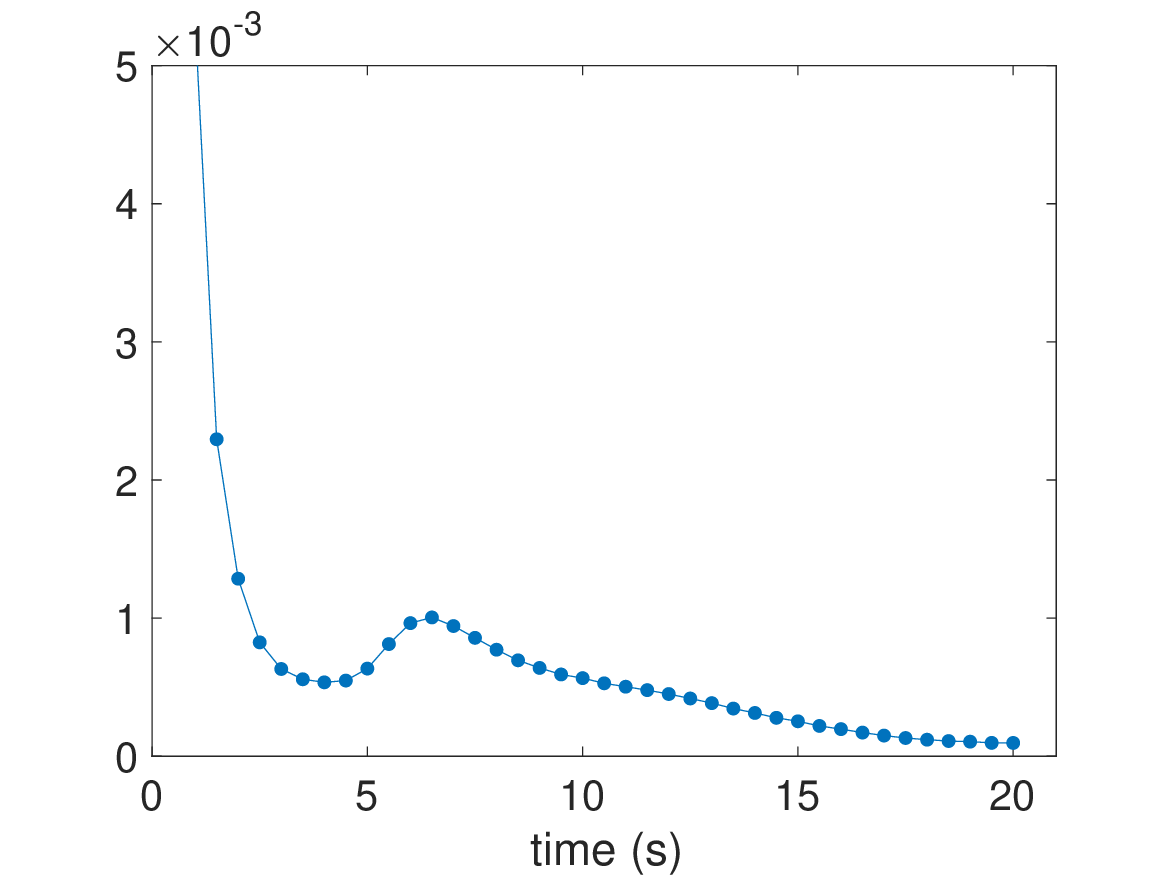}
\end{overpic}
\begin{overpic}[width=0.47\textwidth,grid=false]{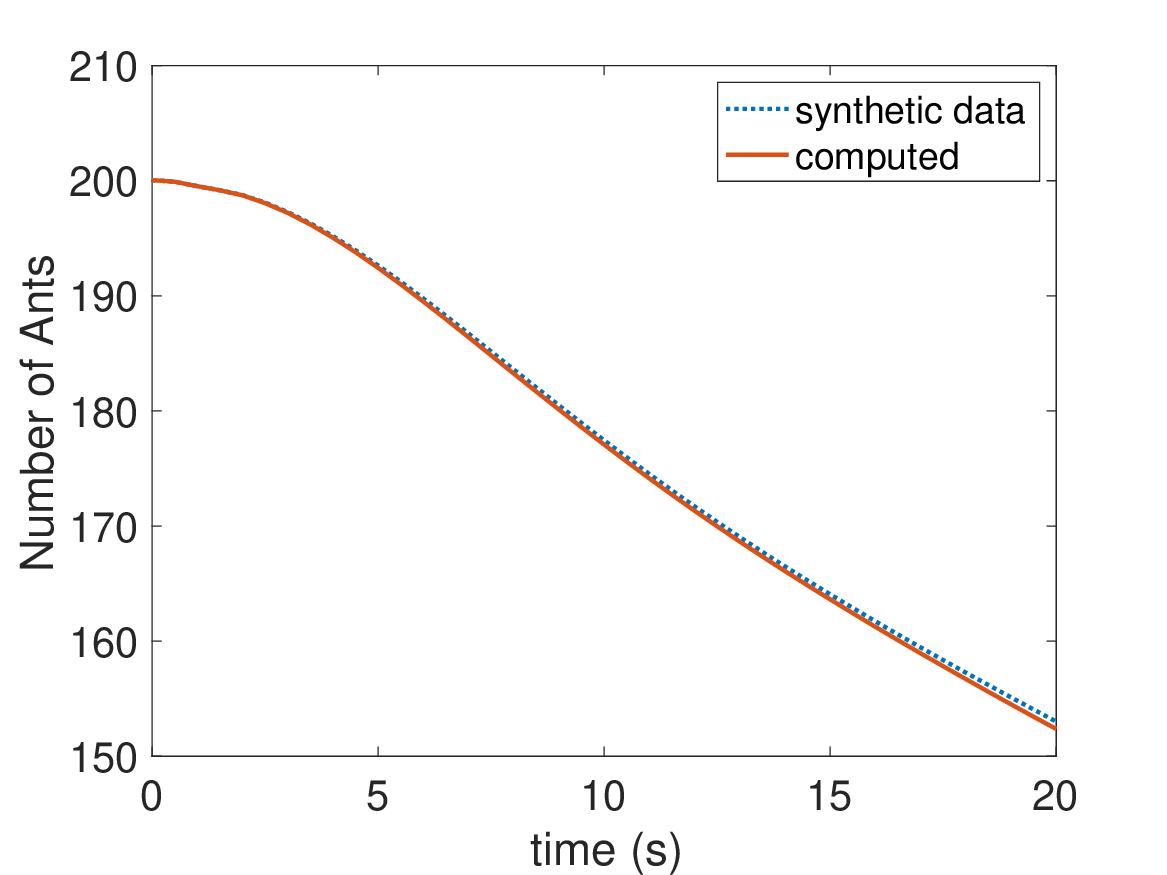}
\end{overpic}
\caption{Left: Functional \eqref{eq:JR_dt}
over time for the circular chamber without column.
Right: Number of ants inside the chamber over time. 
}
\label{fig:circle_J_ants_reg}
\end{figure}

To improve upon the results presented in Fig.~\ref{fig:circleprocess} and \ref{fig:circle_J_ants}, 
we introduce regularization \eqref{eq:td_J_R} with $\varepsilon_{ref}=0.75$ and $\xi=0.1$. 
The values of $\delta$ and $tol$ remain
unchanged.
The comparison with the new optimized
density is shown in 
Fig.~\ref{fig:circleprocess_reg}.
Thanks to the regularization term, the optimized
density now correctly captures the high density
regions around the center of the chamber
at $t = 5, 10$ s, as well as the 
arrow-shaped region at $t = 20$ s.
Note also the improved 
prediction of the number of ants in 
the chamber over time reported in Fig.~\ref{fig:circle_J_ants_reg} (right) with respect to
Fig.~\ref{fig:circle_J_ants} (right). 
By comparing Fig.~\ref{fig:circle_J_ants_reg} (left) with 
Fig.~\ref{fig:circle_J_ants} (left), we see that 
functional \eqref{eq:JR_dt} takes 
smaller values than functional \eqref{eq:td_J}
for $t \geq 5$ s, which reflects the improved
optimized density seen in Fig.~\ref{fig:circleprocess_reg} vs Fig.~\ref{fig:circleprocess}.

\subsection{Circular chamber with partial obstruction}\label{sec:circular_col}

\begin{figure}[htb!]
\center
\includegraphics[width=0.47\textwidth]{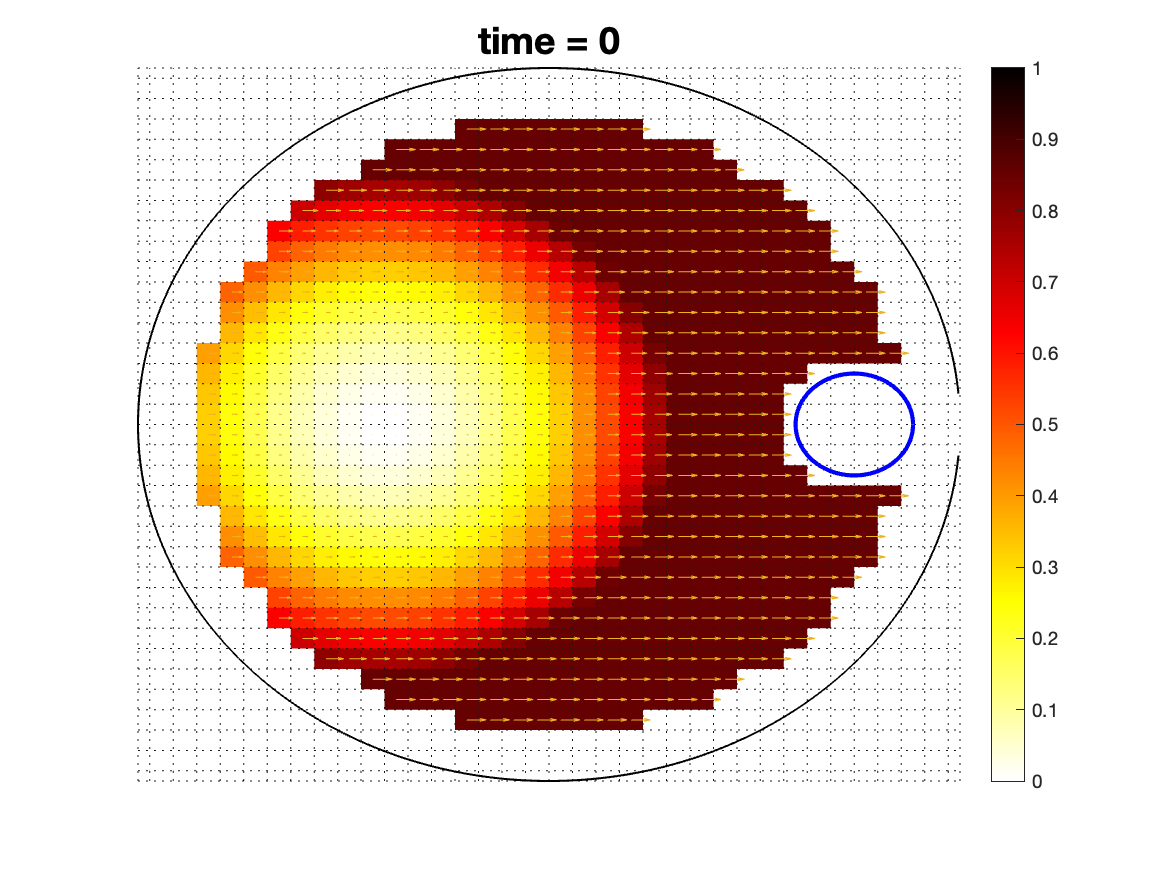}
\begin{overpic}[width=0.33\textwidth,grid=false]{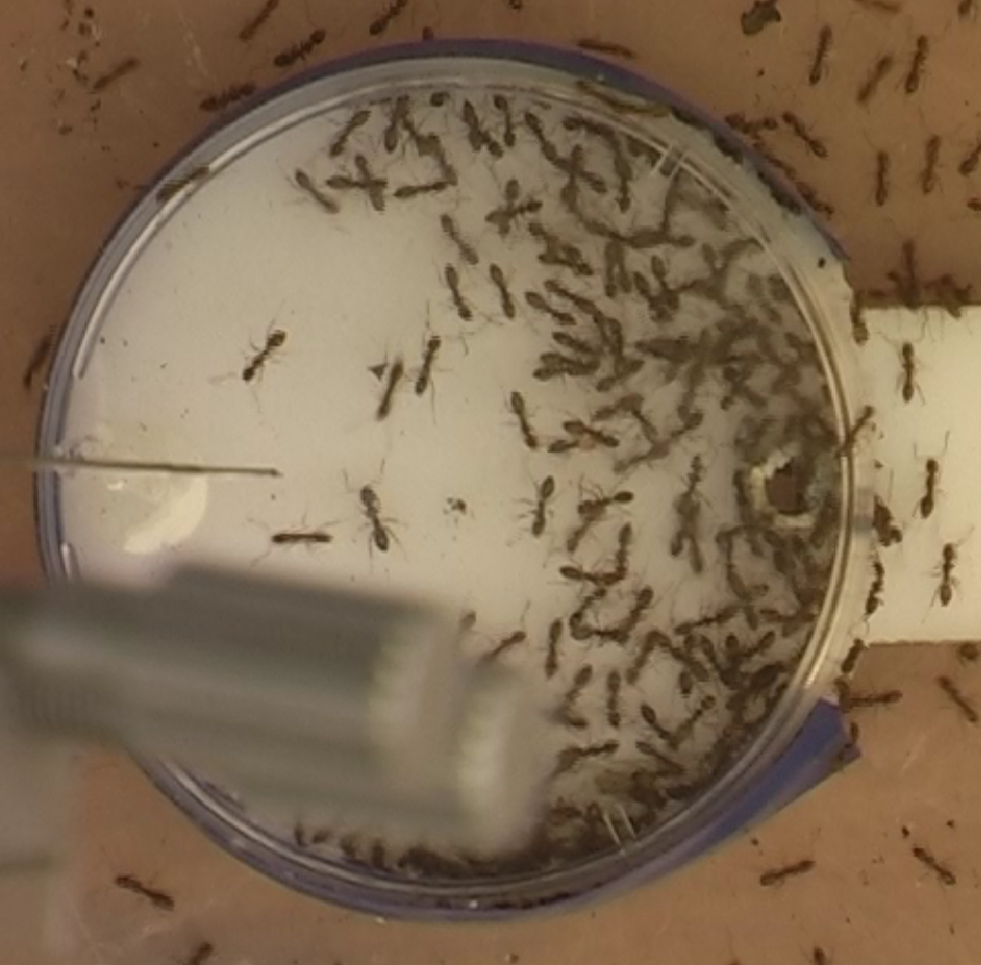}
\end{overpic}
\caption{Left: Computational domain for the circular chamber with column, with initial density and direction for the numerical simulation. Right:
snapshot from repetition 1 out of 30.}
\label{fig:initialcircle2}
\end{figure}

Let us consider the same circular chamber
as in the previous section, but with the addition of a column to obstruct the exit. The column, which has a 5 mm diameter, is located 2 mm to the left of the exit, on the horizontal symmetry lime. 
We place again 200 ants inside the chamber.
The initial density has a crescent-shaped 
high-density region. See Fig.~\ref{fig:initialcircle2} (left).
A similar distribution is observed in one of
the experiments from Ref.~\citenum{SHIWAKOTI20111433}, see Fig.~\ref{fig:initialcircle2} (right).
We assume all the ants initially move with direction $\theta_1$.

\begin{figure}[htb!]
\begin{overpic}[width=0.32\textwidth, grid=false,tics=10]{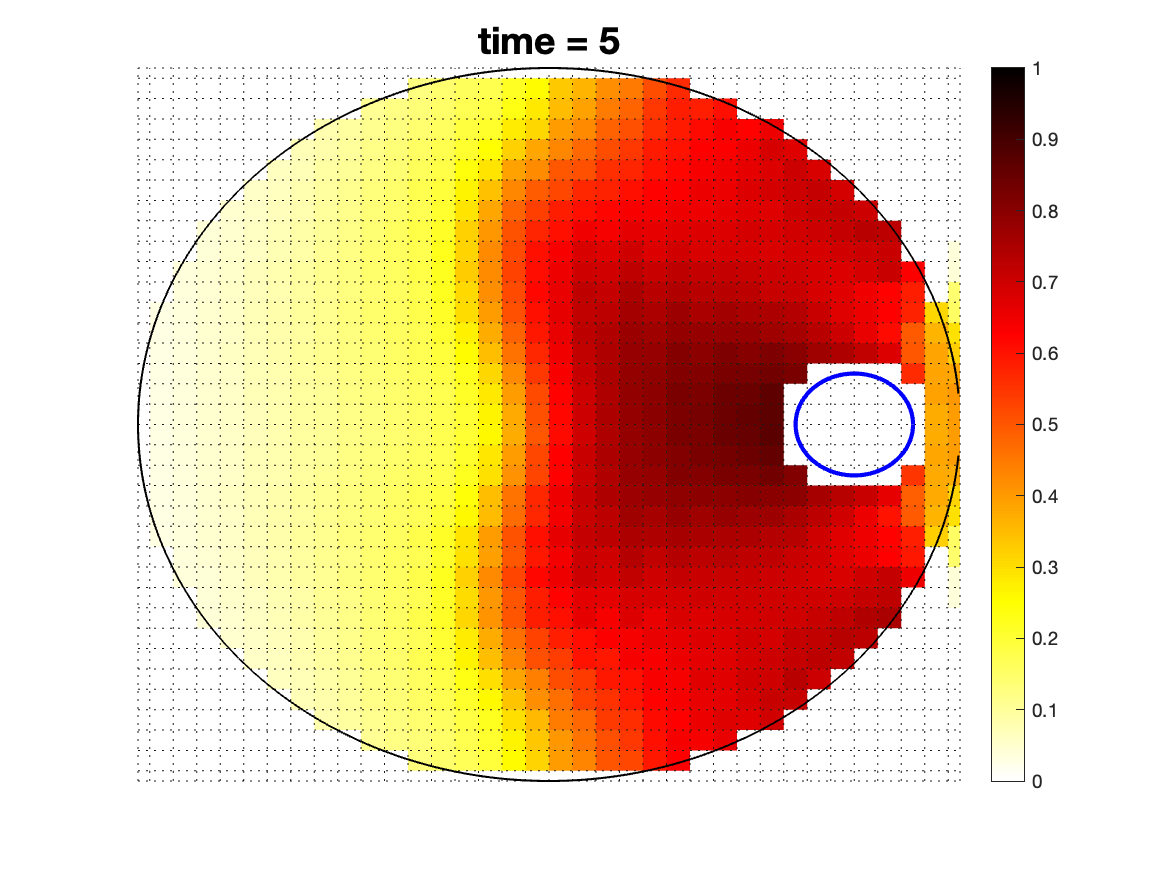}
\end{overpic} 
\begin{overpic}[width=0.32\textwidth, grid=false,tics=10]
{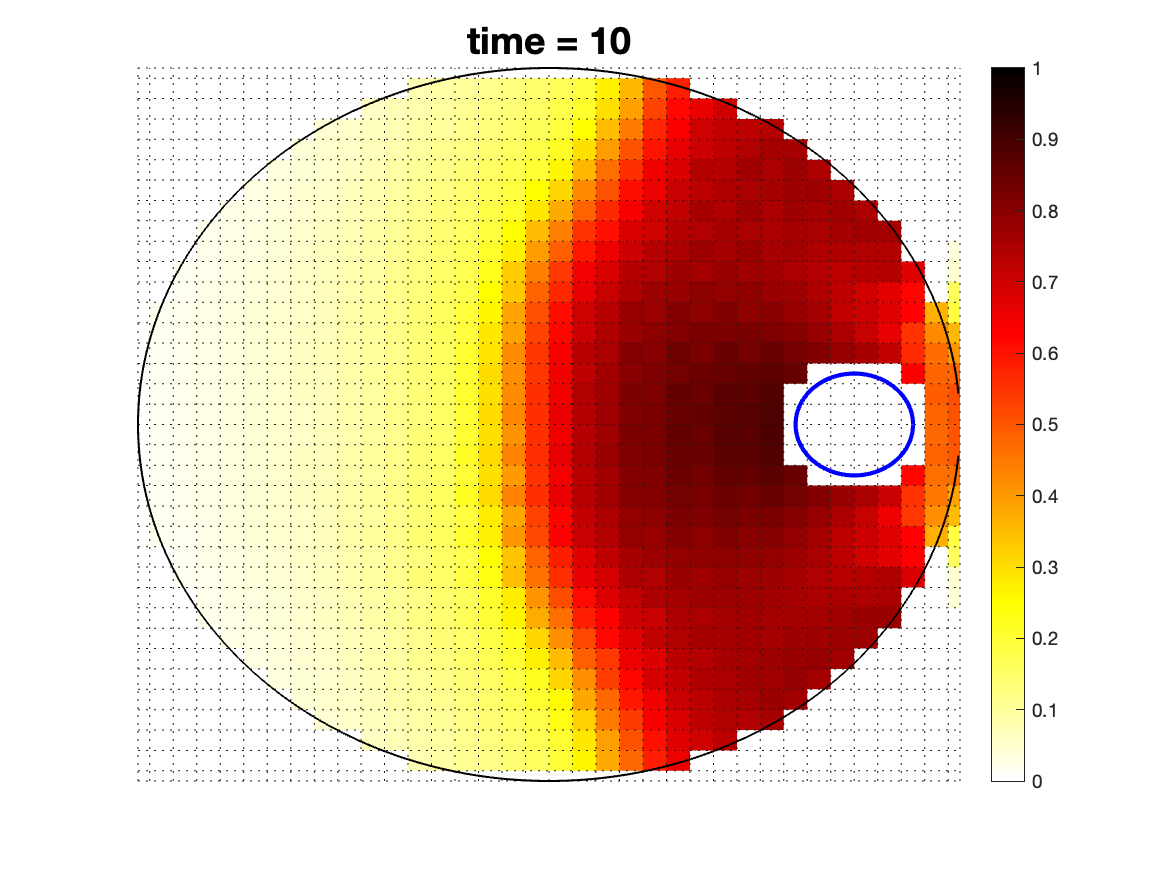}
\put(30,78){\textcolor{black}{$\varepsilon=0.05$}}
\end{overpic}
\includegraphics[width=0.32\textwidth]{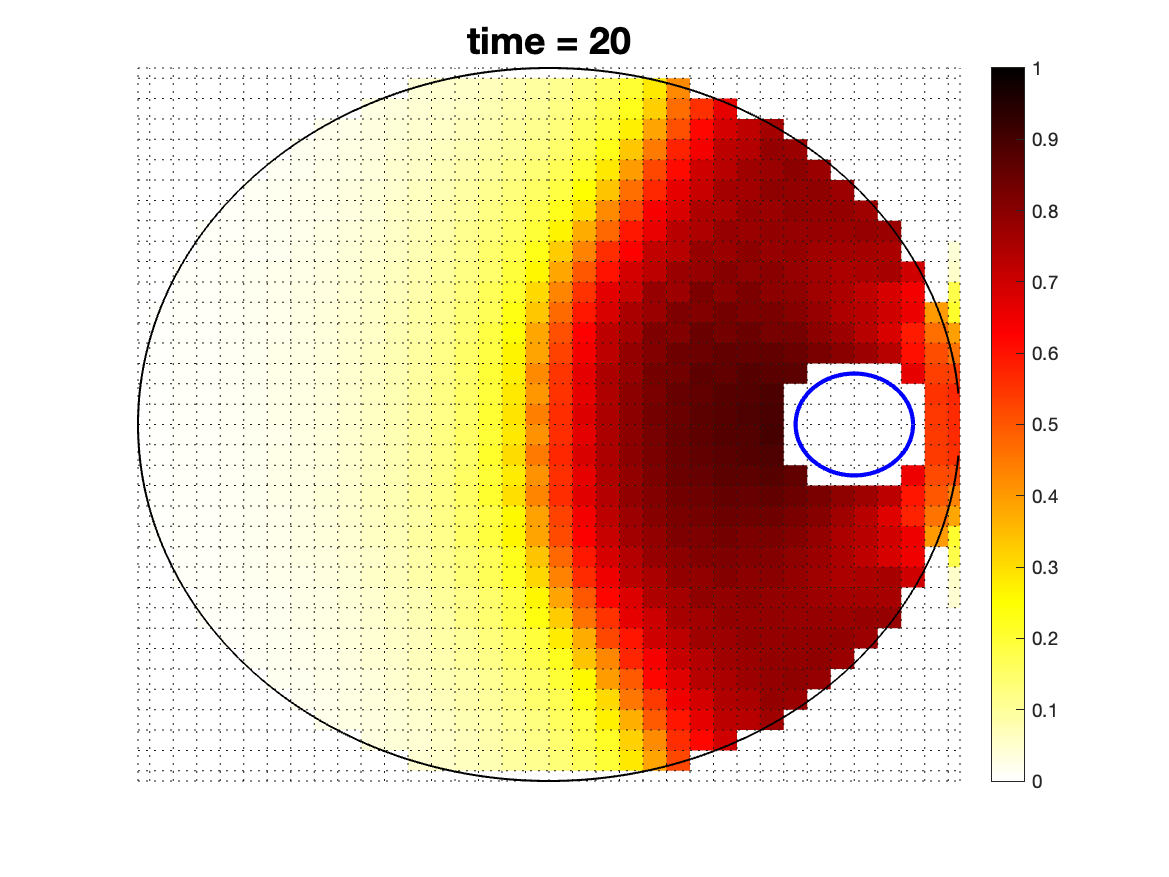} \\
\vskip .2cm
\begin{overpic}[width=0.32\textwidth, grid=false,tics=10]{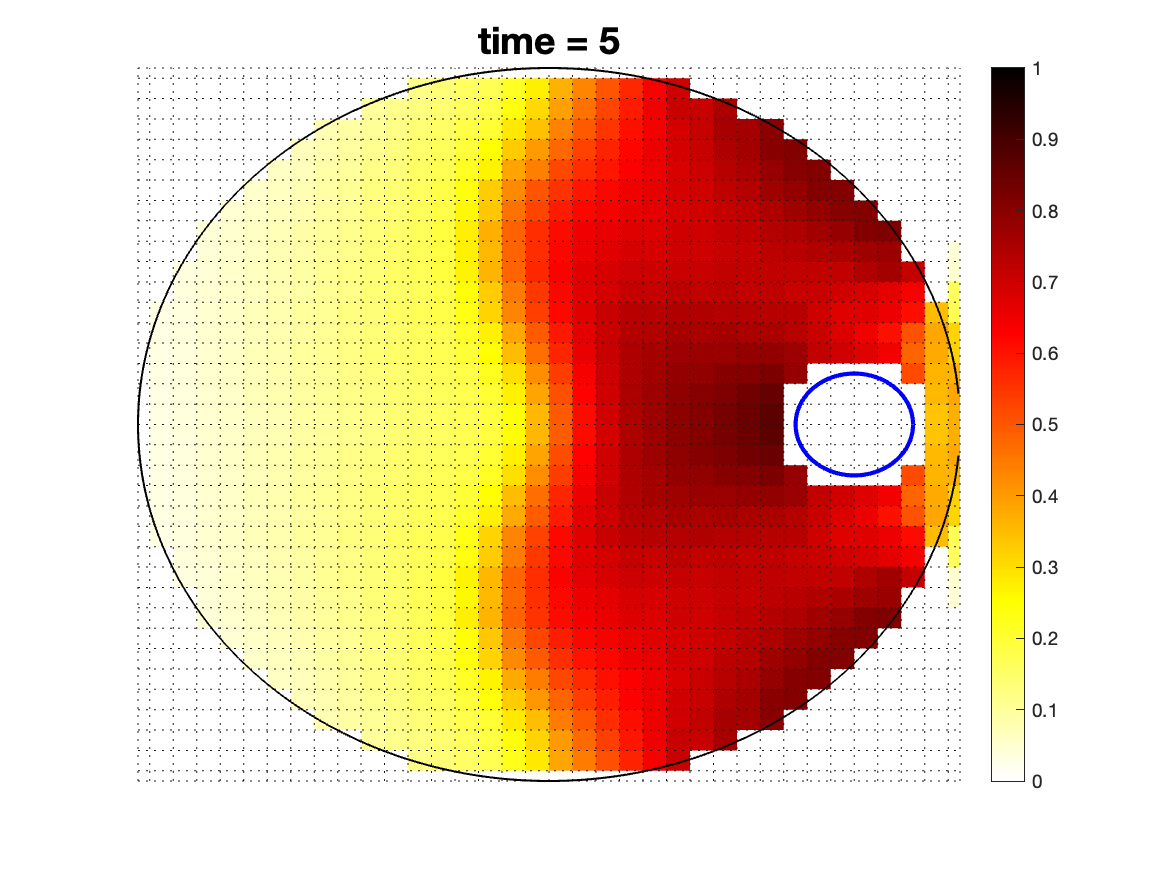}
\end{overpic} 
\begin{overpic}[width=0.32\textwidth, grid=false,tics=10]{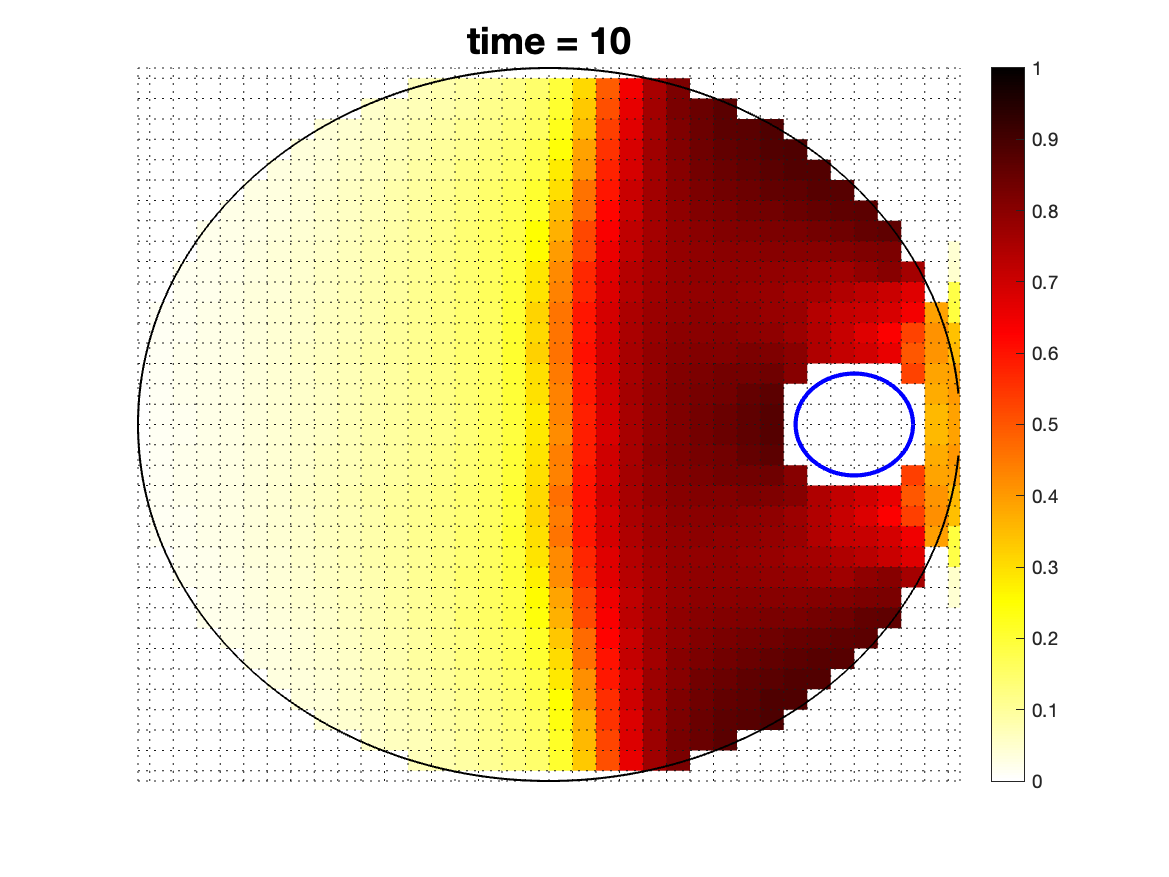}
\put(30,78){\textcolor{black}{$\varepsilon=0.95$}}
\end{overpic} 
\includegraphics[width=0.32\textwidth]
{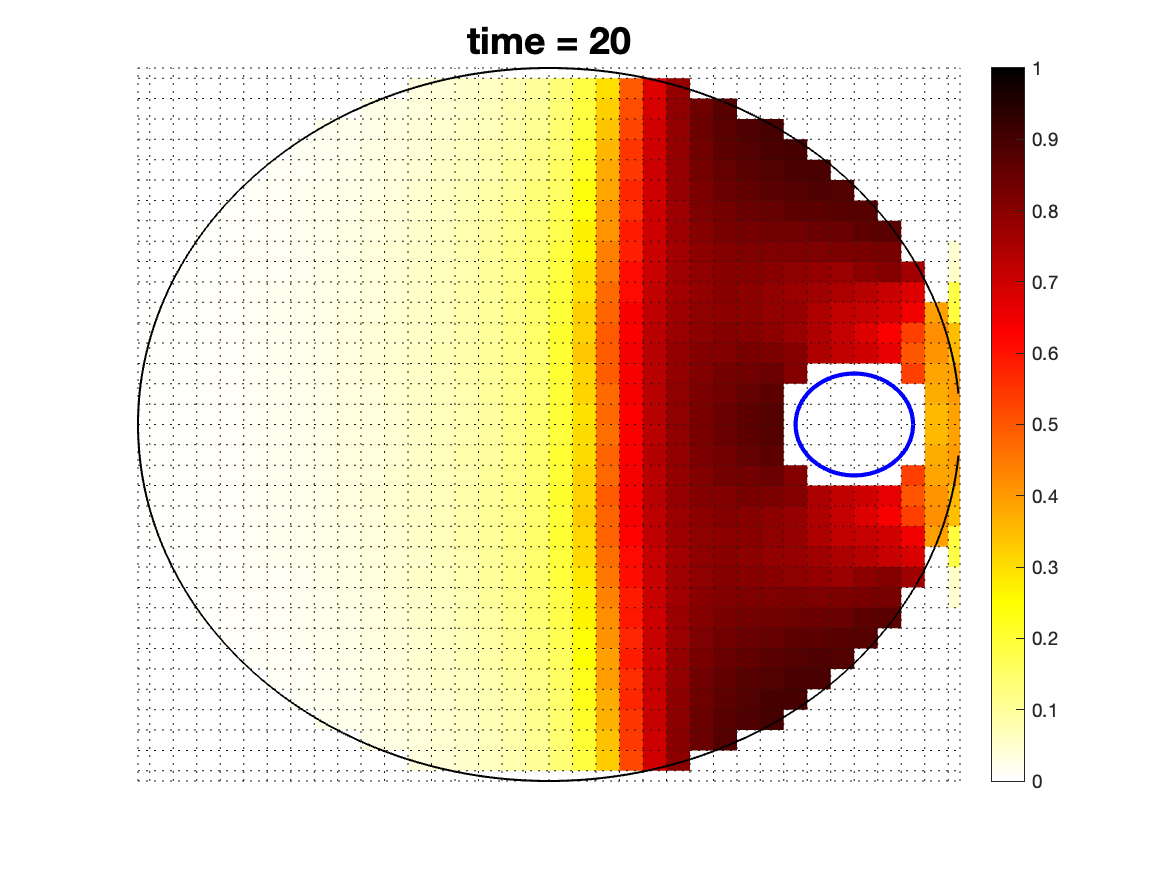}
\caption{Computed density at $t = 5$ s (left), $t = 10$ s (center), and $t = 20$ s (right)
given by the forward problem 
for 200 ants initially placed as
in Fig.~\ref{fig:initialcircle2} (left)
for $\varepsilon=0.05$ (top row) and $\varepsilon=0.95$ (second row).}
\label{fig:circularcol_2eps}
\end{figure}

We consider again conditions of low stress ($\varepsilon = 0.05$)
and high stress ($\varepsilon = 0.95$)
in the initial configuration described above. 
The evacuation dynamics for the two cases
are shown in Fig.~\ref{fig:circularcol_2eps}.
At $t = 5$ s, there is not significant
difference between the computed densities
for $\varepsilon = 0.05$ and $\varepsilon = 0.95$, 
but as time passes the high-density region becomes smaller
and more uniform in space
when the stress level is high. This is because 
the effect of the high stress is to make the active particles follow others, i.e., stick together instead of spreading out, as explained in Sec.~\ref{sec:forward}.

\begin{figure}[htb!]
\begin{overpic}[width=0.32\textwidth, grid=false,tics=10]{Circle_Room_Column_Ped_Eps0_95_Ants_circle_time_10.eps}
\end{overpic} 
\begin{overpic}[width=0.32\textwidth, grid=false,tics=10]{Circle_Room_Column_Ped_Eps0_95_Ants_circle_time_20.eps}
 \put(-10,78){\textcolor{black}{synthetic video data ($\varepsilon=0.95$)}}
\end{overpic} 
\includegraphics[width=0.32\textwidth]{Circle_Room_Column_Ped_Eps0_95_Ants_circle_time_40.eps} \\
\vskip .2cm
\begin{overpic}[width=0.32\textwidth]{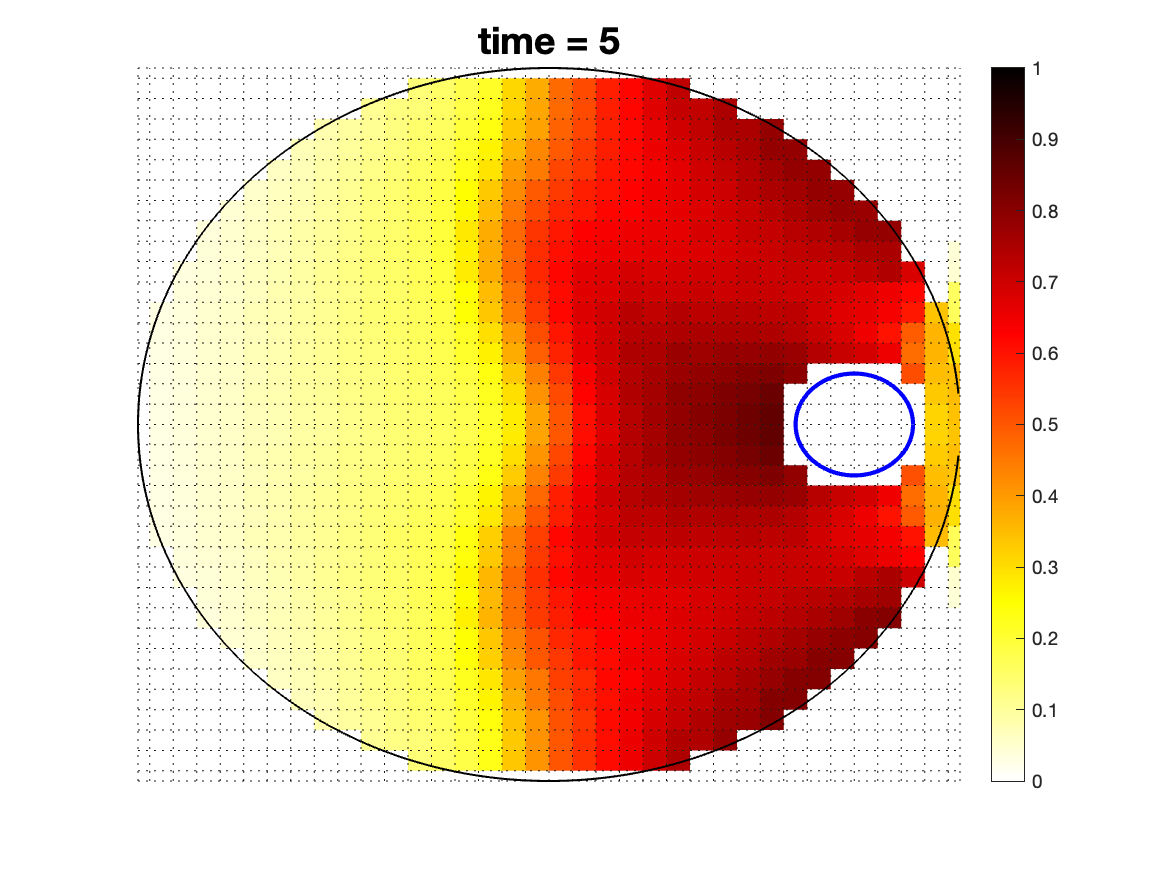}
\end{overpic} 
\begin{overpic}[width=0.32\textwidth]{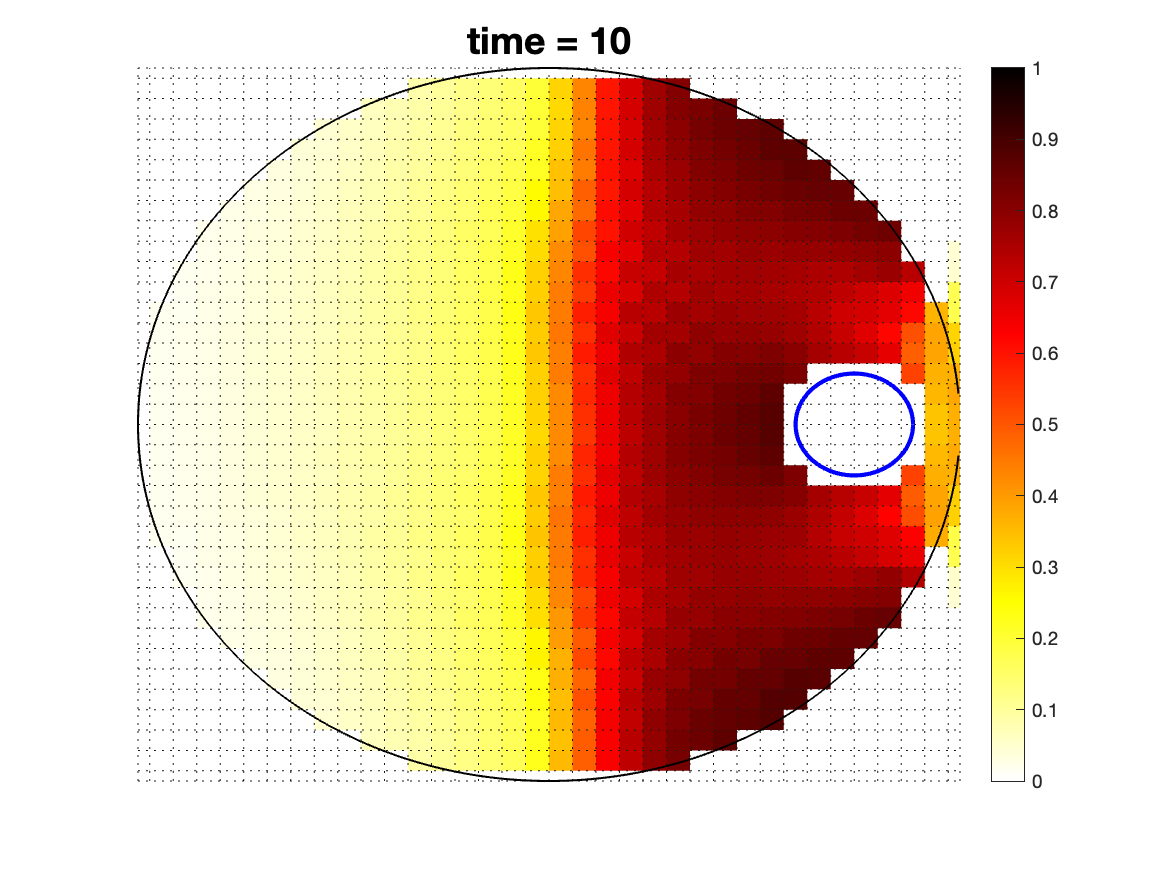}
\put(-15,78){\textcolor{black}{optimized density (from $\varepsilon=0.05$)}}
\end{overpic}
\includegraphics[width=0.32\textwidth]{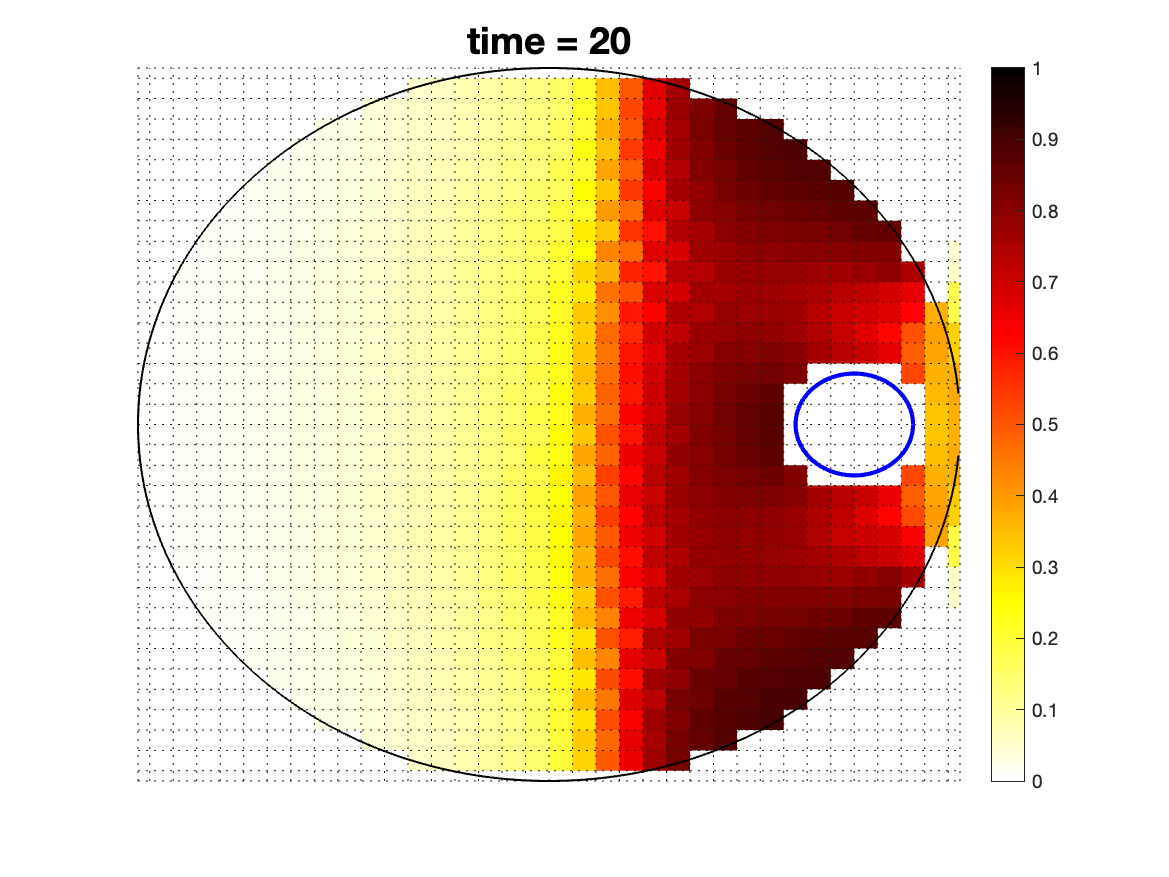}
\\
\vskip .2cm
\begin{overpic}[width=0.32\textwidth]{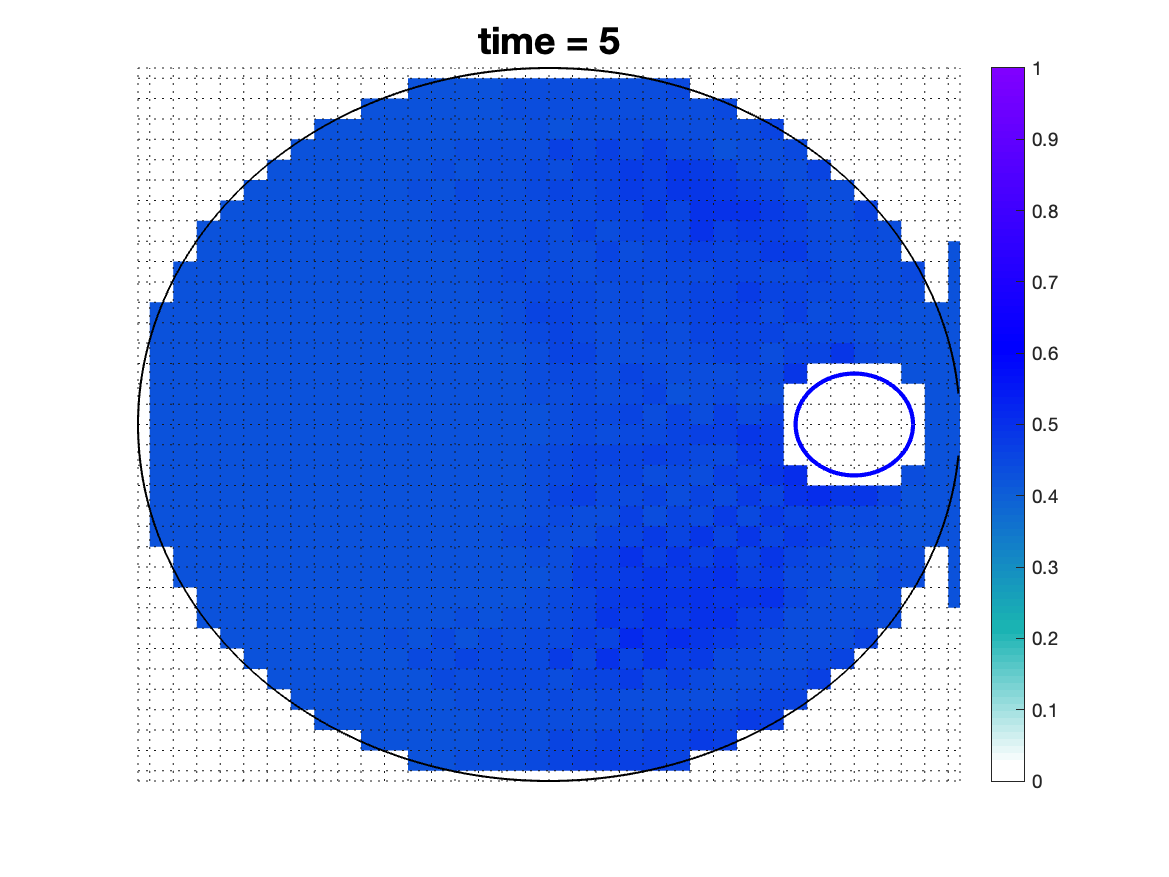}
\end{overpic} 
\begin{overpic}[width=0.32\textwidth]{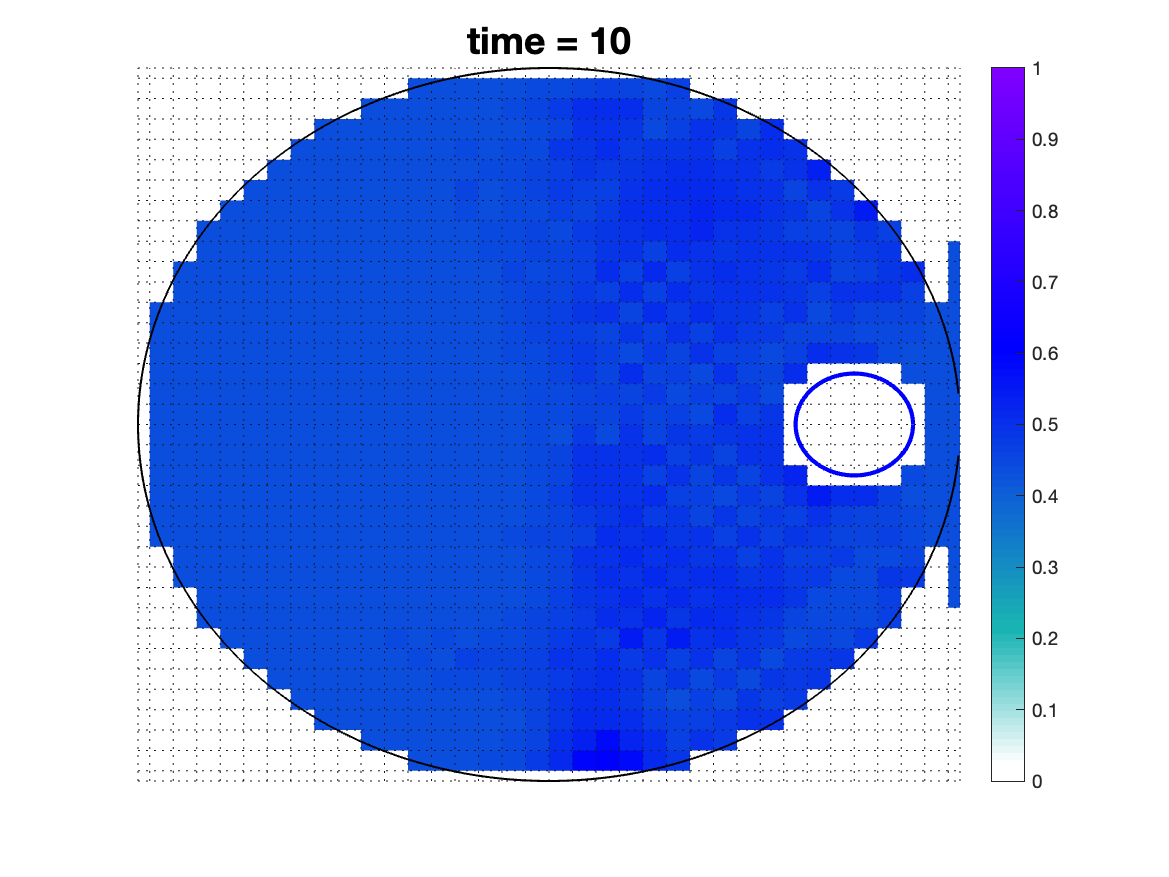}
\put(-20,78){\textcolor{black}{optimized stress level (from $\varepsilon=0.05$)}}
\end{overpic} 
\includegraphics[width=0.32\textwidth]{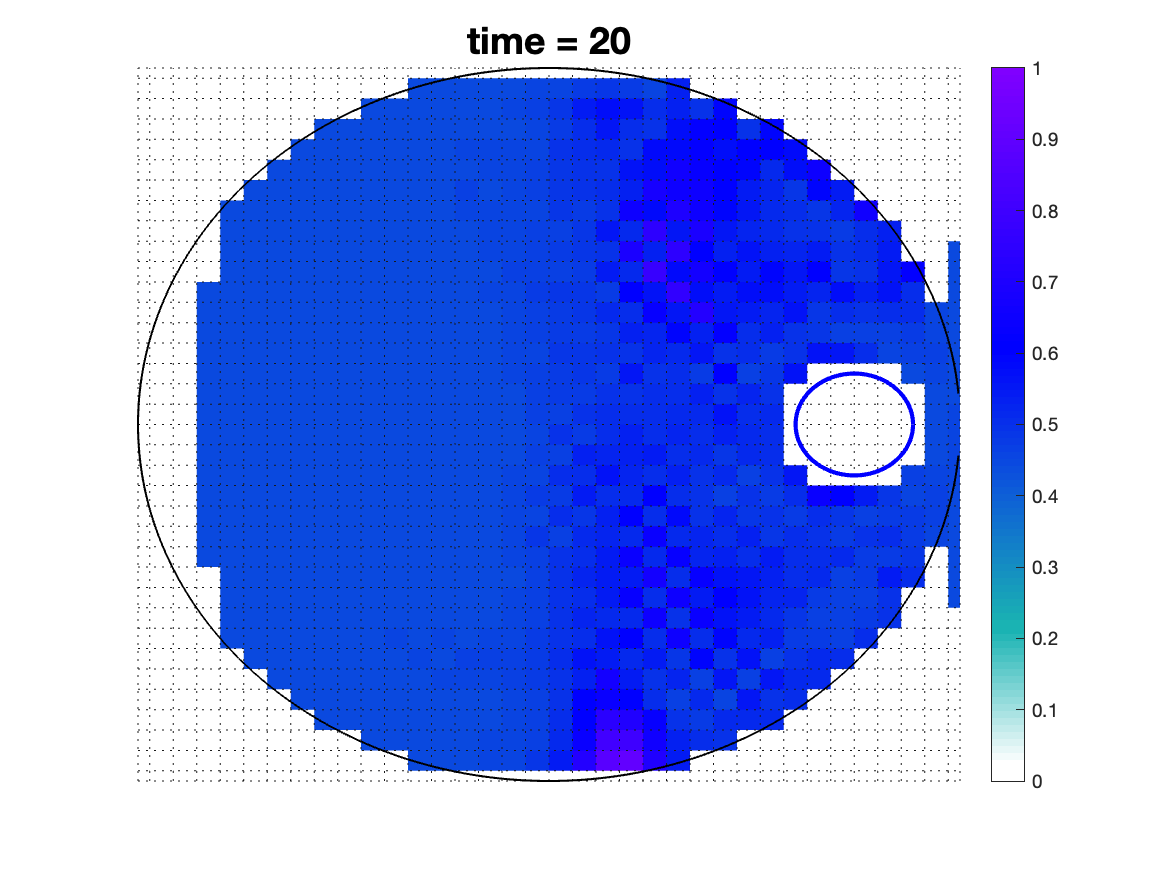}
\caption{Synthetic density data (top), 
optimized density (center), optimized
stress level (bottom) at $t = 5$ s (left), $t = 10$ s (center), and $t = 20$ s (right). Ants are initially placed as
in Fig.~\ref{fig:initialcircle1} (left).
The optimized results were obtained using regularization \eqref{eq:td_J_R} with
$\varepsilon_{ref} = 0.5$.
}
\label{fig:circleprocesswcolumn}
\end{figure}

We use same mesh size and time step as in the previous section, i.e., $\Delta x$ = $\Delta y$ = 1 mm and $\Delta t =0.5$ s. We treat as synthetic video data
the computed results for $\varepsilon=0.95$  and start the optimization procedure with $\varepsilon=0.05$. 
The values of $\delta$ and $tol$ are also the same
as in Sec.~\ref{sec:circular}. However, for this second
test case, we need to regularize the minimization problem
to obtain results that match well the synthetic data.
In particular, we choose regularization \eqref{eq:td_J_R} with $\varepsilon_{ref}=0.5$ and $\xi=0.1$. 
The comparison of synthetic density data and optimized
density is shown in Fig.~\ref{fig:circleprocesswcolumn}, together with the corresponding optimized stress level.
We observe that the synthetic data and the optimized
density match well for the three times shown in  
Fig.~\ref{fig:circleprocesswcolumn}.
As for the optimized stress level, it remains
close to the selected reference value ($\varepsilon_{ref}=0.5$) over the low-density region,
while it is slightly increased in the high-density region. 
This slight increase is sufficient for the optimized
simulation to match the synthetic data.

\begin{figure}[htb!]
\centering
\begin{overpic}[width=0.47\textwidth,grid=false]{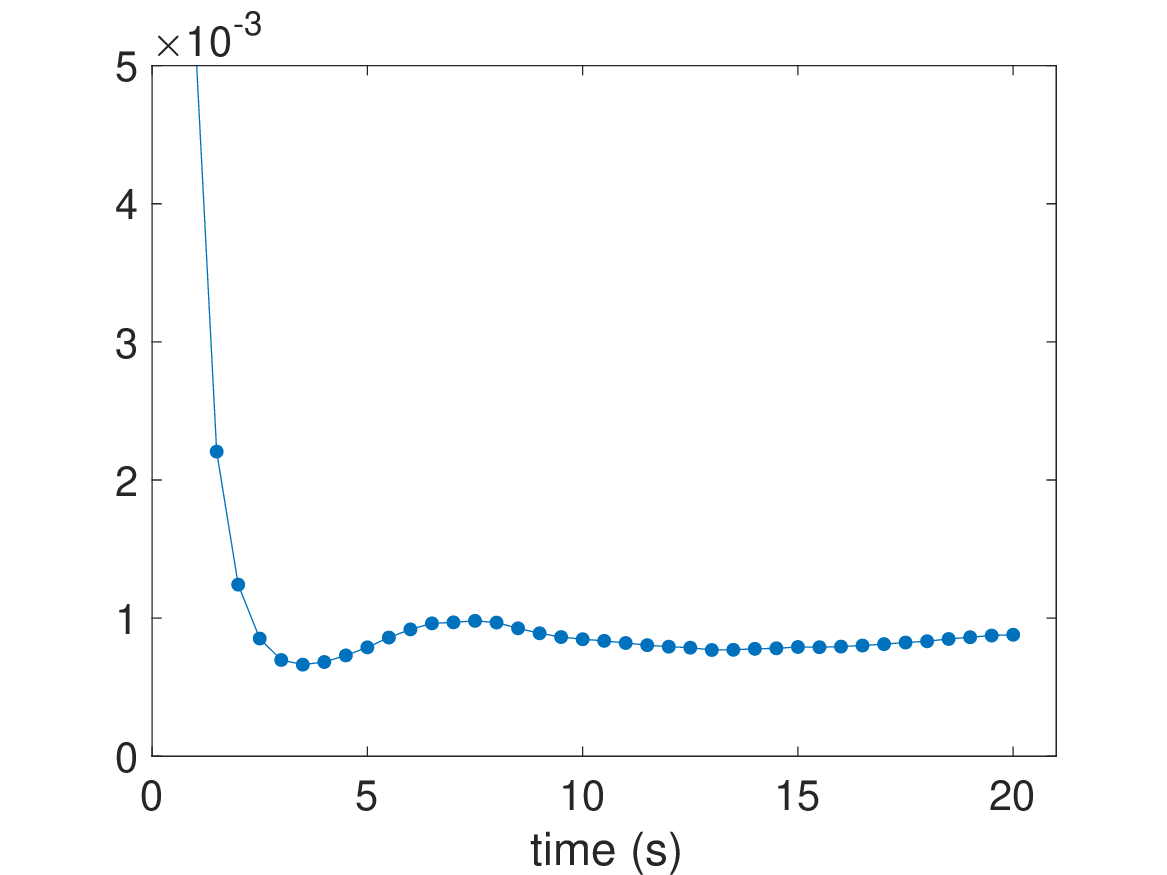}
\end{overpic}
\begin{overpic}[width=0.47\textwidth,grid=false]{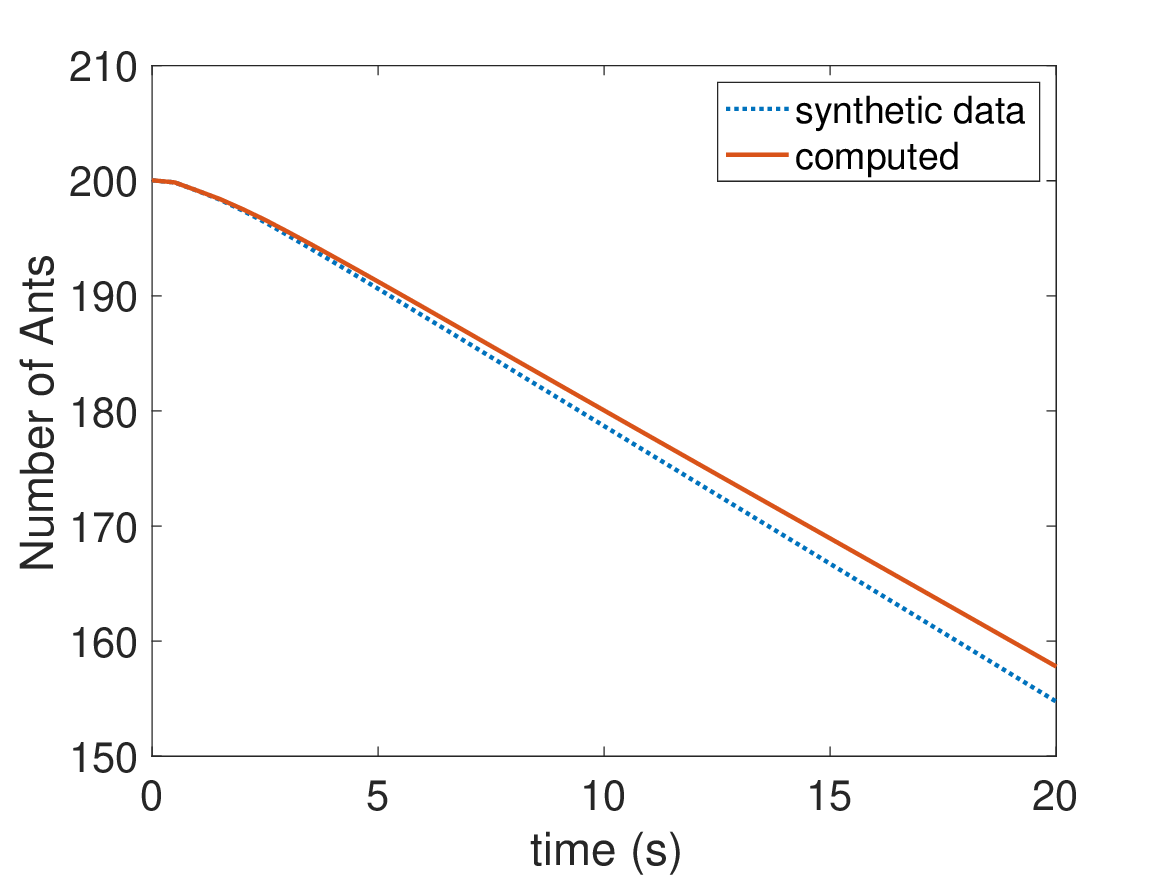}
\end{overpic}
\caption{Left: Functional \eqref{eq:JR_dt}
over time for the circular chamber with column.
Right: Number of ants inside the chamber over time. 
}
\label{fig:circle_col_J_ants}
\end{figure}

To quantify the agreement seen in Sec.~\ref{fig:circleprocesswcolumn}, 
Fig.~\ref{fig:circle_col_J_ants} (left) shows the time evolution of functional \eqref{eq:JR_dt}. 
Aside from the first time few time steps, 
the value of functional \eqref{eq:JR_dt}
is of the order of $10^{-3}$ most of the time.

Fig.~\ref{fig:circle_col_J_ants} (right) 
shows the number of ants inside the chamber over time. 
Although the computed densities in 
Fig.~\ref{fig:circleprocesswcolumn} agree with the synthetic data, we see that at $t = 20$ s there is a small difference 
in the number of ants left in the chamber according
to the synthetic data (158)
and the optimized results (155).

\subsection{Square chamber with exit at the corner of the walls}\label{sec:square}

\begin{figure}[htb!]
\center
\includegraphics[width=0.47\textwidth]{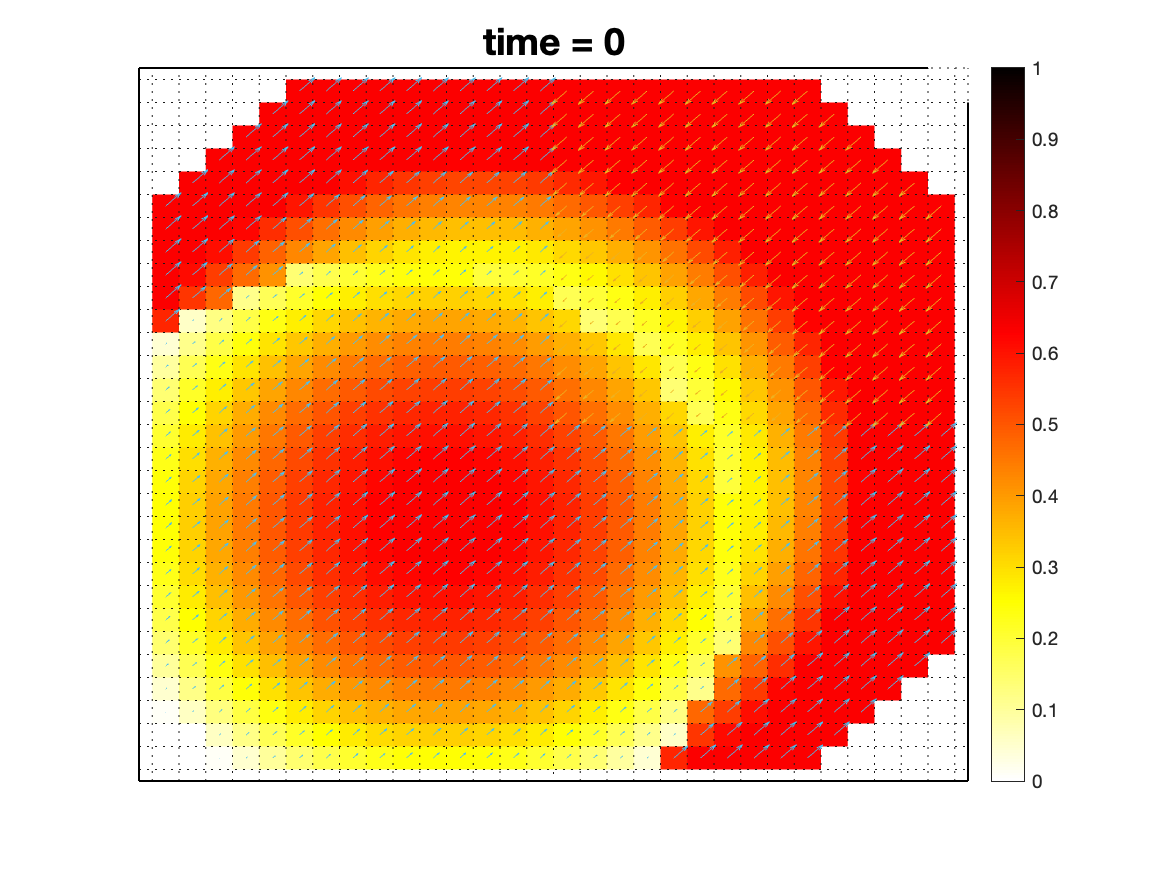}
\caption{Computational domain for the square chamber with exit at the corner of the walls (top right).}
\label{fig:initialsquare}
\end{figure}

The square chamber has a surface area equivalent 
to that of circular chamber, with the side measuring
31 mm. 
The exit is located at the upper right corner and it has size 2.5 mm. Like in the previous tests, we place 200 ants inside the chamber. We choose an initial distribution that mimics
the one used in Sec.~\ref{sec:circular}, i.e., 
a circular group and a crescent-shaped group.
The initial moving direction was assigned as follows:
we divided the square in quadrants and the ants in quadrant
near the exit have initial direction $\theta_6$, while
all others have initial direction $\theta_2$.
See Fig.~\ref{fig:initialsquare}.

\begin{figure}[htb!]
\centering
\begin{overpic}[width=0.32\textwidth, grid=false,tics=10]{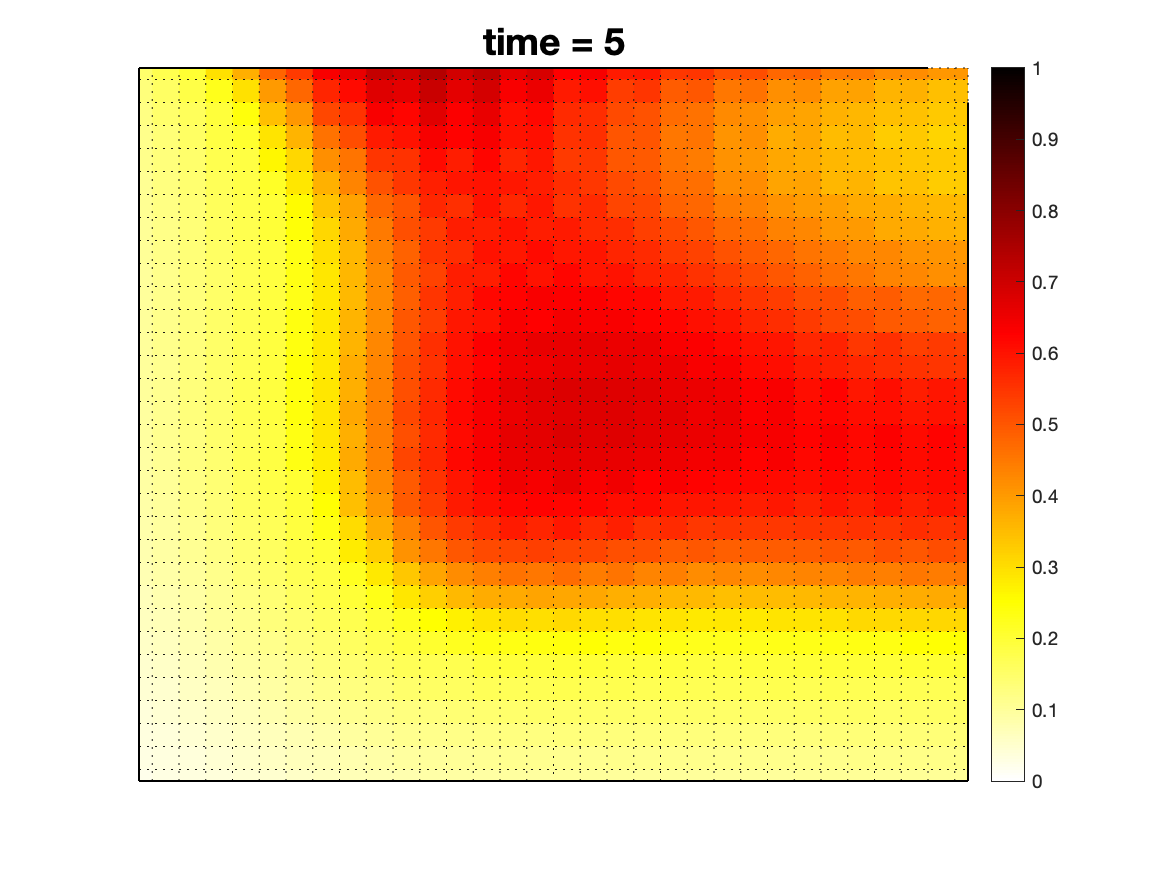}
\end{overpic} 
\begin{overpic}[width=0.32\textwidth, grid=false,tics=10]
{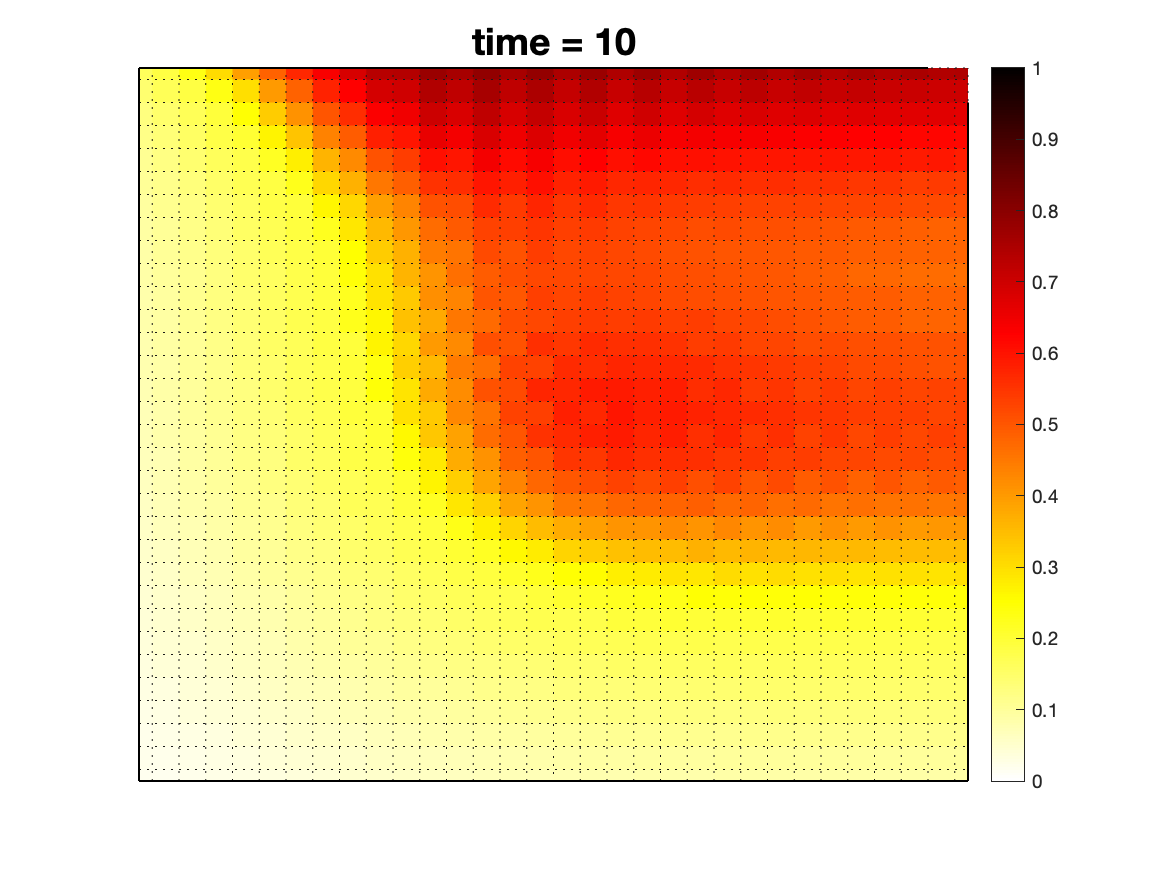}
\put(30,78){\textcolor{black}{$\varepsilon=0.05$}}
\end{overpic} 
\includegraphics[width=0.32\textwidth]{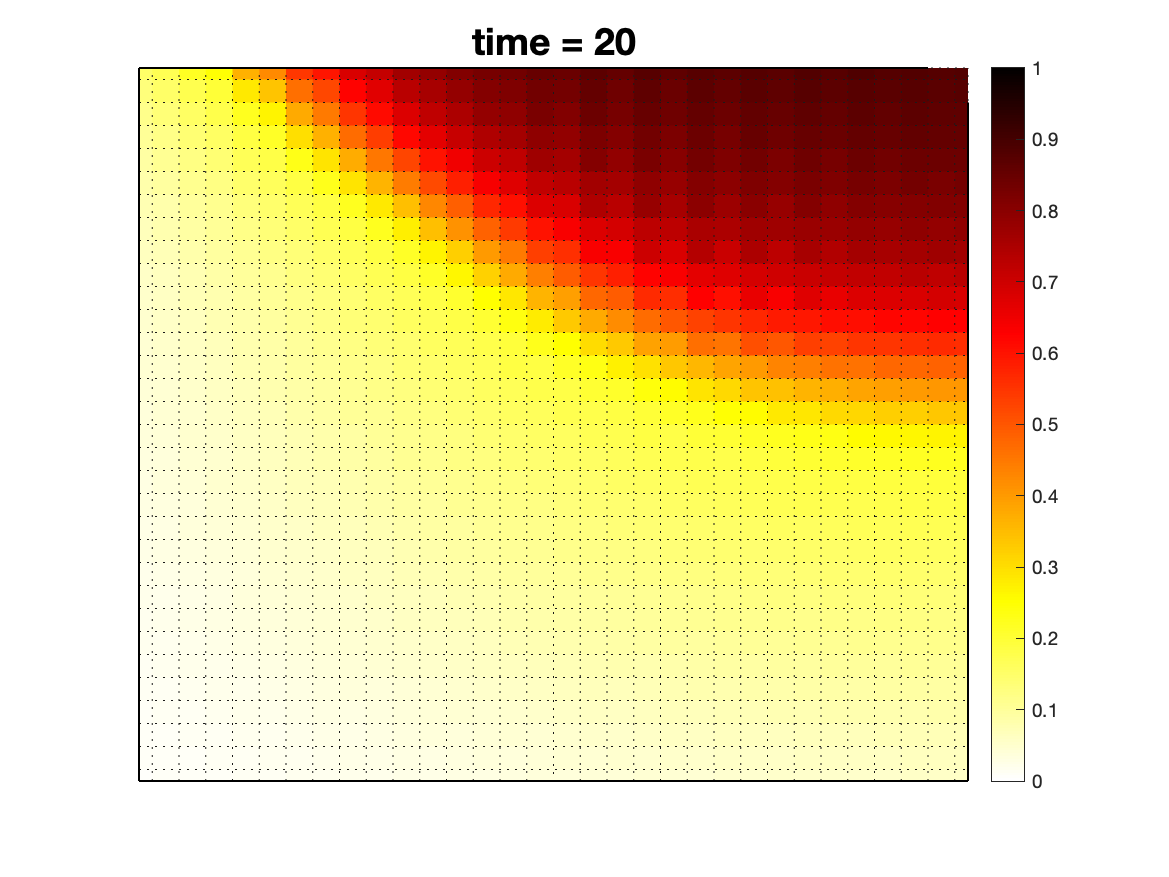}\\
\vskip .2cm
\begin{overpic}[width=0.32\textwidth, grid=false,tics=10]{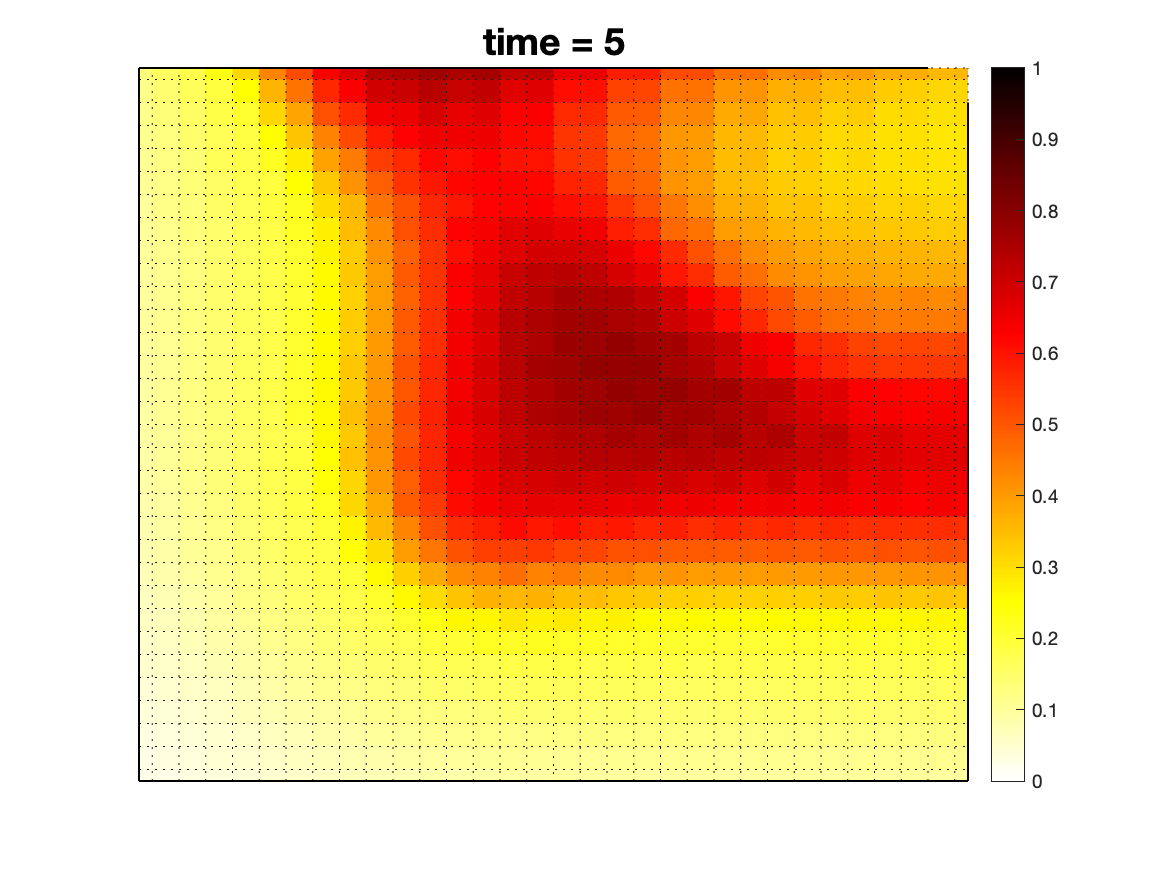}
\end{overpic} 
\begin{overpic}[width=0.32\textwidth, grid=false,tics=10]{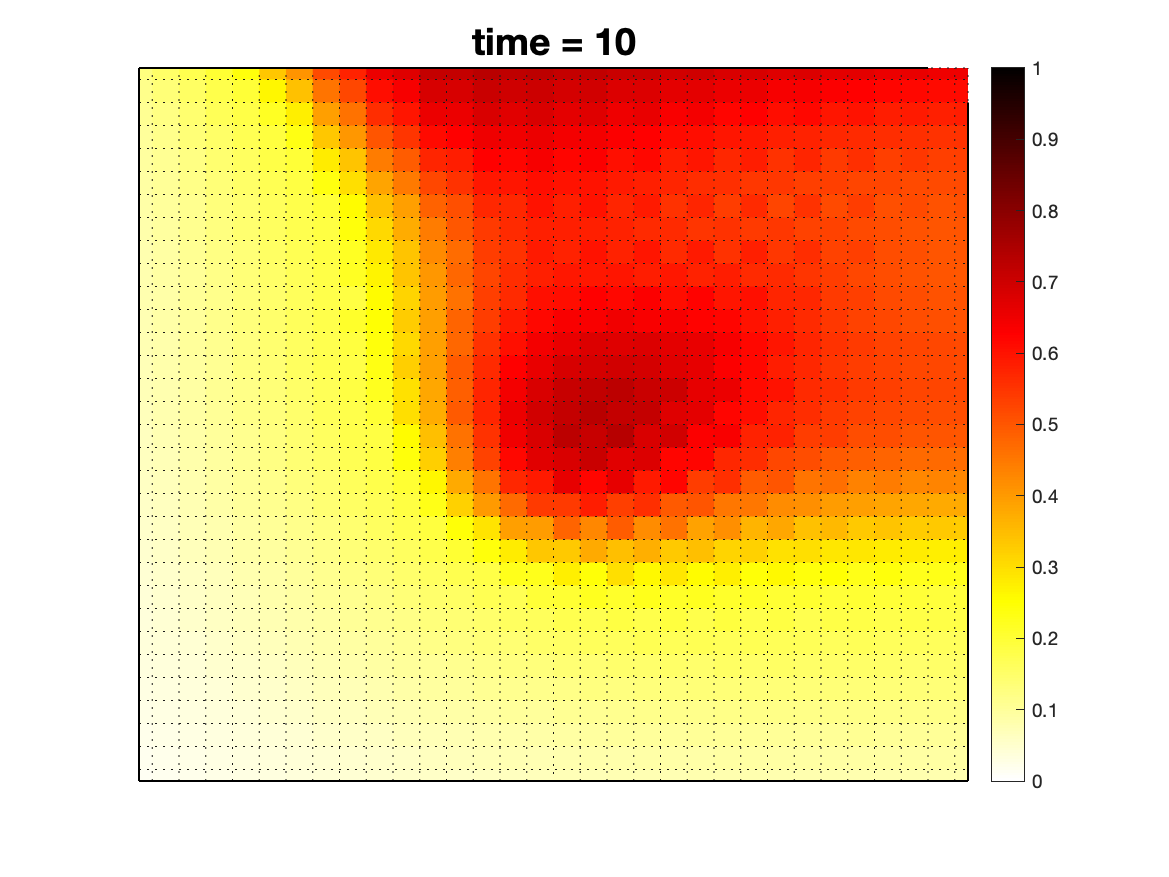}
\put(30,78){\textcolor{black}{$\varepsilon=0.95$}}
\end{overpic} 
\includegraphics[width=0.32\textwidth]{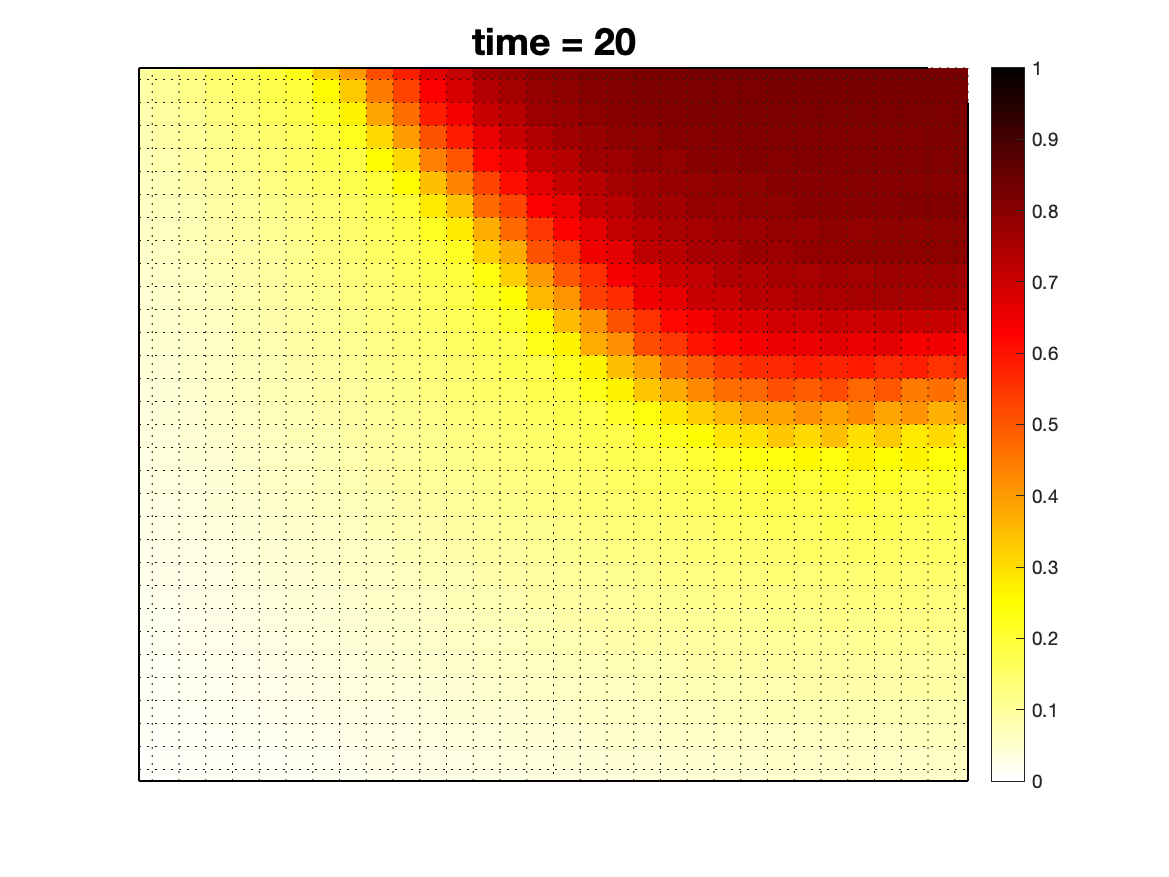}
\caption{Computed density at $t = 5$ s (left), $t = 10$ s (center), and $t = 20$ s (right)
given by the forward problem 
for 200 ants initially placed as
in Fig.~\ref{fig:initialsquare}
for $\varepsilon=0.05$ (top row) and $\varepsilon=0.95$ (second row).}
\label{fig:squareprocess0}
\end{figure}

The evacuation dynamics in 
conditions of low stress ($\varepsilon = 0.05$)
and high stress ($\varepsilon = 0.95$)
are shown in Fig.~\ref{fig:squareprocess0}.
Like in the case of the circular chamber
with unobstructed exit in Sec.~\ref{sec:circular}, 
we see that 
regions of high density form around the center
of the chamber at $t = 5, 10$ s when the stress level is higher.
At $t = 20$ s,  the
high-density region becomes smaller and more uniform in space when the stress
level is $\varepsilon = 0.95$, which is similar to what observed in the circular chamber
with partially obstructed exit in Sec.~\ref{sec:circular_col}.

\begin{figure}[htb!]
\begin{overpic}[width=0.32\textwidth, grid=false,tics=10]{Square_Room_Ped_Eps0_95_Ants_square_time_10.eps}
\end{overpic} 
\begin{overpic}[width=0.32\textwidth, grid=false,tics=10]{Square_Room_Ped_Eps0_95_Ants_square_time_20.eps}
\put(-10,78){\textcolor{black}{synthetic video data ($\varepsilon=0.95$)}}
\end{overpic} 
\includegraphics[width=0.32\textwidth]{Square_Room_Ped_Eps0_95_Ants_square_time_40.eps} \\
\vskip .2cm
\begin{overpic}[width=0.32\textwidth]{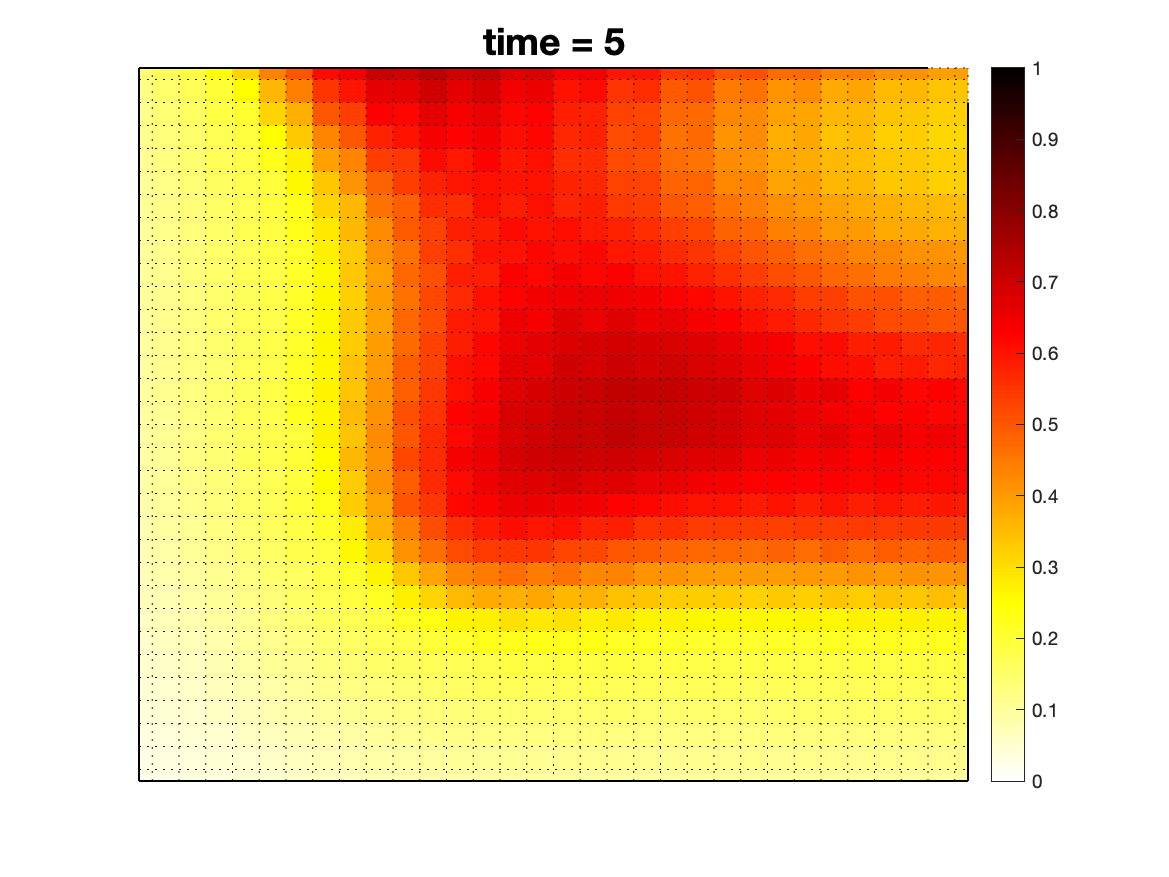}
\end{overpic} 
\begin{overpic}[width=0.32\textwidth]{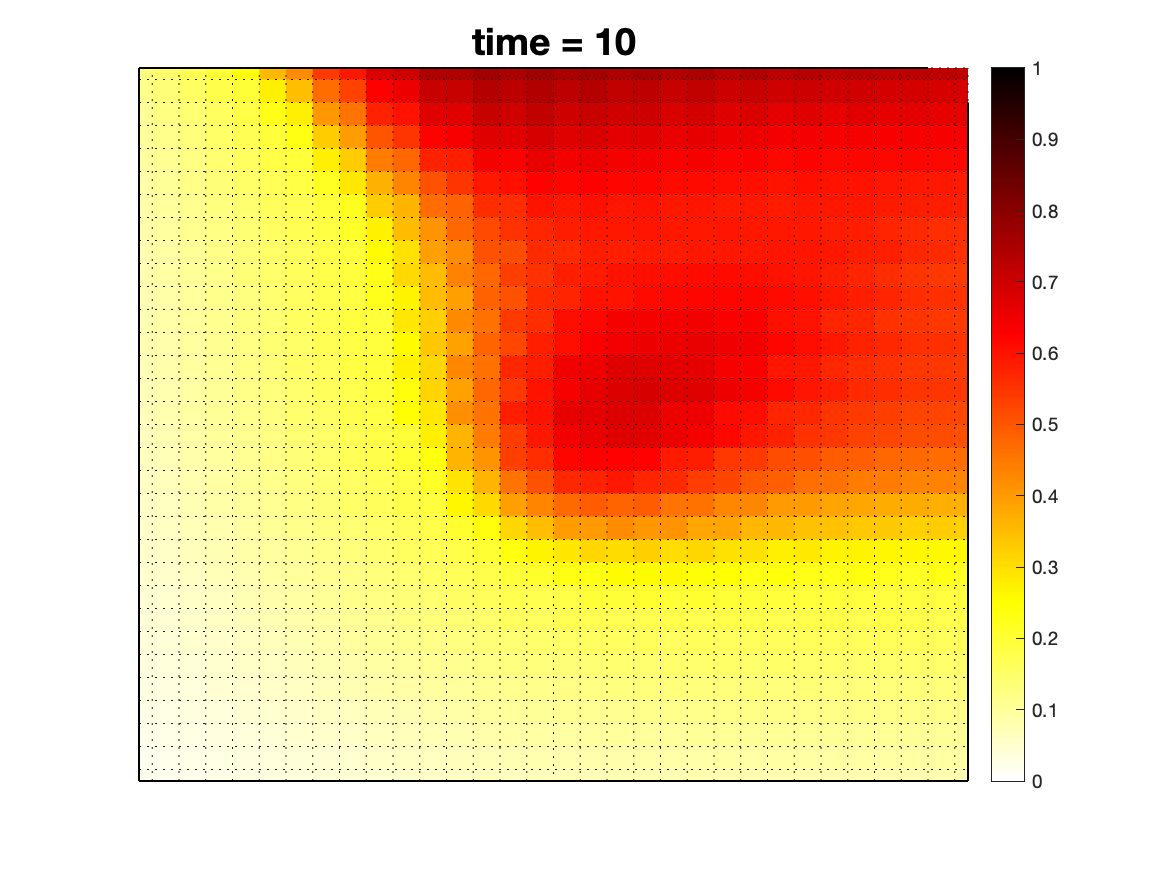}
\put(-15,78){\textcolor{black}{optimized density (from $\varepsilon=0.05$)}}
\end{overpic}
\includegraphics[width=0.32\textwidth]{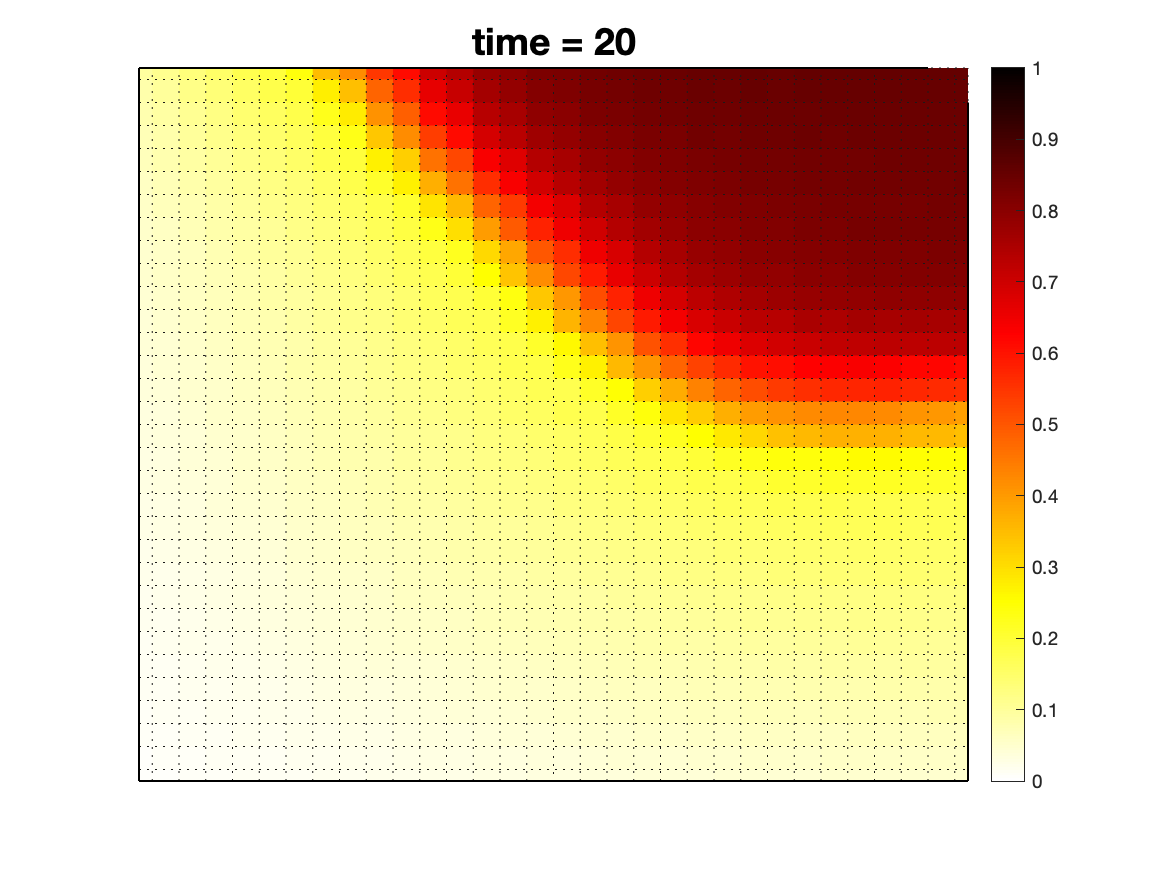}
\\
\vskip .2cm
\begin{overpic}[width=0.32\textwidth]{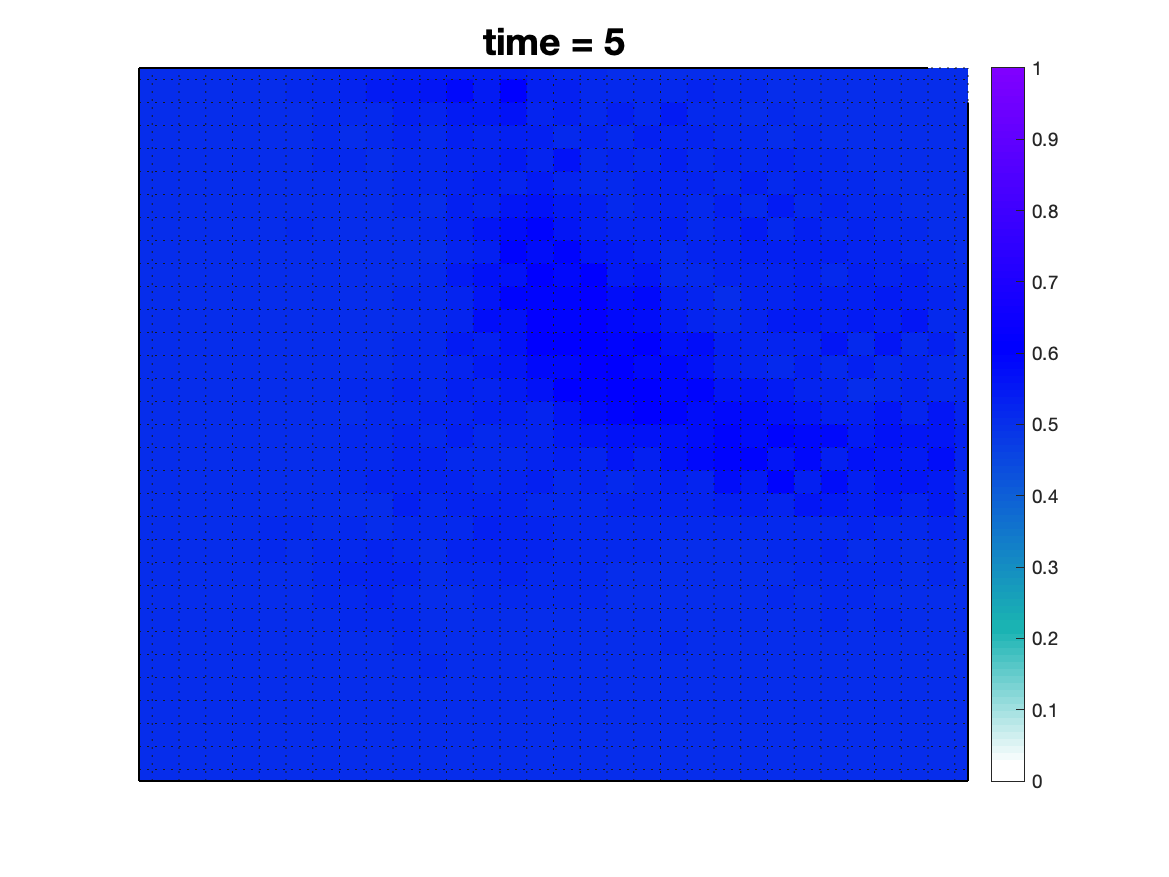}
\end{overpic} 
\begin{overpic}[width=0.32\textwidth]{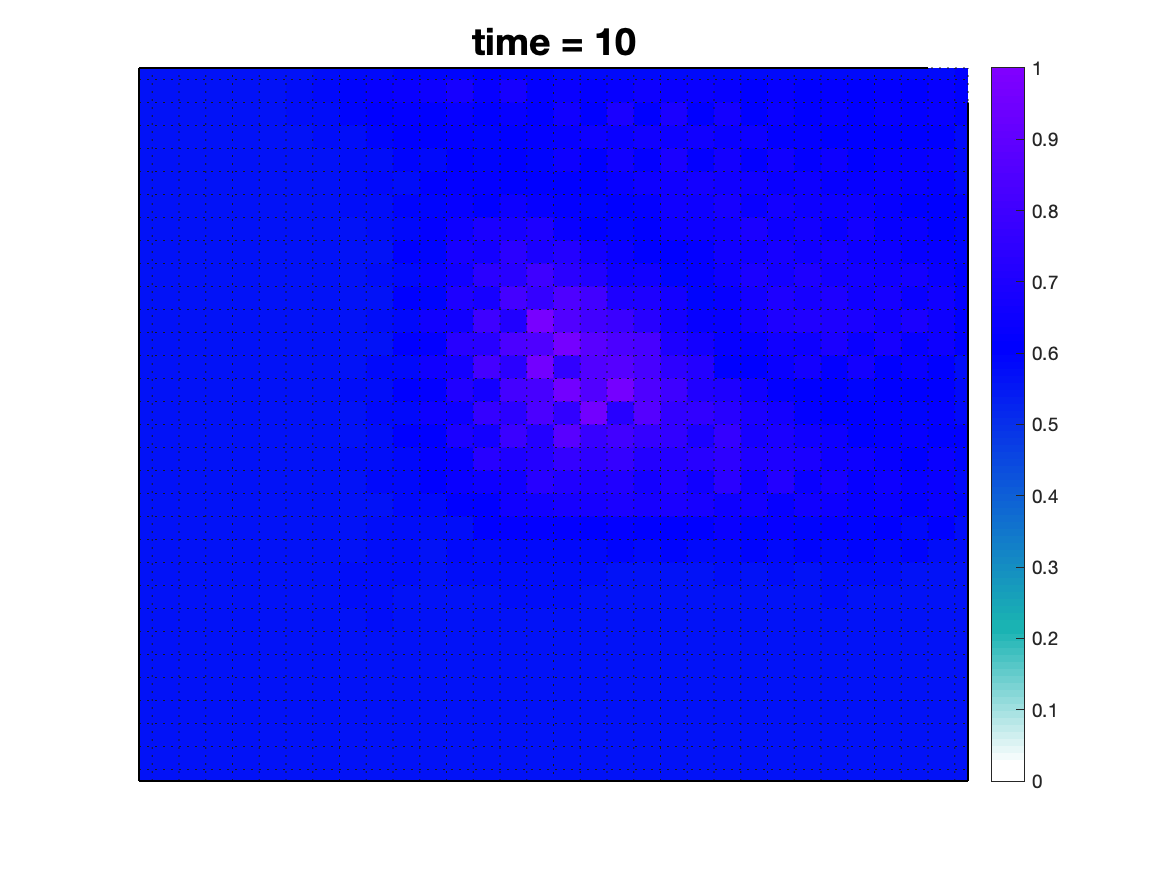}
\put(-20,78){\textcolor{black}{optimized stress level (from $\varepsilon=0.05$)}}
\end{overpic} 
\includegraphics[width=0.32\textwidth]{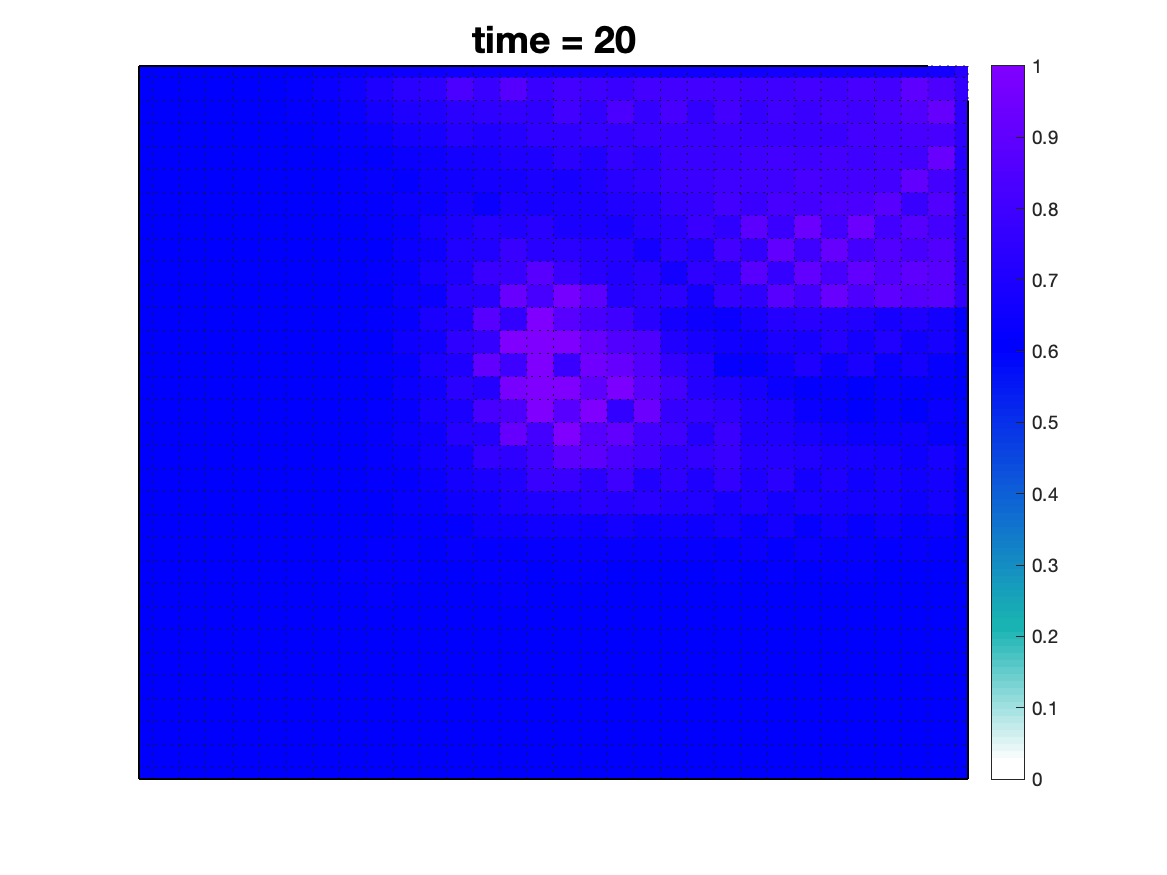}
\caption{Synthetic density data (top), 
optimized density (center), optimized
stress level (bottom) at $t = 5$ s (left), $t = 10$ s (center), and $t = 20$ s (right). Ants are initially placed as
in Fig.~\ref{fig:initialsquare}.
The optimized results were obtained using regularization \eqref{eq:td_J_R} with
$\varepsilon_{ref} = 0.75$.
}
\label{fig:squareprocess}
\end{figure}

We use same mesh size and time step as in the previous sections, and treat as synthetic video data
the computed results for $\varepsilon=0.95$. The 
optimization procedure is started from $\varepsilon=0.05$, with $\varepsilon_{ref}=0.75$ and the same $\xi$, $\delta$, and $tol$
as in Sec.~\ref{sec:circular_col}. 
The comparison of synthetic density data and optimized
density is shown in Fig.~\ref{fig:squareprocess}, together with the corresponding optimized stress level.
Thanks to the use of relaxation, 
the synthetic data and the optimized
density match well. We see that the high-density regions around the center of the 
chamber at $t = 5,10$ s are well captured.
Similarly, the the high-density regions in the synthetic data and optimized density at 
$t = 20$ match well.
The optimized stress level remains around the $\varepsilon_{ref}$ value over most of the domain, 
but it has local peaks close to $0.95$.

\begin{figure}[h]
\centering
\begin{overpic}[width=0.47\textwidth,grid=false]{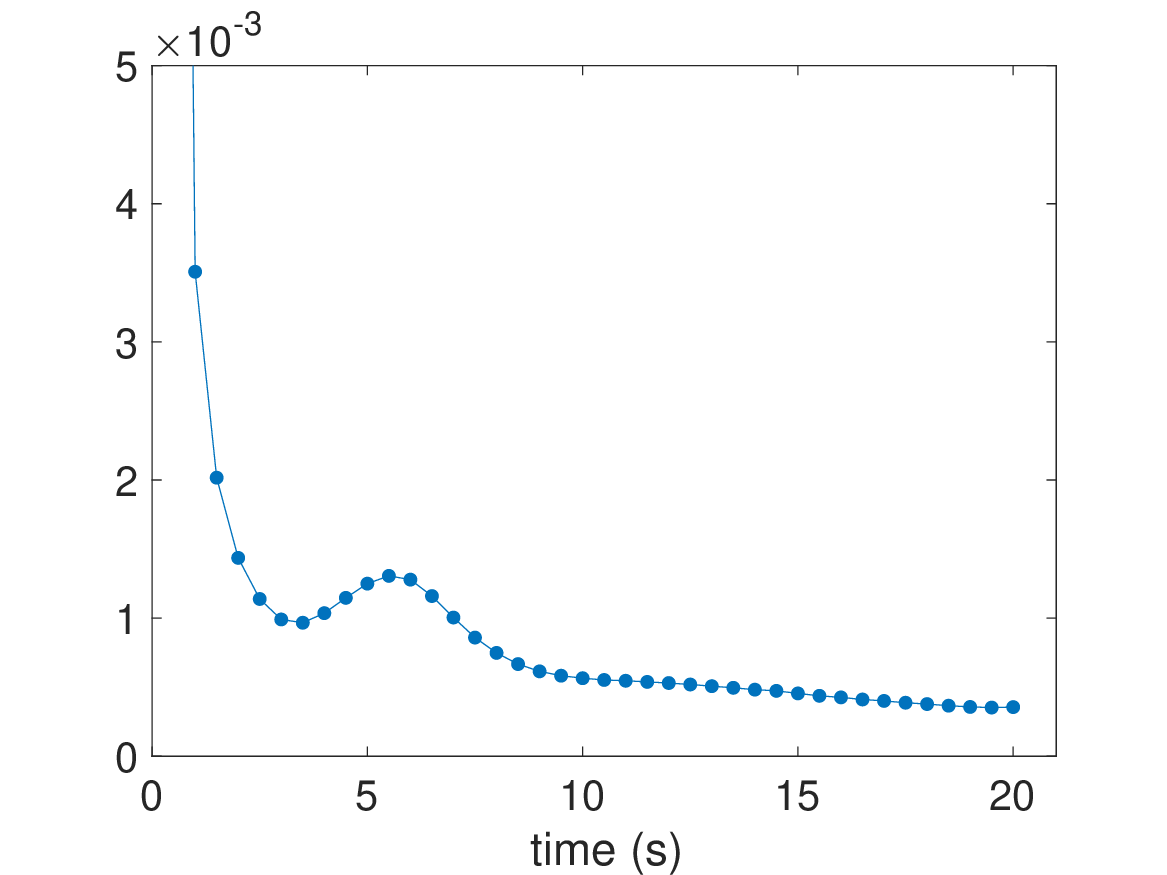}
\end{overpic}
\begin{overpic}[width=0.47\textwidth,grid=false]{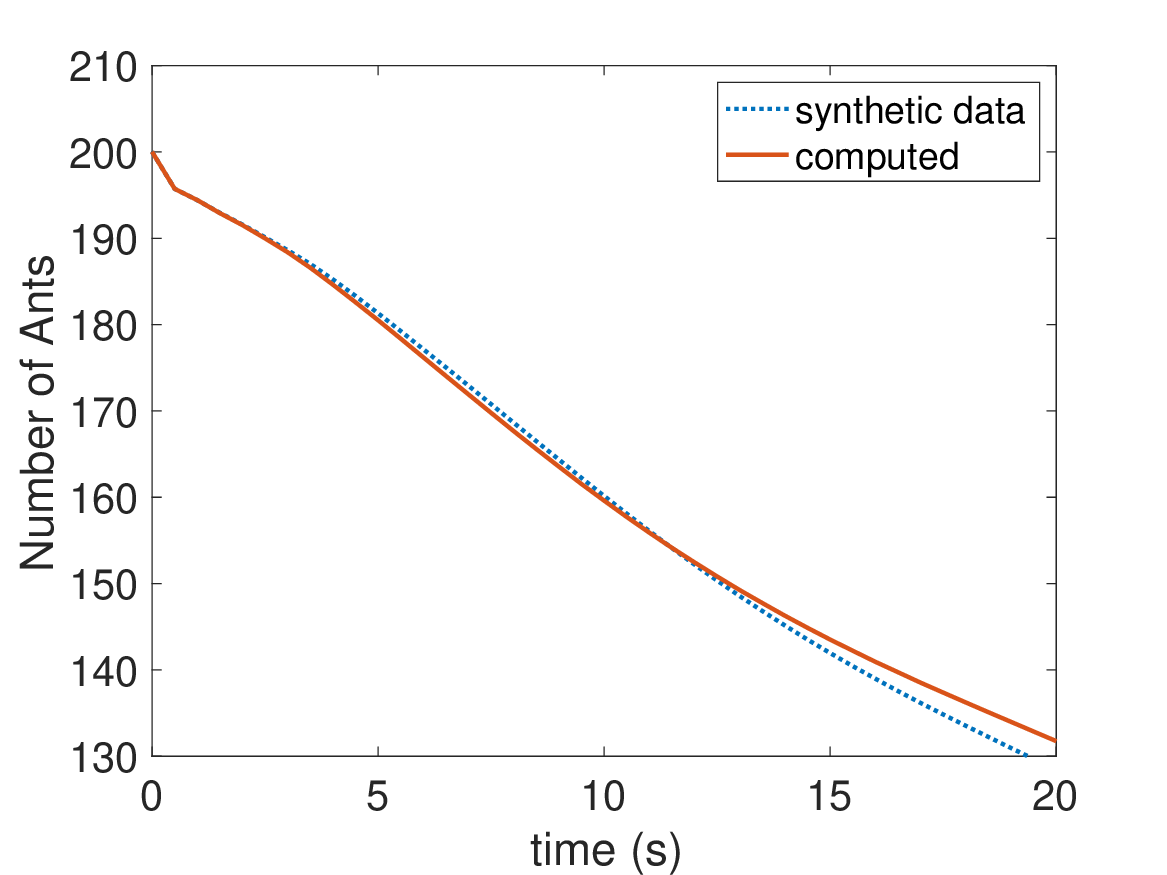}
\end{overpic}
\caption{Left: Functional \eqref{eq:JR_dt}
over time for the square chamber.
Right: Number of ants inside the chamber over time.
}
\label{fig:square_J_ants}
\end{figure}

Fig.~\ref{fig:square_J_ants} (left) shows the time evolution of functional \eqref{eq:JR_dt}. Again, 
we see that the value of functional \eqref{eq:JR_dt}
is of the order of $10^{-3}$ most of the time.

Fig.~\ref{fig:square_J_ants} (right) 
shows the number of ants inside the chamber over time. The numbers of ants given by the synthetic
data and the optimized simulation agree till about 
$t = 14$ s and then they grow slightly apart. 
From Fig.~\ref{fig:square_J_ants} (right), we see that it takes about 12.5 s for the first 50 ants
to leave the room. This is also consistent with
the measured mean escape time (out of 10 repetitions in 
Ref.~\citenum{SHIWAKOTI20111433}) for the first 50 ants, which is 12.9 s $\pm$ 0.5 s standard deviation.

\section{Conclusions and perspectives}\label{sec:concl}
The large majority of existing mathematical models for crowd dynamics assume that pedestrians do not alter their, typically
rational, walking behavior during the time interval of interest. Such models
cannot be used to simulate emergency situations, like forced evacuations, when people may behave irrationally
in response to fear. In this paper, we consider a kinetic model for crowd dynamics that features a parameter representing the stress level
and we present a data-driven technique to set this parameter. 

The data-driven technique is a minimization approach, based on an inverse crowd dynamics problem
and a suitable functional measuring the distance between given
crowd density data and
the numerical solution. We applied this new methodology to model problems inspired from the experiments on panicked ants in Ref.~\citenum{SHIWAKOTI20111433}, with
synthetic data generated by numerical simulations of the forward problem.
Our numerical experiments indicate that regularization is needed to obtain computed results that match well the density data. We used a regularization term that requires a reference
value for the stress level parameter, i.e., one needs a good guess for
the stress level on order for the proposed approach to yield accurate results.

The next step is to apply our methodology with the data extracted from the video recordings of the experiments in Ref.~\citenum{SHIWAKOTI20111433}
and further investigate the sensitivity to the reference
stress level.



\section*{Acknowledgments}
We thank Dr.~Martin Burd (Monash University) and 
Dr.~Nirajan Shiwakoti (RMIT University) for sharing with us the video
recordings of their experiments on ants in the circular chamber with and without partial obstruction near the exit. D.L. acknowledges support of Simons Foundation grant MPS-TSM-00002738.



\bibliographystyle{plain}
\bibliography{contagion_model_update}

\end{document}